\documentclass[11pt, reqno]{amsart}
\usepackage{mathrsfs, mathpazo}
\usepackage{amsfonts, amsbsy}
\usepackage{ dsfont}
\usepackage{amssymb, mathabx, mathtools}
\usepackage{xcolor}
\usepackage{cancel}

\newcommand{\OLD}[1]{\ifmmode{\color{red}\cancel{#1}}\else{\color{red}#1}\fi}

\usepackage[colorlinks,citecolor=blue,urlcolor=blue]{hyperref}
\usepackage[english]{babel}
\usepackage{amsmath}
\theoremstyle{remark}
\newcommand{\re}[1]{(\ref{#1})}

\newcommand{\Sup} {\displaystyle \sup}

\newcommand{\Int }    {\displaystyle \int}

\newcommand{\abs}[1]{\left\vert#1\right\vert}
\newcommand{\norm}[1]{\left\Vert#1\right\Vert}

\selectfont

\theoremstyle{plain}
\newtheorem{theorem}{\textbf{Theorem}}[section]
\newtheorem{lemma}[theorem]{\textbf{Lemma}}
\newtheorem{claim}[theorem]{\textbf{Claim}}

\newtheorem{proposition}[theorem]{\textbf{Proposition}}
\newtheorem{cor}[theorem]{\textbf{Corollary}}
\theoremstyle{definition}
\newtheorem{assumption}[theorem]{\textbf{Assumption}}

\newtheorem{remark}[theorem]{\textbf{Remark}}

\newtheorem{definition}[theorem]{\textbf{Definition}}

\numberwithin{equation}{section}
\numberwithin{figure}{section}
\newcommand{\bz}{\mathcal{Z}}
\newcommand{\ba}{\mathcal{A}}
\newcommand{\bx}{\mathcal{X}}
\newcommand{\by}{\mathcal{Y}}

\newcommand{\be}{\mathbb{E}}

\newcommand{\bj}{\mathcal{J}}

\newcommand{\bo}{\mathcal{O}}
\newcommand{\bk}{\mathcal{K}}
\newcommand{\br}{\mathbb{R}}

\newcommand{\eps}{\varepsilon}

\title[Well-posedness of 2D stochastic  Keller-Segel system with logistic source]{Well-posedness of 2D stochastic parabolic-parabolic Keller-Segel system with logistic source}
\author{Erika Hausenblas$^{1}$,  Boris Jidjou Moghomye$^{1}$, Paul Andr\'e Razafimandimby$^2$}
\dedicatory{\vspace{-10pt}\normalsize{$^{1}$  Department of Mathematics and Information Technology, Technical University of  Leoben, Leoben Franz Josef Strasse 18, 8700 Leoben, Austria \\
$^{2}$ School of Mathematical Sciences, Dublin City University, Collins Avenue, Dublin 9, Ireland		
}}

\begin{document}

	\maketitle

\begin{abstract}
We consider a stochastic, fully parabolic Keller--Segel system with multiplicative noise in non-divergence form and a logistic source in a two-dimensional bounded domain, subject to random initial data. We construct a global-in-time, non-negative martingale solution by means of a suitable regularisation procedure, combined with Banach fixed-point arguments, a localisation procedure, and a stochastic compactness method. We further establish pathwise uniqueness and uniqueness in law for martingale solutions, and prove the existence of a unique global-in-time, non-negative weak solution without imposing any smallness condition on the initial data. 
\end{abstract}

\section{Introduction}
\label{sec:intro}


The Keller--Segel system, introduced in \cite{Keller,Patlak}, is one of the
cornerstones of mathematical biology. It describes the chemotactic aggregation
of cells, for instance, the amoeba \emph{Dictyostelium discoideum}, in
response to a chemical signal that the cells themselves secrete. In two space dimensions, the classical deterministic system exhibits a
striking mass-dependent dichotomy driven by the competition between
diffusion and aggregation: in the absence of logistic terms, sufficiently
small initial mass ensures global existence, whereas sufficiently large
initial mass may lead to finite-time blow-up; see, for instance,
\cite{Herrero,Horstmann1,Horstmann,Winkler1}.
The addition of a logistic source with quadratic damping will change fundamentally
 this picture: in the deterministic setting, it suppresses
finite-time blow-up and ensures global existence even for arbitrarily
large initial mass; see, e.g., \cite{Lankeit,Tello,Winkler,Zheng}.
This damping mechanism models the saturation of the cell population
due to limited resources. A central question is whether its stabilizing
effect persists under random perturbations. Such perturbations arise
naturally in biological and physical environments and can influence
both cell-population dynamics and the evolution of the chemical signal,
motivating the study of stochastic Keller--Segel systems. An early
stochastic formulation was proposed by Chavanis \cite{Chavanis},
building on the approach introduced by Dean \cite{Dean}.
To examine this question, we introduce stochastic forcing into both
the cell-density and chemical-concentration equations. To be more precise,
we consider the following stochastic Keller--Segel system, inspired
by the formulation of Chavanis \cite{Chavanis} we get
\begin{equation}\label{equa1}
	\begin{split}
	&	dn(t,x)-\Delta n(t,x)\,dt+\chi\nabla\cdot(n(t,x)\nabla c(t,x))\,dt
		=g_1(n(t,x))\,dW_1(t), \ t\in [0,T], \ x\in \bo\\[2mm]
	&\kappa	dc(t,x)-\Delta c(t,x)\,dt=n(t,x)\,dt-c(t,x)\,dt+\kappa g_2(c(t,x))\,dW_2(t),  \ t\in [0,T], \ x\in \bo,
	\end{split}
\end{equation}
where $T>0$ is a fixed final observation time and $\mathcal{O}\subset\mathbb{R}^{d}$, with $d\in{\{}1,2,3{\}}$,  is a smooth, bounded domain with boundary $\partial\mathcal{O}$, in which the chemotactic process takes place. The unknowns
$n$ and $c$ represent the density of \emph{Dictyostelium discoideum} and the
concentration of the chemical attractant, respectively.  The constant
$\chi>0$ is the chemotactic sensitivity, while  the functions $g_1$ and $g_2$  characterize the  random perturbations, and $W_i$, with $i\in \{1,2\}$ are two independent cylindrical Wiener processes defined  on a given filtered probability space $(\Omega,\mathcal{F},\mathbb{P},\mathbb{F}_T)$ with filtration
$\mathbb{F}_T:=\{\mathcal{F}_t\}_{t\in[0,T]}$ satisfying the usual conditions.  The parameter $\kappa\in\{0,1\}$ distinguishes between the two diffusion regimes for the chemical concentration. The parameter $\kappa$ determines the structure of the stochastic
Keller--Segel system \eqref{equa1}. For $\kappa=0$, the chemical
concentration is governed by an elliptic equation, corresponding to
a quasi-steady approximation in which the chemical equilibrates much
faster than the cell density evolves. This yields the stochastic
parabolic--elliptic system. For $\kappa=1$, the chemical concentration
retains its time-dependent dynamics, giving rise to the stochastic
parabolic--parabolic system.
The mathematical analysis of stochastic Keller--Segel systems has
attracted considerable attention, with well-posedness and blow-up
among the central questions. In the parabolic--elliptic setting,
Chevyrev, Hambly and Mayorcas \cite{Chevyrev} established global
well-posedness and the existence of a unique invariant measure
for a model with chemorepulsion on the one-dimensional torus.
Misiats, Stanzhytskyi
and Topaloglu \cite{Misiats} showed that in two dimensions, divergence-form
perturbations admit global weak solutions for sufficiently small initial data,
whereas non-divergence-form noise leads to finite-time blow-up with positive
probability. Mayorcas and Toma\v{s}evi\'c \cite{Mayorcas} provided blow-up
criteria involving the chemotactic sensitivity.  For the fully parabolic case, Hausenblas, Mukherjee, and Tran \cite{Hausenblas} established the existence of a martingale solution in one spatial dimension. More recently, Choi et al. \cite{Choi} proved the existence of a local solution in three spatial dimensions. The corresponding two-dimensional case appears to be less understood, despite the particular importance of this dimension for the Keller--Segel system. Indeed, in two dimensions, the competition between diffusion and chemotactic aggregation leads to the well-known critical behavior of the deterministic system. This observation naturally raises the question of whether the stabilizing effect of a logistic source persists under stochastic perturbations. In this direction, Chen, Zhai, and Zhang \cite{Chen} considered the stochastic parabolic--elliptic Keller--Segel system with a logistic source and nonlinear noise, and established the existence of mild solutions. In contrast, we consider the following two-dimensional stochastic parabolic--parabolic Keller--Segel system with a logistic source, driven by multiplicative noise in non-divergence form:
\begin{equation}\label{equa1.1}
	\begin{split}
	&	dn(t,x)-\Delta n(t,x)\,dt+\chi\nabla\cdot(n(t,x)\nabla c(t,x))\,dt\\
&		=\alpha n(t,x)\,dt-\mu n^{2}(t,x)\,dt+g_1(n(t,x))\,dW_1(t), \ t\in [0,T], \ x\in \bo,\\[2mm]
	&	dc(t,x)-\Delta c(t,x)\,dt=n(t,x)\,dt-c(t,x)\,dt+g_2(c(t,x))\,dW_2(t),  \ t\in [0,T], \ x\in \bo.
	\end{split}
\end{equation}
Here, $\alpha\in\mathbb{R}$ represents the difference between the birth and death rates, while $\mu{>}0$ measures the strength of the density-dependent damping. Accordingly, the terms $\alpha n$ and $-\mu n^2$ describe the logistic growth and damping of the cell population, respectively. The stochastic perturbations are driven by the independent cylindrical Wiener processes  defined  for any $i\in\{1,2\}$ by
\begin{equation}\label{equ1.3}
W_{i}(t):=\sum_{j=1}^{\infty}
\beta_j^{i}(t)e_j,
\qquad t\in[0,T],
\end{equation}
where $(\beta_j^{i}(t))_{t\in[0,T]}$, $j\in\mathbb{N}$ and $i\in\{1,2\}$, is a family of independent standard Brownian motions defined on the filtered probability space $(\Omega,\mathcal{F},\mathbb{P},\mathbb{F}_T)$. Moreover, $\{e_j\}_{j=1}^{\infty}$ denotes an orthonormal basis of $L^{2}(\mathcal{O})$ consisting of eigenfunctions of the Laplacian subject to homogeneous Neumann boundary conditions.

To complete the formulation of the system, we supplement the above equations with suitable boundary and initial conditions, given by 
\begin{equation}\label{equa1.2}
	\begin{cases}
		\mathbf{n}\cdot\nabla n(t,x)=\mathbf{n}\cdot\nabla c(t,x)=0
		& \text{on } [0,T]\times\partial\mathcal{O},\\[1mm]
		n(0,x)=n_0(x),\quad c(0,x)=c_0(x)
		& \text{on }  \mathcal{O},
	\end{cases}
\end{equation}
where $\mathbf{n}$ is the outward normal vector to $\partial\mathcal{O}$.  

Our main goal is to establish the existence and uniqueness of global weak solutions to the system \eqref{equa1.1}--\eqref{equa1.2}, showing that, in two spatial dimensions and under multiplicative noise in non-divergence form, the quadratic logistic damping is sufficiently strong to prevent the finite-time blow-up observed in the absence of logistic source \cite{Mayorcas,Misiats}, without imposing any smallness condition on the initial data.
Our contributions are threefold.

$\bullet$ \textbf{Global existence of martingale solutions in the fully parabolic two-dimensional case without smallness conditions.}
We prove that, for any Borel probability measures $\vartheta_1$ and $\vartheta_2$ on suitable metric spaces, respectively, with finite moments of order $32$, the system \eqref{equa1.1}--\eqref{equa1.2} admits at least one global-in-time martingale solution, which is non-negative almost surely. The result holds for arbitrary initial data and for all parameters $\chi>0$, $\mu>0$, and $\alpha\in\mathbb R$, without imposing any smallness condition. This appears to be the first global	existence result for the two-dimensional fully parabolic stochastic
Keller--Segel system with logistic source and non-divergence-form noise. We construct a martingale solution in four steps. First, we introduce
	a regularised system, parameterised by $\epsilon\in(0,1)$, in which the
	cross-diffusion term and the logistic term are truncated in a way that
	preserves non-negativity. The well-posedness of the regularised system is
	established via the Banach fixed-point theorem in suitable Banach spaces of
	path-continuous   functions. A localisation procedure then
	yields a maximal local solution, and uniform a priori estimates upgrade it to
	a global one.
	Second, we prove that the approximate solutions are non-negative. This is
	achieved by approximating the negative-part functional by a family of $C^2$
	functions and applying the It\^o formula; the specific structure of the noise
	coefficient $g_2$ (see Assumption~\ref{assum2.2}(A$_3$) below for the precise formulation) is used here to control
	the stochastic contribution.
	Third, we derive uniform-in-$\epsilon$ bounds. Because of the strong
	nonlinearity of the cross-diffusion term, the regularity obtained is
	insufficient to apply the Markov inequality and the Aldous condition directly
	to the density. We therefore develop a tightness criterion based on
	\emph{uniform bounds in probability} rather than uniform bounds of norms. This
	is the technical heart of the paper.
	Fourth, we apply the Jakubowski--Skorokhod representation theorem
	\cite{Jakubowski} to pass to the limit along a subsequence and identify the
	limit as a non-negative martingale solution. 
	criterion is what allows us to preserve the non-negativity of the
	approximation in the limit. 
	
$\bullet$  \textbf{A compactness argument that replaces uniform moment bounds by
		uniform bounds in probability.} The obstruction described above rules out the
	classical route to compactness, which relies on uniform-in-$\epsilon$ moment
	estimates in high norms for the density. We circumvent it by establishing
	\emph{uniform bounds in probability} for the approximate solutions, which are
	sufficient to apply the Jakubowski--Skorokhod representation theorem
	\cite{Jakubowski} and to pass to the limit. 
	We expect
	the technique to be applicable to other stochastic nonlinear systems 
	for which global moment estimates fail.
	
$\bullet$  \textbf{Pathwise uniqueness, uniqueness in law, and a unique weak
		solution.} We establish pathwise uniqueness, uniqueness in law, and
	the existence of a unique global non-negative weak solution with prescribed
	initial data. 
Pathwise uniqueness is obtained via a Gronwall inequality applied
to the difference of two solutions, with the cross-diffusion term controlled
by the Gagliardo--Nirenberg interpolation inequality and the Lipschitz
properties of the noise coefficients. Uniqueness in law then follows from
\cite[Theorem~5.3]{Kunze}, and the existence of a unique weak solution from
\cite[Corollary~5.4]{Kunze}.


The rest of the paper is organised as follows. Section \ref{sec:prelim} introduces the notation, standing assumptions, the notions of martingale and weak solutions, and the precise statements of our main results. Section \ref{sec:proof1} is devoted to the proof of the first main result, while the proof of the second main result is presented in Section \ref{sec:proof2}. Appendix~\ref{app:reg} contains the construction of the solutions of the regularised system (proof of Proposition~\ref{propo3.5}).
	
\section{Notations, preliminaries, standing assumptions and main results}
\label{sec:prelim}

\subsection{Preliminaries and standing assumptions}
\label{subsec:preliminaries}

For $a,b\in\mathbb{R}$, we set $a\wedge b:=\min(a,b)$ and
$a\vee b:=\max(a,b)$. For $m\in\mathbb{N}$ and $1\le p\le\infty$,
$L^p(\mathcal{O})$ and $W^{m,p}(\mathcal{O})$ denote the standard Lebesgue
and Sobolev spaces on $\mathcal{O}$, endowed with the standard norms denoted by
$\|\cdot\|_{L^p}$ and $\|\cdot\|_{W^{m,p}}$, respectively. For non-integer
$\theta>0$, we use $W^{\theta,p}(\mathcal{O})$ 
to denote the Slobodeckij 
space, equipped with the classical norm $\|\cdot\|_{W^{\theta,p}}$. For
$\theta>0$, $1<p<\infty$, and $1\le q\le\infty$, let
$B^\theta_{p,q}(\mathcal{O})$ be the corresponding Besov space with the usual norm denoted by
$\|\cdot\|_{B^\theta_{p,q}}$; see for instance, \cite{Triebel}. Finally, for
$\theta\in(1+\frac{1}{p},2)$, we denote by
$B^\theta_{p,q,\mathbf n}(\mathcal{O})$ the Besov space satisfying the
homogeneous Neumann boundary condition associated with the outward unit
normal vector $\mathbf n$, i.e.
$
B^\theta_{p,q,\mathbf{n}}(\mathcal{O})
:=\bigl\{\psi\in B^\theta_{p,q}(\bo):\mathbf{n}\cdot\nabla\psi|_{\partial\mathcal{O}}=0\bigr\}$; see, e.g., \cite[Example~A.5]{Agresti}.
We recall from \cite[P.~310, and  P.~328, Remark~1]{Triebel} that for $0\le\theta\le\beta<\infty$ and $1<p\le q<\infty$,
\begin{equation}\label{equa1.3}
\begin{split}
&W^{\theta,q}(\mathcal{O})=B^\theta_{q,q}(\mathcal{O})
\quad\text{if}\quad\theta\in(0,\infty)\setminus\mathbb N\text{ and }1<q<\infty,\\
&W^{\beta,p}(\mathcal{O})\hookrightarrow W^{\theta,q}(\mathcal{O})
\quad\text{if}\quad\beta-\tfrac{2}{p}>\theta-\tfrac{2}{q},
\end{split}
\end{equation}
where $\hookrightarrow$ denotes continuous embedding. As a consequence of
\eqref{equa1.3}, we have
\begin{equation}\label{equa1.5}
	B^\frac{3}{2}_{4,4}(\mathcal{O})=W^{\frac{3}{2},4}(\mathcal{O})\hookrightarrow W^{1,4}(\mathcal{O})\hookrightarrow W^{1,2}(\mathcal{O})\hookrightarrow L^p(\mathcal{O})\qquad\forall p{\in[2,\infty)}.
\end{equation}
Let $A_p$ be the operator $-\Delta+\mathrm{Id}$ in
$L^p(\mathcal{O})$, $p>1$,  with Neumann boundary condition, i.e.
\[
A_p\varphi:=-\Delta\varphi+\varphi,\qquad
D(A_p):=\bigl\{\psi\in W^{2,p}(\mathcal{O}):\mathbf{n}\cdot\nabla\psi|_{\partial\mathcal{O}}=0\bigr\}.
\]
Since {$A_p$ is sectorial with $0\in\rho(A_p)$,} 
fractional powers $A_p^\theta$, $\theta\ge 0$, whose domain satisfies the
embedding
\begin{equation}\label{equa2.1*}
	D(A_p^\theta)\hookrightarrow W^{1,p}(\mathcal{O})\quad\text{if }\theta>\tfrac12,
\end{equation}
for every $p>1$ (see, e.g., \cite[Section~2]{Horstmann}). By
\cite[Theorem~2.3]{Denk}, $A_p$ admits a bounded $H^\infty$-calculus in
$L^p(\mathcal{O})$; consequently, by \cite[Theorem~16.11]{Yagi},
\begin{equation}\label{equa1.7}
	\begin{aligned}
		D(A_p^\frac{\theta}{2})&={H^{\theta,p}}(\mathcal{O})
		&&\text{if }0\le\theta<1+\tfrac1p,\\
		D(A_p^\frac{\theta}{2})&=\bigl\{\psi\in {H^{\theta,p}}(\mathcal{O}):
		\mathbf{n}\cdot\nabla\psi|_{\partial\mathcal{O}}=0\bigr\}
		&&\text{if }1+\tfrac1p<\theta\le 2.
	\end{aligned}
\end{equation}
{Here, $H^{\theta,p}(\mathcal O)$ denotes the Bessel potential space; see, e.g., \cite[Definition 4.2.1 and Definition 2.3.1]{Triebel} for its definition. We recall from \cite[Theorem 4.6.1]{Triebel} that  for $1<p<\infty$ and $0\leq \theta<\infty$,
\begin{equation}\label{equ2.5}
B^\theta_{p,\min(2,p)}(\bo)\hookrightarrow H^{\theta,p}(\bo)\hookrightarrow B^\theta_{p,\max(2,p)}(\bo).
\end{equation}
}
For $t\ge 0$, we denote by $e^{-tA_p}$ the semigroup generated by $-A_p$ and
by $e^{t\Delta}$ the heat semigroup generated by $\Delta$. We recall from
\cite[Section~2]{Horstmann} that $e^{-tA_p}$ is independent of $p$, in the
sense that
\begin{equation}\label{equa1.8*}
	e^{-tA_p}\varphi=e^{-tA_q}\varphi
	\quad\forall\varphi\in L^p(\mathcal{O})\cap L^q(\mathcal{O}),\ p,q>1,
\end{equation}
and that there exist constants $\bk>0$ and $\lambda>0$, depending only on $\mathcal{O}${, $\theta$, $p$, and $T$}  such that, for every $\theta\ge 0$, $p>1$ and $t{\in(0,T]}$,
\begin{equation}\label{equa2.4*}
		\bigl\|A_p^\theta e^{-tA_p}\varphi\bigr\|_{L^p}
		\le\bk t^{-\theta}e^{-\lambda t}\|\varphi\|_{L^p},\text{ and }
		\bigl\|A_p^\theta e^{t\Delta}\varphi\bigr\|_{L^p}
		\le\bk t^{-\theta}\|\varphi\|_{L^p}.
\end{equation}
We shall also make use of the following Gagliardo--Nirenberg interpolation inequality
\begin{equation}\label{equa2.3}
	\|u\|_{L^4}\le\bk\|u\|_{W^{1,2}}^\frac{1}{2}\|u\|_{L^2}^\frac{1}{2},\qquad u\in W^{1,2}(\mathcal{O}).
\end{equation}
For any Lebesgue, Sobolev or Besov space $E$, we set
$
E_+:=\{u\in E:u(x)\ge 0\ \text{for almost every }x\in\mathcal{O}\}.
$
Since $E_+$ is a closed subset of $E$, it is a complete metric space with
respect to the metric $d_E(u,v):=\|u-v\|_E$.
Let $X$ and $Y$ be Banach spaces. We denote by $\mathcal{L}(Y;X)$ the space of
bounded linear operators from $Y$ to $X$, and by $\|\cdot\|_X$ the norm on
$X$. The duality pairing between $X'$ and $X$ is denoted by
${}_{X'}\langle\cdot,\cdot\rangle_X$. If $X$ is a Hilbert space, we write
$(\cdot,\cdot)_X$ for its inner product; when $X=L^2(\mathcal{O})$ we simply
write $(\cdot,\cdot)$. For $\gamma$-radonifying operators from
$L^2(\mathcal{O})$ to $X$ we use the notation $\gamma(L^2(\mathcal{O}),X)$;
recall that $\Psi\in\gamma(L^2(\mathcal{O}),X)$ if
\begin{equation*}
\|\Psi\|_{\gamma(L^2,X)}^2
:=\tilde{\mathbb{E}}\Bigl\|\sum_{k=1}^\infty\gamma_k\Psi(e_k)\Bigr\|_X^2<\infty,
\end{equation*}
where $\{\gamma_k\}_{k\ge 1}$ are independent standard normal random variables
on some probability space $(\tilde\Omega,\tilde{\mathcal{F}},\tilde{\mathbb{P}})$.
We refer to \cite[Section~2]{Neerven} for further details. 
For $1\le p<\infty$, $L^p(0,T;X)$ denotes the space of (equivalence classes of)
Bochner measurable functions $u:[0,T]\to X$ with
$\Int_0^T\|u(s)\|_X^p\,ds<\infty$. The space $W^{1,p}(0,T;X)$ consists of all
$u\in L^p(0,T;X)$ with $\frac{du}{dt}\in L^p(0,T;X)$. For $\theta\in(0,1)$,
$W^{\theta,p}(0,T;X)$ consists of all $u\in L^p(0,T;X)$ with finite norm
\[
\|u\|^p_{W^{\theta,p}(0,T;X)}
:=\int_0^T\|u(s)\|_X^p\,ds
+\int_0^T\int_0^T\frac{\|u(t)-u(s)\|_X^{{p}}}{|t-s|^{1+\theta p}}\,dt\,ds.
\]
We write $C([0,T];X)$ for the space of continuous functions
$u:[0,T]\to X$, and $X_{\text{weak}}$ for $X$ endowed with the weak topology. 
The space $C([0,T];X_{\text{weak}})$ consists of all 
{function $u:[0,T]\to X$} such
that $t\mapsto{}_{X'}\langle\psi,u(t)\rangle_X$ is continuous for every
$\psi\in X'$. Given a filtered probability space
$(\tilde\Omega,\tilde{\mathcal{F}},\tilde{\mathbb{P}},\tilde{\mathbb{F}}_T)$
and $1\le q<\infty$, $L^q_{\tilde{\mathbb{F}}_T}(\tilde\Omega;L^p(0,T;X))$
(resp.\ $L^q_{\tilde{\mathbb{F}}_T}(\tilde\Omega;C([0,T];X))$) denotes the
space of $\tilde{\mathbb{F}}_T$-
{progressively measurable} processes $(u(t))_{t\in[0,T]}$ with
\[
\tilde{\mathbb{E}}\Bigl(\int_0^T\|u(s)\|_X^p\,ds\Bigr)^{q/p}<\infty
\quad\Bigl(\text{resp. }\tilde{\mathbb{E}}\sup_{s\in[0,T]}\|u(s)\|_X^q<\infty\Bigr).
\]
We denote by
$\mathcal{L}_2^0(\mathcal{O})$ (resp.\ $\mathcal{L}_2^1(\mathcal{O})$) the
space of Hilbert--Schmidt operators from $L^2(\mathcal{O})$ to $L^2(\mathcal{O})$
(resp.\ to $W^{1,2}(\mathcal{O})$); see \cite[Appendix~C]{DaPrato}. {It is well known that, in general, the cylindrical Wiener process defined by \eqref{equ1.3} is not an \(L^2(\bo)\)-valued process. Nevertheless, there always exists a larger Hilbert space $U_1$ such that $L^2(\bo)\hookrightarrow U_1$ with Hilbert--Schmidt embedding $\bj$, allowing us to identify $W_i$, $i\in\{1,2\}$, as a $Q$-Wiener process on $U_1$ for some trace-class operator $Q$; see, for instance, \cite[Subsection 2.5.1]{Liu}. Throughout the rest of the paper, we implicitly work with this extension and simply refer to $W_i$, $i\in\{1,2\}$, as cylindrical Wiener processes on $L^2(\bo)$. The same convention applies to stochastic integration with respect to $W_i$. More precisely, for every $\Phi\in L^2_{\mathbb{F}_T}(\Omega;L^2(0,T;\mathcal{L}_2^\kappa(\mathcal{O})))$, $\kappa\in\{0,1\}$, we use the notation
\begin{equation*}
\int_0^t\Phi(s)dW_i(s):=\int_0^t\Phi(s)\circ \bj^{-1}dW_i(s), \quad t\in [0,T].
\end{equation*}
This definition is well defined and independent of the choice of $U_1$ and $\bj$; see, e.g., \cite[Subsection 2.5.2]{Liu}.
In the remainder of the paper,} $\mathcal{K}$ denotes a positive constant that may change
from line to line.

We now collect the assumptions on the initial data and noise coefficients used throughout the paper. Before doing so, we recall from \cite[Theorem~3.5 and Section~4]{Neerven} that the operator $A_4$ has $L^4$-maximal regularity. Consequently, the stochastic convolution map
\[
\Phi:G\mapsto\Phi(G):=\int_0^\cdot e^{-(\cdot-s)A_4}G(s)\,dW_2(s)
\]
belongs to
$\mathcal{L}(L^4_{\mathbb{F}_T}(\Omega;L^4(0,T;\gamma(L^2(\mathcal{O});D(A_4^\frac{1}{2}))));
L^4_{\mathbb{F}_T}(\Omega;L^4(0,T;D(A_4))))$, and following the notation of \cite[Section~4]{Neerven}, we denote by
$\mathcal{K}_4^\Diamond$ its operator norm, i.e.
\begin{equation}\label{equa1.9*}
	\mathcal{K}_4^\Diamond
	:=\|\Phi\|_{\mathcal{L}(L^4_{\mathbb{F}_T}(\Omega;L^4(0,T;\gamma(L^2(\mathcal{O});D(A_4^\frac{1}{2}))));
		L^4_{\mathbb{F}_T}(\Omega;L^4(0,T;D(A_4))))}.
\end{equation}

\begin{assumption}\label{assum2.2}
$ $\\
$\bullet$(A$_1$) The initial data are non-negative $\mathcal{F}_0$-measurable random
		variables $n_0:\Omega\to (L^2(\mathcal{O}))_+$ and
		$c_0:\Omega\to (B_{4,4,\mathbf{n}}^\frac{3}{2}(\mathcal{O}))_+$ such that
		\begin{equation}\label{equa1.8}
			n_0\in L^{32}(\Omega{,\mathcal F_0};(L^2(\mathcal{O}))_+) \text{ and }
			c_0\in L^{32}(\Omega{,\mathcal F_0};(B_{4,4,\mathbf{n}}^\frac{3}{2}(\mathcal{O}))_+),
		\end{equation}
$\bullet$(A$_2$) The map $g_1:L^2(\mathcal{O})\to\mathcal{L}_2^0(\mathcal{O})$ is
		such that there exist   $L_0>0$ and	$L_1>0$ satisfying
		\begin{equation}\label{equa2.3*}
			\|g_1(u)\|_{\mathcal{L}_2^0}\le L_0\quad\forall u\in L^2(\mathcal{O}),
			\qquad
			\|g_1(u)-g_1(v)\|_{\mathcal{L}_2^0}\le L_1\|u-v\|_{L^2}\quad\forall u,v\in L^2(\mathcal{O}).
		\end{equation}
$\bullet$(A$_3$) The map $g_2:W^{1,2}(\mathcal{O})\to\mathcal{L}_2^1(\mathcal{O})$
		is such that there exists $L_2>0$ satisfying
		\begin{equation}\label{equa2.4}
			\begin{aligned}
				\|g_2(u)-g_2(v)\|_{\mathcal{L}_2^0}&\le L_2\|u-v\|_{L^2}
				&&\forall u,v\in W^{1,2}(\mathcal{O}),\\
				\|g_2(u)-g_2(v)\|_{\mathcal{L}_2^1}&\le L_2\|u-v\|_{W^{1,2}}
				&&\forall u,v\in W^{1,2}(\mathcal{O}),
			\end{aligned}
		\end{equation}
		and there exists a square-summable real sequence $\{z_k\}_{k\ge 1}$ such that,
		for all $u\in W^{1,2}(\mathcal{O})$ and any $k\ge 1$,
		\begin{equation}\label{equa2.6}
			\begin{aligned}
			&	|[g_2(u)(e_k)](x)|\le|z_k|\,|u(x)|
				&&\text{for almost every }x\in\mathcal{O},\\
			\end{aligned}
		\end{equation}
		Moreover, there exists $L_3>0$ such that the restriction
		$g_2:D(A_4^\frac{1}{2})\to\gamma(L^2(\mathcal{O}),D(A_4^\frac{1}{2}))$ satisfies
		\begin{equation}\label{equa2.4**}
			\begin{aligned}
			&	\|g_2(u)\|_{\gamma(L^2,D(A_4^\frac{1}{2}))}
				\le L_3\|u\|_{D(A_4^\frac{1}{2})}
				&&\forall u\in D(A_4^\frac{1}{2}),\\
			&	\|g_2(u)-g_2(v)\|_{\gamma(L^2,D(A_4^\frac{1}{2}))}
				\le\frac{1}{\mathcal{K}_4^\Diamond+1}\|u-v\|_{D(A_4)}
				&&\forall u,v\in D(A_4),
			\end{aligned}
		\end{equation}
		where $\mathcal{K}_4^\Diamond>0$ is given by \eqref{equa1.9*}.
\end{assumption}

\begin{remark}\label{rem:example}
	Inequality \re{equa2.6} is essential for the non-negativity argument in Lemma~\ref{lem4.3} below. {As a direct consequence of \re{equa2.6}, for every $u\in L^2(\mathcal O)$ and $k\geq1$,
	\begin{equation}\label{equ2.14}
		[g_1(u)(e_k)](x)=0
		\quad\text{for almost every }x\in\{y\in \bo: u(y)=0\}.
	\end{equation}}
 An example of maps $g_1$ and $g_2$ satisfying (A$_2$) and (A$_3$) is obtained by fixing $f\in L^2(\bo)$, a square-summable real sequence $\{z_k\}_{k\ge 1}$, and setting
	\[
	g_1(u)(v):=\frac{(f,v)}{\sqrt{\norm{u}_{L^2}^2+1}}\,u,
	\qquad
	g_2(u)(v):=\sum_{k=1}^\infty z_k\,(v,e_k)_{L^2}\,u.
	\]
	The verification of (A$_2$)--(A$_3$) is elementary: the Hilbert--Schmidt norms follow from the summability of $\{z_k\}$ and the boundedness of the scalar factor, the Lipschitz bounds follow from the Lipschitz continuity of $s\mapsto(1+s^2)^{-1/2}$ and the inequality $\abs{(v,e_k)}\le\norm{v}_{L^2}$, and the pointwise inequalities \re{equa2.6} hold by construction. The smallness condition in \re{equa2.4**}$_2$ is satisfied by choosing $\norm{z}_{\ell^2}$ small enough. 
\end{remark}

\subsection{Main results}
\label{subsec:mainresults}

We now introduce the notions of martingale and weak solution for
\eqref{equa1.1}--\eqref{equa1.2}. Throughout, we use the abbreviations
\begin{equation}\label{equa1.16*}
	\begin{aligned}
		\mathcal{Z}_1&:=C([0,T];(W^{3,2}(\mathcal{O}))')\cap C([0,T];L^2_{\mathrm{weak}}(\mathcal{O}))\cap L^2(0,T;W^{1,2}(\mathcal{O})),\\
		\mathcal{Z}_2&:=C([0,T];L^2(\mathcal{O}))\cap C([0,T];W^{1,2}_{\mathrm{weak}}(\mathcal{O}))\cap L^2(0,T;D(A_2)),
	\end{aligned}
\end{equation}
which will appear repeatedly below.

\begin{definition}\label{defi1.3}
	Let $\vartheta_1$ and $\vartheta_2$ be Borel probability measures on the
	closed cones $(L^2(\mathcal{O}))_+$ and $(B^\frac{3}{2}_{4,4,\mathbf{n}}(\mathcal{O}))_+$,
	respectively, with finite $32$-nd moments, i.e.
	\begin{equation}\label{equa1.17}
		\int\|u\|_{L^2}^{32}\,d\vartheta_1<\infty,
		\qquad
		\int\|v\|^{32}_{B^\frac{3}{2}_{4,4}}\,d\vartheta_2<\infty.
	\end{equation}
	A \emph{martingale solution} to \eqref{equa1.1}--\eqref{equa1.2} with initial
	law $(\vartheta_1,\vartheta_2)$ is a quadruple
	\[
	\bigl((\bar{\Omega},\bar{\mathcal{F}},\bar{\mathbb{P}},\bar{\mathbb{F}}_T),
(	\bar{W}_1,\bar{W}_2),\bar{n},\bar{c}\bigr),
	\]
such that:
	\begin{itemize}
		\item $(\bar{\Omega},\bar{\mathcal{F}},\bar{\mathbb{P}},\bar{\mathbb{F}}_T)$
		is a filtered probability space whose filtration $\bar{\mathbb{F}}_T$ satisfies
		the usual conditions;
		\item $\bar{W}_1$ and $\bar{W}_2$ are independent $L^2(\mathcal{O})$-cylindrical
		Wiener processes over
		$(\bar{\Omega},\bar{\mathcal{F}},\bar{\mathbb{P}},\bar{\mathbb{F}}_T)$;
		\item $\bar{n}$ and $\bar{c}$ are $\bar{\mathbb{F}}_T$-adapted, with
		$\bar{\mathbb{P}}$-a.s. paths in $\mathcal{Z}_1$ and $\mathcal{Z}_2$,
		respectively, and such that
		\begin{equation}\label{equa1.18}
			\mathrm{Law}_{\bar{\mathbb{P}}}(\bar{n}(0))=\vartheta_1,
			\qquad
			\mathrm{Law}_{\bar{\mathbb{P}}}(\bar{c}(0))=\vartheta_2,
		\end{equation}
		and
		\begin{equation}\label{equa1.19}
			\bar{n}\in L^4_{\bar{\mathbb{F}}_T}(\bar{\Omega};L^2(0,T;L^2(\mathcal{O}))),
			\qquad
			\bar{c}\in L^4_{\bar{\mathbb{F}}_T}(\bar{\Omega};L^2(0,T;D(A_2)));
		\end{equation}
		\item for all $t\in[0,T]$,
		\begin{equation}\label{equa1.20}
			\begin{split}
			&	(\bar{n}(t),\varphi)+\int_0^t(\nabla\bar{n}(s),\nabla\varphi)\,ds
				-\chi\int_0^t(\bar{n}(s)\nabla\bar{c}(s),\nabla\varphi)\,ds\\
				&=(\bar{n}(0),\varphi)
				+\int_0^t\bigl(\alpha\bar{n}(s)-\mu\bar{n}^2(s),\varphi\bigr)\,ds+\int_0^t(g_1(\bar{n}(s)),\varphi)\,d\bar{W}_1(s)\quad \forall \varphi\in D(A_2),\\[1mm]
				&(\bar{c}(t),\psi)+\int_0^t(A_2\bar{c}(s),\psi)\,ds\\
				&=(\bar{c}(0),\psi)
				+\int_0^t(\bar{n}(s),\psi)\,ds+\int_0^t(g_2(\bar{c}(s)),\psi)\,d\bar{W}_2(s)\quad \forall\psi\in L^2(\mathcal{O}),
			\end{split}
		\end{equation}
	$\bar{\mathbb{P}}$-almost surely.
	\end{itemize}
\end{definition}

\begin{definition}\label{defi1.4}
Let $(n_0,c_0)$ be a random
	variable satisfying (A$_1$) of Assumption~\ref{assum2.2}. A \emph{weak
		solution} of \eqref{equa1.1}--\eqref{equa1.2} with initial data $(n_0,c_0)$ is
	a pair $(n,c)$ of $\mathbb{F}_T$-adapted processes, with $\mathbb{P}$-almost surely paths
	in $\mathcal{Z}_1$ and $\mathcal{Z}_2$, respectively, such that
	\begin{equation}\label{equa1.21}
		n\in L^4_{\mathbb{F}_T}(\Omega;L^2(0,T;L^2(\mathcal{O}))),
		\qquad
		c\in L^4_{\mathbb{F}_T}(\Omega;L^2(0,T;D(A_2))),
	\end{equation}
and satisfying the stochastic variational system \eqref{equa1.20} for all
$t\in[0,T]$, $\mathbb{P}$-almost surely, with $\bar n(0)$, $\bar c(0)$,
$\bar W_1$, and $\bar W_2$ replaced by $n_0$, $c_0$, $W_1$, and $W_2$,
respectively.
\end{definition}

We now state the main results of the paper. The first is the existence of a
martingale solution.

\begin{theorem}\label{theo1.5}
	{Let $\chi>0$, $\alpha\in\mathbb R$ and $\mu>0$, and} {a}ssume (A$_2$) and (A$_3$) of Assumption~\ref{assum2.2}. For any Borel
	probability measures $\vartheta_1$ and $\vartheta_2$ on
	$(L^2(\mathcal{O}))_+$ and $(B^\frac{3}{2}_{4,4,\mathbf{n}}(\mathcal{O}))_+$ with finite
	$32$-nd moments, i.e. satisfying \eqref{equa1.17}, there exists at least one
	martingale solution
	\[
	\bigl((\bar{\Omega},\bar{\mathcal{F}},\bar{\mathbb{P}},\bar{\mathbb{F}}_T),
	(\bar{W}_1,\bar{W}_2),\bar{n},\bar{c}\bigr)
	\]
	to \eqref{equa1.1}--\eqref{equa1.2} which satisfies the non-negativity
	property
	\begin{equation}\label{equa1.23}
	\begin{split}
&	\bar{\mathbb{P}}\bigl\{\omega\in\bar{\Omega}:
\bar{n}(\omega,t,x)\ge 0
\ \text{for almost every  }(t,x)\in[0,T]\times\mathcal{O}\bigr\}=1,\\
&	\bar{\mathbb{P}}\bigl\{\omega\in\bar{\Omega}:
 \bar{c}(\omega,t,x)\ge 0
\ \text{for almost every }(t,x)\in[0,T]\times\mathcal{O}\bigr\}=1.
	\end{split}
	\end{equation}
\end{theorem}

The second result concerns pathwise uniqueness, uniqueness in law, and the
existence of a unique weak solution.

\begin{theorem}\label{theo1.6}
	{Let $\chi>0$, $\alpha\in\mathbb R$ and $\mu>0$, and s}uppose (A$_2$) and (A$_3$) of Assumption~\ref{assum2.2} hold. Let
	$\vartheta_1$ and $\vartheta_2$ be Borel probability measures on
	$(L^2(\mathcal{O}))_+$ and $(B^\frac{3}{2}_{4,4,\mathbf{n}}(\mathcal{O}))_+$,
	respectively, with finite $32$-nd moments. Let
	\[
	\bigl((\bar{\Omega},\bar{\mathcal{F}},\bar{\mathbb{P}},\bar{\mathbb{F}}_T),
(	\bar{W}_1,\bar{W}_2),\bar{n}_i,\bar{c}_i\bigr),\qquad i=1,2,
	\]
	be martingale solutions to \eqref{equa1.1}--\eqref{equa1.2} {on the same stochastic basis and with the same Wiener processes, such that $\bar n_1(0)=\bar n_2(0)$ and $\bar c_1(0)=\bar c_2(0)$ $\bar{\mathbb P}$-a.s.\ (in particular} with the same
	initial law $(\vartheta_1,\vartheta_2)${)}. Then: 
	\begin{itemize}
		\item[(i)] (\emph{Pathwise uniqueness}) The two solutions coincide, i.e.
		\[
		\bar{\mathbb{P}}\bigl\{\omega\in\bar{\Omega}:
		(\bar{n}_1(\omega,t),\bar{c}_1(\omega,t))
		=(\bar{n}_2(\omega,t),\bar{c}_2(\omega,t))
		\ \forall t\in[0,T]\bigr\}=1.
		\]
		\item[(ii)] (\emph{Uniqueness in law}) If
		$((\bar{\Omega}_i,\bar{\mathcal{F}}^i,\bar{\mathbb{P}}_i,\bar{\mathbb{F}}_T^i),
		(\bar{W}_1^i,\bar{W}_2^i),\bar{n}_i,\bar{c}_i)$, $i=1,2$, are martingale
		solutions with the same initial  {law $(\vartheta_1,\vartheta_2)$}, then
		\[
		\mathrm{Law}_{\bar{\mathbb{P}}_1}(\bar{n}_1,\bar{c}_1)
		=\mathrm{Law}_{\bar{\mathbb{P}}_2}(\bar{n}_2,\bar{c}_2)
		\quad\text{on }\mathcal{Z}_1\times\mathcal{Z}_2,
		\]
		where $\mathcal{Z}_1$ and $\mathcal{Z}_2$ are defined in \eqref{equa1.16*}.
		\item[(iii)] (\emph{Existence of a unique weak solution}) If, in addition,
		$(n_0,c_0)$ satisfies (A$_1$) of Assumption~\ref{assum2.2}, then
		\eqref{equa1.1}--\eqref{equa1.2} admits a unique weak solution $(n,c)$ with
		initial data $(n_0,c_0)$. Moreover, $(n,c)$ inherits the non-negativity
		property, i.e.
		\begin{equation}\label{equa1.28}
\begin{split}
&			\mathbb{P}\bigl\{\omega\in\Omega:
n(\omega,t,x)\ge 0
 \text{ for almost every  }(t,x)\in[0,T]\times\mathcal{O}\bigr\}=1,\\
		&	\mathbb{P}\bigl\{\omega\in\Omega:
c(\omega,t,x)\ge 0
\ \text{for almost every  }(t,x)\in[0,T]\times\mathcal{O}\bigr\}=1.
\end{split}
		\end{equation}
	\end{itemize}
\end{theorem}
\section{Proof of Theorem \ref{theo1.5}}{\label{sec:proof1}}
This section is devoted to the proof of the first main result of the paper. It is organised as follows. In Subsection \ref{subsec:reg}, we introduce the regularised stochastic variational system with a regularisation parameter $\eps\in(0,1)$  and state the global existence of its solution. The full construction of the regularised solution is deferred to Appendix \ref{app:reg}. In Subsection \ref{subsec:nonneg}, we establish the non-negativity of the regularised solutions. In Subsection \ref{subsec:apriori}, we derive $\eps$-independent uniform estimates and establish the tightness of the probability laws of the regularised solutions. Finally, in Subsection \ref{subsec:limit}, we pass to the limit along a suitable sequence $\eps_\ell\xrightarrow{\ell\uparrow\infty}0$, thereby completing the proof of Theorem \ref{theo1.5}.
To begin, let $\vartheta_1$ and $\vartheta_2$ be two fixed Borel probability measures on ${(}L^2(\bo){)_+}$ and ${(}B^\frac{3}{2}_{4,4,\mathbf{n}}(\bo){)_+}$, respectively, satisfying \re{equa1.17}.

\subsection{The regularised system}
\label{subsec:reg}

On the probability space $(\Omega,\mathcal{F}_0,\mathbb{P})$, we consider two random variables $u_0:(\Omega,\mathcal{F}_0)\to ((L^2(\bo))_+,\mathcal{B}((L^2(\bo))_+))$, and $v_0:(\Omega,\mathcal{F}_0)\to ((B_{4,4,\mathbf{n}}^\frac{3}{2}(\bo))_+,\mathcal{B}((B_{4,4,\mathbf{n}}^\frac{3}{2}(\bo))_+))$ such that
\begin{equation}\label{equa2.2}
	\text{Law}_{\mathbb{P}}(u_0)=\vartheta_1\  \text{ and } \  \text{Law}_{\mathbb{P}}(v_0)=\vartheta_2.
\end{equation}
Since the thirty-second moments of $\vartheta_1$ and $\vartheta_2$ are respectively finite on $(L^2(\bo))_+$ and $(B^\frac{3}{2}_{4,4,\mathbf{n}}(\bo))_+$, we infer that
\begin{equation}\label{equa2.1}
	u_0\in L^{32}_{\mathcal{F}_0}(\Omega;L^2(\bo)), \quad\text{ and }\quad v_0\in L^{32}_{\mathcal{F}_0}(\Omega; B_{4,4,\mathbf{n}}^\frac{3}{2}(\bo)).
\end{equation}
Let $\eps\in (0,1)$ be arbitrary fixed.  The regularised stochastic variational problem consists in finding, on the filtered probability space $(\Omega, \mathcal{F}, \mathbb{P}, \mathbb{F}_T)$, a $\mathbb{F}_T$-progressively measurable, $L^2(\bo)\times B_{4,4,\mathbf{n}}^\frac{3}{2}(\bo)$-valued process $\{(n_\eps(t),c_\eps(t))\}_{t\in [0,T]}$ satisfying for all $t\in [0,T]$,
\begin{equation}\label{equa3.1}
	\begin{split}
		&	(n_\eps(t),\varphi)+\int_0^t( \nabla n_\eps(s),\nabla\varphi)ds-\chi \int_0^t(F_\eps(n^+_\eps(s))\nabla c_\eps(s),\nabla \varphi )ds-(u_0,\varphi)\\
		&=\int_0^t(\alpha n_\eps(s)-\mu [H_\eps(n^+_\eps(s))]^2,\varphi)ds+\int_0^t(g_1(n^+_\eps(s)),\varphi)dW_1(s) \quad\forall \varphi\in W^{1,2}(\bo), \  \mathbb{P}\text{-a.s.}, \\
		&	(c_\eps(t),\psi)+\int_0^t(A_2 c_\eps(s),\psi)ds=(v_0,\psi)+\int_0^t(H_\eps(n^+_\eps(s)),\psi)ds+\int_0^t(g_2(c_\eps(s)),\psi)dW_2(s),
	\end{split}
\end{equation}
for all $ \psi\in L^2(\bo)$ and $\mathbb{P}$-almost surely, where $F_\eps:{[0,\infty)}\to \br$ and $H_\eps:\br\to\br_+$ are functions defined respectively for any $z\in \br$ by
\begin{equation*}
	F_\eps(z):=\frac{z}{1+\eps z} \quad \text{ and }\quad 	H_\eps(z):=\begin{cases}
		0\text{ if } z\leq 0,\\
		-\frac{1}{\eps^2}z^3+\frac{2}{\eps}z^2\text{ if } 0<z<\eps,\\
		z\text{ if }\eps\leq z<\frac{1}{\eps},\\
		z-(z-\frac{1}{\eps})^3+(z-\frac{1}{\eps})^2\text{ if }\frac{1}{\eps}\leq z<\frac{1}{\eps}+1,\\
		\frac{1}{\eps}+1\text{ if } z{\geq} \frac{1}{\eps}+1,
	\end{cases}
\end{equation*}
and $z^+:=\max(0,z)$ for any $z\in \br$. Let us note that the choice of {the $C^1$} functions $F_\eps$ and $H_\eps$ ensures that 
\begin{itemize}
	\item For any $z{\geq 0}$ {(for $F_\eps$) and any $z\in\br$ (for $H_\eps$)}
	\begin{equation}\label{equa3.2}
		0<F'_\eps(z)\leq 1 \ \text{ and }\ 	\abs{H'_\eps(z)}\leq \frac{4}{3}.
	\end{equation}
	\item For any $z\in \br$,
	\begin{equation}\label{equa3.3}
		\lim_{\eps\longrightarrow 0}H_\eps(z)= z^+{,\ 0\le H_\eps(z)\le z^+}\ \text{ and }\	{0\le}F_\eps(z)\leq {\min\bigl(z,}\frac{1}{\eps}{\bigr)}\ \text{ if }   z\geq 0.
	\end{equation}
	\item Finally, there exists a constant $\bk_\eps>0$ such that
	\begin{equation}\label{equa3.4}
		\abs{H_\eps(z)}\leq \bk_\eps, \text{ for all } z\in \br.
	\end{equation}
\end{itemize}
We are now in  position to state the  existence of a solution to the regularised system \re{equa3.1}. 

\begin{proposition}\label{propo3.5}
	Assume that items (A$_2$) and (A$_3$) of Assumption \ref{assum2.2} hold, and let $u_0$ and $v_0$ satisfy \re{equa2.1}. Then, for every $\eps\in(0,1)$, the regularised stochastic variational system \re{equa3.1} admits a global solution $(n_\eps,c_\eps)$, i.e.\ $(n_\eps,c_\eps)$ is $\mathbb{F}_T$-progressively measurable with $\mathbb{P}$-a.s.\ paths in $L^2(\bo)\times B_{4,4,\mathbf{n}}^\frac{3}{2}(\bo)$, satisfies \re{equa3.1} for all $t\in [0,T]$, and enjoys the regularity
	\begin{equation}\label{equa3.7*}
		\begin{split}
			&n_\eps\in L^4_{\mathbb{F}_T}(\Omega; C([0,T]; L^2(\bo)))\cap L^4_{\mathbb{F}_T}(\Omega; L^2(0,T; W^{1,2}(\bo))),\text{ and}\\
			&c_\eps\in L^4_{\mathbb{F}_T}(\Omega; C([0,T]; B_{4,4}^\frac{3}{2}(\bo)))\cap L^4_{\mathbb{F}_T}(\Omega; L^4(0,T; D(A_2))).
		\end{split}
	\end{equation}
\end{proposition}

Since the proof of Proposition~\ref{propo3.5} is rather technical and lengthy, we defer it to Appendix~\ref{app:reg} for the sake of readability and to keep the main arguments focused. The proof relies on the construction of a maximal local solution via the Banach fixed-point theorem, followed by a priori estimates that rule out finite-time blow-up. 

\subsection{Non-negativity of the regularised solutions}
\label{subsec:nonneg}

In this subsection we show that the solutions of the regularised system are non-negative. We fix a sequence $\{\eps_\ell\}_{\ell\geq 1}$ of points of $(0,1)$ such that $\eps_\ell\xrightarrow{\ell\uparrow\infty}0$ and, for any $\ell\geq 1$, denote by $(n_{\eps_\ell},c_{\eps_\ell})$ the global solution to the system \re{equa3.1} obtained by Proposition \ref{propo3.5}.

\begin{lemma}\label{lem4.3}
	Suppose the assumptions of Proposition \ref{propo3.5} hold and let $\{(n_{\eps_\ell},c_{\eps_\ell})\}_{\ell\geq 1}$ be the family of global solutions to the system \re{equa3.1} obtained from Proposition \ref{propo3.5}. Then, for any $\ell\geq 1$, the couple $(n_{\eps_\ell},c_{\eps_\ell})$ satisfies
	\begin{equation}\label{equa4.19}
		\begin{split}
			&\mathbb{P}\left\{ \omega\in \Omega: n_{\eps_\ell}(t,\omega, x)\geq 0\ \text{ for all } t\in [0,T] \text{ and  for  almost every } x\in  \bo \right\}=1,\\
			&\mathbb{P}\left\{ \omega\in \Omega: c_{\eps_\ell}(t,\omega, x)\geq 0\ \text{ for all } t\in [0,T] \text{ and for  almost every } x\in  \bo \right\}=1.
		\end{split}
	\end{equation}
\end{lemma}

\begin{proof}
	Let   $\ell\geq 1$ and $t \in [0,T]$ be arbitrary fixed. For any $v\in \{n_{\eps_\ell} ,c_{\eps_\ell} \}$,	we define
	$
	v^-:=\max(-v,0).
	$
	To prove \re{equa4.19}, it is sufficient to show that
	$\norm{v^-(t)}_{L^2}=0$ almost surely. Since the functional $\Psi:L^2(\bo)\to \br$ defined by
	\begin{equation*}
		\Psi(u):=\int_{\bo}(u^-(x))^2dx, \ u\in L^2(\bo),
	\end{equation*}
	is not {twice} Fr{\'e}chet-differentiable{,} so that the It\^o formula {cannot} be applied {directly}, let us introduce suitable regular functionals approximating $\Psi$ by considering a family of $C^2$-functions $\{\Theta_\eta:\br\to \br\}_{\eta\in (0,1)}$ defined for any $\eta\in (0,1)$ by
	\begin{equation*}
		\Theta_\eta(s):=\begin{cases}
			-s, \hspace{3.6cm}\text{ if } s<-\eta,\\
			-\frac{3}{\eta^4}s^5-\frac{8}{\eta^3}s^4-\frac{6}{\eta^2}s^3,\hspace{1cm}\text{ if } -\eta\leq s<0,\\
			0, \hspace{3.9cm} \text{ if } s\geq 0.
		\end{cases}
	\end{equation*}
	Note that for any $\eta\in (0,1)$ and $s\in \br$,
	\begin{equation}\label{5.2*}
		\Theta^{''}_\eta(s)\xrightarrow{\eta{\downarrow} 0} 0, \ \Theta^{'}_\eta(s)\xrightarrow{\eta{\downarrow} 0}-1_{(-\infty,0)}(s),\text{ and }\Theta_\eta(s)\xrightarrow{\eta{\downarrow} 0} \Theta(s),
	\end{equation}
	where $1_{(-\infty,0)}$ is the indicator function of $(-\infty,0)$ and $\Theta(s):=s^-=\max(-s,0)$. In addition, there exists $\bk>0$ such that for any $\eta\in (0,1)$,
	\begin{equation}\label{5.2**}
		\abs{\Theta_\eta(s)}\leq \bk \abs{s}, \quad \abs{\Theta'_\eta(s)}\leq \bk\quad\text{ and }\quad \abs{\Theta_\eta(s)\Theta_\eta^{''}(s)}\leq \bk, \quad \text{ for any } s\in \br.
	\end{equation}
	Let $\eta\in (0,1)$ be fixed and let us introduce the functional $\Psi_\eta:L^2(\bo)\to \br$ defined by
	\begin{equation*}
		\Psi_\eta(u):=\int_{\bo}(\Theta_\eta(u(x)))^2dx,\quad  \forall u\in L^2(\bo).
	\end{equation*}
	It is twice Fr{\'e}chet-differentiable and for any $u,h\in L^2(\bo)$,
	\begin{equation*}
		\Psi'_\eta(u):=2\Theta_\eta(u)\Theta_\eta'(u),\text{ and }\Psi^{''}_\eta(u)(h):=2\Theta_\eta(u)\Theta_\eta^{''}(u)h+2(\Theta'_\eta(u))^2h \text{ in $L^2(\bo)$}.
	\end{equation*}
	By using \re{5.2**}, one may verify that $\Psi_\eta$ satisfies conditions (i)-(v) of \cite[Theorem 1.2]{Pardoux}. Then, recalling that the expectation of the stochastic integral vanishes, we apply the Itô formula in the form of \cite[Theorem~1.2]{Pardoux} to $\Psi_\eta(n_{\eps_\ell})$ to derive that
	\begin{equation*}
		\begin{split}
			\be\Psi_\eta(n_{\eps_\ell} (t))&=\be\Psi_\eta(u_0)-2\be\Int_0^t\left(\abs{\nabla n_{\eps_\ell}(s)}^2,\Theta_\eta(n_{\eps_\ell}(s))\Theta_\eta^{''}(n_{\eps_\ell}(s))\right)ds\\
			&\qquad-2\be\Int_0^t\left(\abs{ \nabla n_{\eps_\ell}(s)}^2,(\Theta_\eta^{'}(n_{\eps_\ell}(s)))^2\right)ds\\
			& \qquad+2\chi\be\Int_0^t \left(F_{\eps_\ell}(  n^+_{\eps_\ell}(s))\nabla  c_{\eps_\ell}(s)\cdot\nabla  n_{\eps_\ell}(s),\Theta_\eta(n_{\eps_\ell}(s))\Theta_\eta^{''}(n_{\eps_\ell}(s))\right)ds\\
			&\qquad+2\chi\be\Int_0^t\left( F_{\eps_\ell}(  n^+_{\eps_\ell}(s))\nabla  c_{\eps_\ell}(s)\cdot\nabla  n_{\eps_\ell}(s),(\Theta_\eta^{'}(n_{\eps_\ell}(s)))^2\right)ds\\
			&\qquad+2\be\int_0^{t}(\alpha n_{\eps_\ell}(s)-\mu [H_\eps(  n^+_{\eps_\ell}(s))]^2,\Theta_\eta(n_{\eps_\ell}(s))\Theta_\eta'(n_{\eps_\ell}(s)))ds\\
			& \qquad+\be\int_0^t\sum_{k=1}^\infty\left(g_1(  n^+_{\eps_\ell}(s))e_k,\Theta_\eta(n_{\eps_\ell}(s))\Theta_\eta^{''}(n_{\eps_\ell}(s))g_1(  n^+_{\eps_\ell}(s))e_k\right)ds\\
			&\qquad+\be\int_0^t\sum_{k=1}^\infty\left(g_1(  n^+_{\eps_\ell}(s))e_k,(\Theta_\eta^{'}(n_{\eps_\ell}(s)))^2g_1(  n^+_{\eps_\ell}(s))e_k\right)ds\\
			&:=B_1(\eta)+B_2(t,\eta)+B_3(t,\eta)+B_4(t,\eta)+B_5(t,\eta)+B_6(t,\eta)+B_7(t,\eta)+B_8(t,\eta).
		\end{split}
	\end{equation*}
	By definition of $\Theta_\eta$, we have
	\begin{equation*}
		\Theta_\eta(n_{\eps_\ell})\Theta_\eta^{''}(n_{\eps_\ell})=0\text{ and }\Theta_\eta^{'}(n_{\eps_\ell})=0\text{ on the set } \left\{ x\in \bo:  n_{\eps_\ell}(x)\geq 0\right\}.
	\end{equation*}
	This implies that $B_4(t,\eta)=0$ and $B_5(t,\eta)=0$. {In addition, keeping \eqref{equ2.14} and the definition of $H_{\eps_\ell}$ in mind, we derive that $B_7(t,\eta)=0$, $B_8(t,\eta)=0$}, 
	 and
	\begin{equation*}
		B_6(t,\eta)=2\alpha \be\int_0^{t}(n_{\eps_\ell}(s),\Theta_\eta(n_{\eps_\ell}(s))\Theta_\eta'(n_{\eps_\ell}(s)))ds.
	\end{equation*}
	By definition of $B_3(t,\eta)$, we have $B_3(t,\eta)\leq 0$; {moreover $n_{\eps_\ell}\,\Theta_\eta(n_{\eps_\ell})\Theta_\eta'(n_{\eps_\ell})\ge0$, so $B_6\le 2|\alpha|\,\be\Int_0^t\Int_\bo n_{\eps_\ell}(s,x)\Theta_\eta(n_{\eps_\ell}(s,x))\Theta_\eta'(n_{\eps_\ell}(s,x))dxds$}, {and} we infer that
	\begin{equation}\label{equa4.26}
		\be\Psi_\eta(n_{\eps_\ell} (t))\leq B_1(\eta)+B_2(t,\eta)+{2|\alpha|}\be\int_0^{t}\int_\bo n_{\eps_\ell}(s,x)\Theta_\eta(n_{\eps_\ell}(s,x))\Theta_\eta'(n_{\eps_\ell}(s,x))dxds.
	\end{equation}
	It follows from \re{5.2*} that for any arbitrary $x\in \bo$, $(\Theta_\eta(n_{\eps_\ell}(t,x)))^2\xrightarrow{\eta{\downarrow} 0} (n_{\eps_\ell}^-(t,x) )^2$. In addition, we deduce from \re{5.2**} that for any $\eta\in (0,1)$, $\abs{\Theta_\eta(n_{\eps_\ell}(t,x))}^2\leq\bk \abs{n_{\eps_\ell}(t,x)}^2$ and $n_{\eps_\ell}^2(t,\cdot)\in L^1(\bo)$. Hence by applying the dominated convergence theorem, we get
	\begin{equation*}
		\be\Psi_\eta(n_{\eps_\ell}(t))\xrightarrow{\eta{\downarrow} 0}\be\Psi(n_{\eps_\ell}(t))=\be\norm{n_{\eps_\ell}^-(t)}^2_{L^2}.
	\end{equation*}
Using a similar argument and the fact that $u_0\geq 0$, we also have
	$
	B_1(\eta)\xrightarrow{\eta{\downarrow} 0}{\be}\Psi(u_0)=0.
	$
	Noting that
	\begin{equation*}
		\begin{split}
			B_2(t,\eta)&={-2}\be\int_0^t\int_{\bo}1_{\{ -\eta\leq n_{\eps_\ell}(s)<0\}}\Theta_\eta(n_{\eps_\ell} (s,x))\Theta_\eta^{''}(n_{\eps_\ell}(s,x))\abs{\nabla n_{\eps_\ell}(s,x)}^2dxds,
		\end{split}
	\end{equation*}
it follows from \re{5.2**}, \re{5.2*} and the dominated convergence theorem that
	$
	B_2(t,\eta)\xrightarrow{\eta{\downarrow} 0}0
	$
	and
	\begin{equation*}
		\be\int_0^{t}\int_\bo n_{\eps_\ell}(s,x)\Theta_\eta(n_{\eps_\ell}(s,x))\Theta_\eta'(n_{\eps_\ell}(s,x))dxds\xrightarrow{\eta{\downarrow} 0}\be\int_0^{t}\norm{n_{\eps_\ell}^-(s)}^2_{L^2}ds.
	\end{equation*}
	Combining these calculations, we can pass to the limit in \re{equa4.26} when $\eta$ goes to $0$ and obtain
	\begin{equation*}
		\be\norm{n_{\eps_\ell}^-(t)}^2_{L^2}\leq 	{2|\alpha|}\int_0^{t}\be\norm{n_{\eps_\ell}^-(s)}^2_{L^2}ds,
	\end{equation*}
	and the Gronwall inequality yields $\be\norm{n_{\eps_\ell}^-(t)}^2_{L^2}=0$. This implies that $\norm{n_{\eps_\ell}^-(t)}^2_{L^2}=0, \  {\mathbb{P}}\text{-a.s}$, and \re{equa4.19}$_1$ follows. It remains now to prove \re{equa4.19}$_2$. It is very similar and follows the same steps as the proof of \re{equa4.19}$_1$, but the stochastic impact is treated differently. In fact the It\^o formula from \cite[Theorem 4.32]{DaPrato} implies
	\begin{equation*}
		\begin{split}
			\be\Psi_\eta(c_{\eps_\ell} (t))&=\be\Psi_\eta(v_0)-2\be\Int_0^t\left(\abs{\nabla c_{\eps_\ell}(s)}^2,\Theta_\eta(c_{\eps_\ell}(s))\Theta_\eta^{''}(c_{\eps_\ell}(s))\right)ds\\
			&\quad-2\be\Int_0^t\left(\abs{ \nabla c_{\eps_\ell}(s)}^2,(\Theta_\eta^{'}(c_{\eps_\ell}(s)))^2\right)ds-2\be\int_0^{t}(c_{\eps_\ell}(s),\Theta_\eta(c_{\eps_\ell}(s))\Theta_\eta'(c_{\eps_\ell}(s)))ds\\
			&\quad+2\be\int_0^{t}(H_{\eps_\ell}(  n^+_{\eps_\ell}(s)),\Theta_\eta(c_{\eps_\ell}(s))\Theta_\eta'(c_{\eps_\ell}(s)))ds\\
			& \quad+\be\int_0^t\sum_{k=1}^\infty\left(g_2(c_{\eps_\ell}(s))e_k,\Theta_\eta(c_{\eps_\ell}(s))\Theta_\eta^{''}(c_{\eps_\ell}(s))g_2(c_{\eps_\ell}(s))e_k\right)ds\\
			&\quad+\be\int_0^t\sum_{k=1}^\infty\left(g_2(c_{\eps_\ell}(s))e_k,(\Theta_\eta^{'}(c_{\eps_\ell}(s)))^2g_2(c_{\eps_\ell}(s))e_k\right)ds\\
			&:=D_1(\eta)+D_2(t,\eta)+D_3(t,\eta)+D_4(t,\eta)+D_5(t,\eta)+D_6(t,\eta)+D_7(t,\eta).
		\end{split}
	\end{equation*}
	Let $k\geq 1$ be an integer. By using \re{equa2.6}, we get
	\begin{equation*}
		\begin{split}
			\left(g_2(c_{\eps_\ell})e_k,(\Theta_\eta^{'}(c_{\eps_\ell}))^2g_2(c_{\eps_\ell})e_k\right)&\leq \int_\bo [g_2(c_{\eps_\ell}(x))e_k(x)]^2(\Theta_\eta^{'}(c_{\eps_\ell}(x)))^2dx\\
			&\leq\abs{z_k}^2 \int_\bo \abs{c_{\eps_\ell}(x)}^2(\Theta_\eta^{'}(c_{\eps_\ell}(x)))^2dx.
		\end{split}
	\end{equation*}
	This implies that
	\begin{equation*}
		D_7(t,\eta)\leq \sum_{k=1}^\infty\abs{z_k}^2\be\int_0^t\int_\bo \abs{c_{\eps_\ell}(x,s)}^2(\Theta_\eta^{'}(c_{\eps_\ell}(x,s)))^2dxds.
	\end{equation*}
	Following the same procedure, we can also prove that
	\begin{equation*}
		D_6(t,\eta)\leq \sum_{k=1}^\infty\abs{z_k}^2\be\int_0^t\int_\bo \abs{c_{\eps_\ell}(x,s)}^2\abs{\Theta_\eta(c_{\eps_\ell}(x,s))\Theta_\eta^{''}(c_{\eps_\ell}(x,s))}dxds.
	\end{equation*}
	Since 
	{$D_3(t,\eta)\leq 0$}, we infer from the above inequalities that
	\begin{equation*}
		\begin{split}
			\be\Psi_\eta(c_{\eps_\ell} (t))&\leq \be\Psi_\eta(v_0)-2\be\Int_0^t\left(\abs{\nabla c_{\eps_\ell}(s)}^2,\Theta_\eta(c_{\eps_\ell}(s))\Theta_\eta^{''}(c_{\eps_\ell}(s))\right)ds\\
			&\quad-2\be\int_0^{t}(c_{\eps_\ell}(s),\Theta_\eta(c_{\eps_\ell}(s))\Theta_\eta'(c_{\eps_\ell}(s)))ds+2\be\int_0^{t}(H_{\eps_\ell}(  n^+_{\eps_\ell}(s)),\Theta_\eta(c_{\eps_\ell}(s))\Theta_\eta'(c_{\eps_\ell}(s)))ds\\
			&\quad+\bk\be\int_0^t\int_\bo \abs{c_{\eps_\ell}(x,s)}^2(\Theta_\eta^{'}(c_{\eps_\ell}(x,s)))^2dxds\\
			&\quad+\bk\be\int_0^t\int_\bo \abs{c_{\eps_\ell}(x,s)}^2\abs{\Theta_\eta(c_{\eps_\ell}(x,s))\Theta_\eta^{''}(c_{\eps_\ell}(x,s))}dxds.
		\end{split}
	\end{equation*}
	By using \re{5.2**}, \re{5.2*} and the dominated convergence theorem, we can pass to the limit in this inequality when $\eta$ goes to $0$ and obtain
	\begin{equation*}
		\begin{split}
			\be\norm{c_{\eps_\ell}^-(t)}^2_{L^2}
			&\leq 2\be\int_0^{t}(c_{\eps_\ell}(s),c_{\eps_\ell}^-(s))ds+	\bk\int_0^{t}\be\norm{c_{\eps_\ell}^-(s)}^2_{L^2}ds-2\be\int_0^{t}(H_{\eps_\ell}(  n^+_{\eps_\ell}(s)),c^-_{\eps_\ell}(s))ds\\
			&= 2\be\int_0^{t}(c^+_{\eps_\ell}(s),c_{\eps_\ell}^-(s))ds-2\be\int_0^{t}(c^-_{\eps_\ell}(s),c_{\eps_\ell}^-(s))ds+	\bk\int_0^{t}\be\norm{c_{\eps_\ell}^-(s)}^2_{L^2}ds\\
			&\qquad -2\be\int_0^{t}(H_{\eps_\ell}(  n^+_{\eps_\ell}(s)),c^-_{\eps_\ell}(s))ds.
		\end{split}
	\end{equation*}
	Since $(H_{\eps_\ell}(  n^+_{\eps_\ell}),c^-_{\eps_\ell})\geq 0$ and $(c^+_{\eps_\ell},c_{\eps_\ell}^-)=0$, we deduce that
	\begin{equation*}
		\be\norm{c_{\eps_\ell}^-(t)}^2_{L^2}+2\be\int_0^{t}(H_{\eps_\ell}(  n^+_{\eps_\ell}(s)),c^-_{\eps_\ell}(s))ds+2\be\int_0^{t}(c^-_{\eps_\ell}(s),c_{\eps_\ell}^-(s))ds\leq 	\bk\int_0^{t}\be\norm{c_{\eps_\ell}^-(s)}^2_{L^2}ds.
	\end{equation*}
	By the Gronwall lemma, we get $\be\norm{c_{\eps_\ell}^-(t)}^2_{L^2}=0$ and \re{equa4.19}$_2$ follows, completing the proof.
\end{proof}

\begin{remark}\label{remark4.2}
	Owing to Lemma \ref{lem4.3}, it holds that for any $t\in [0,T]$, $n_{\eps_\ell}^+(t)=n_{\eps_\ell}(t)$ and $c_{\eps_\ell}^+(t)=c_{\eps_\ell}(t)$ almost surely.
\end{remark}

\subsection{Uniform a priori estimates and tightness}
\label{subsec:apriori}

We now proceed to prove uniform bounds for the family of global solutions $\{(n_{\eps_\ell},c_{\eps_\ell})\}_{\ell\geq 1}$ to the regularised  system \re{equa3.1}.

\begin{lemma}\label{lem4.3*}
	For any $p\in [2,32]$, there exists a constant $\bk>0$, {depending only on $T$, $\chi$, $\alpha$, $\mu$, $L_0$, $L_2$, $L_3$, $\mathcal O$, and the moments of $u_0$ and $v_0$}, such that for any $\ell\geq 1$,
	\begin{equation}\label{equa4.5}
		\begin{split}
			&\be\sup_{s\in[0,T]}\norm{n_{\eps_\ell}(s)}_{L^1}^\frac{p}{2}+\mu^{{\frac p2}} \be\left(\int_0^T\norm{H_{\eps_\ell}(n_{\eps_\ell}(s))}^2_{L^2}ds\right)^\frac{p}{2}\leq \left(\be\norm{u_0}^\frac{p}{2}_{L^1}+\bk L_0^\frac{p}{2}\abs{\bo}^{{\frac p4}}\right)e^{\bk T},\\
			&\be	\sup_{s\in[0,T]}\norm{ c_{\eps_\ell}(s)}_{L^2}^p+\be\left(\int_0^T\norm{ c_{\eps_\ell}(s)}_{W^{1,2}}^2ds\right)^\frac{p}{2}\leq \left(\bk\be\norm{v_0}^p_{L^{2}}+\bk\right)e^{\bk T},\\
			&\be\sup_{s\in[0,T]}\norm{\nabla c_{\eps_\ell}(s)}_{L^2}^p+\be\left(\int_0^T\norm{c_{\eps_\ell}(s)}^2_{W^{2,2}}ds\right)^\frac{p}{2}\leq \bk\be\norm{v_0}^p_{W^{1,2}}+\bk\\
			&\be\left(\Int_0^T \norm{\nabla\sqrt{n_{\eps_\ell}(s)+1}}^2_{L^2}ds\right)^4\leq \left(\be\norm{(u_0+1)\ln(u_0+1)}^4_{L^1}+\bk\right)e^{\bk T}\\
			&\be\left(\int_0^T\norm{n_{\eps_\ell}(s)}_{L^2}^2ds\right)^2\leq \bk \left(\be\norm{u_0}^4_{L^1}+\be\norm{(u_0+1)\ln(u_0+1)}^4_{L^1}+1\right)e^{\bk T}.
		\end{split}
	\end{equation}
	In addition, there exists a deterministic constant $\bk>0$ such that for any $\ell\geq 1$,
	\begin{equation}\label{equa4.6*}
		\begin{split}
			&\sup_{r\in[0,T]}\norm{n_{\eps_\ell}(r)}_{L^2}^2+\int_0^{T}\norm{\nabla n_{\eps_\ell}(s)}^2_{L^2}ds\\
			&\leq \left(\norm{u_0}_{L^2}^2+TL_0^2\right)\exp\left(\bk\Int_0^T\left(\norm{\nabla c_{\eps_\ell}(s)}^4_{L^4}+1\right)ds\right)\\
			&	\quad+2\sup_{r\in[0,T]}\abs{\int_0^{r}(g_1(n_{\eps_\ell}(s)),n_{\eps_\ell}(s))dW_1(s)}\exp\left(\bk\Int_0^T\left(\norm{\nabla c_{\eps_\ell}(s)}^4_{L^4}+1\right)ds\right)\quad \mathbb{P}\text{-a.s.}.
		\end{split}
	\end{equation}
\end{lemma}

\begin{proof}
	Let $p\in [2,32]$ be arbitrary fixed. Let us also fix arbitrary $t\in [0,T]$. By Proposition \ref{propo3.5}, we can take $\varphi:=1$ in \re{equa3.1}$_1$ and use Remark \ref{remark4.2}, to get for any $r\in [0,t]$,
	\begin{equation*}
		\norm{n_{\eps_\ell}(r)}_{L^1}=\norm{u_0}_{L^1}+\alpha\int_0^r\norm{n_{\eps_\ell}(s)}_{L^1}ds-\mu \int_0^{{r}}\norm{H_{\eps_\ell}(n_{\eps_\ell}(s))}^2_{L^2}ds+\int_0^r(g_1(n_{\eps_\ell}(s)),1)dW_1(s),\  \mathbb{P}\text{-a.s.},
	\end{equation*}
	This implies that
	\begin{equation*}
		\begin{split}
			&\be\sup_{r\in[0,t]}\norm{n_{\eps_\ell}({r})}_{L^1}^\frac{p}{2}+\mu^{{\frac p2}} \be\left(\int_0^t\norm{H_{\eps_\ell}(n_{\eps_\ell}(s))}^2_{L^2}ds\right)^\frac{p}{2}\\
			&\leq \be\norm{u_0}^\frac{p}{2}_{L^1}+\bk{|\alpha|^{\frac p2}}\be\left(\int_0^t\norm{n_{\eps_\ell}(s)}_{L^1}ds\right)^\frac{p}{2}+\bk\be \sup_{r\in [0,t]}\abs{\int_0^r(g_1(n_{\eps_\ell}(s)),1)dW_1(s)}^\frac{p}{2}.
		\end{split}
	\end{equation*}
	Thanks to the H\"older inequality, we get for $p>2$,
	\begin{equation*}
		\be\left(\int_0^t\norm{n_{\eps_\ell}(s)}_{L^1}ds\right)^\frac{p}{2}\leq \bk t^{{\frac{p-2}{2}}}\be\int_0^t\norm{n_{\eps_\ell}(s)}^\frac{p}{2}_{L^1}ds\leq  \bk T^{{\frac{p-2}{2}}}\be\int_0^t\norm{n_{\eps_\ell}(s)}^\frac{p}{2}_{L^1}ds.
	\end{equation*}
	This implies that for $p\in [2,32]$,
	\begin{equation*}
		\be\left(\int_0^t\norm{n_{\eps_\ell}(s)}_{L^1}ds\right)^\frac{p}{2} \leq  \bk  \be\int_0^t\norm{n_{\eps_\ell}(s)}^\frac{p}{2}_{L^1}ds.
	\end{equation*}
	The Burkholder–Davis–Gundy inequality and \re{equa2.3*} give
	\begin{equation*}
		\be \sup_{r\in [0,t]}\abs{\int_0^r(g_1(n_{\eps_\ell}(s)),1)dW_1(s)}^\frac{p}{2}\leq\bk \be \left(\int_0^t\abs{\bo}\norm{g_1(n_{\eps_\ell}(s))}_{\mathcal{L}_2^0}^2ds\right)^\frac{p}{4}\leq \bk L_0^\frac{p}{2}\abs{\bo}^{{\frac p4}}T^\frac{p}{4}.
	\end{equation*}
	It then follows from the above inequalities that
	\begin{equation*}
		\be\sup_{s\in[0,t]}\norm{n_{\eps_\ell}(s)}_{L^1}^\frac{p}{2}+\mu^{{\frac p2}} \be\left(\int_0^t\norm{H_{\eps_\ell}(n_{\eps_\ell}(s))}^2_{L^2}ds\right)^\frac{p}{2}\leq \be\norm{u_0}^\frac{p}{2}_{L^1}+\bk L_0^{{\frac p2}}\abs{\bo}^{{\frac p4}}+\bk\int_0^t\be\norm{n_{\eps_\ell}(s)}^\frac{p}{2}_{L^1}ds,
	\end{equation*}
	for all $t\in [0,T]$, and the Gronwall inequality yields \re{equa4.5}$_1$. Next, applying the It\^o formula with $ \norm{c_{\eps_\ell}}_{L^2}^2$ and using Remark \ref{remark4.2} and an integration-by-parts (note that the boundary integral will vanish because of $c_\eps\in D(A_2)$), it follows that
	\begin{equation*}
		\begin{split}
			\norm{ c_{\eps_\ell}(r)}_{L^2}^2+2\int_0^r\norm{ c_{\eps_\ell}({s})}_{W^{1,2}}^2ds
			&=\norm{v_0}^2_{L^{2}}+2\int_0^r(H_{\eps_\ell}(n_{\eps_\ell}(s)),c_{\eps_\ell}(s))ds\\
			&\qquad +\int_0^r\norm{g_2(c_{\eps_\ell}(s))}^2_{\mathcal{L}_2^0}ds+2\int_0^r(g_2(c_{\eps_\ell}(s)),c_{\eps_\ell}(s))dW_2(s).
		\end{split}
	\end{equation*}
	By using \re{equa2.4}, we get
	\begin{equation*}
		2(H_{\eps_\ell}(n_{\eps_\ell}),c_{\eps_\ell})\leq2 \norm{H_{\eps_\ell}(n_{\eps_\ell})}_{L^2}\norm{c_{\eps_\ell}}_{L^2}\leq 2\norm{H_{\eps_\ell}(n_{\eps_\ell})}_{L^2}^2+\frac{1}{2}\norm{c_{\eps_\ell}}_{W^{1,2}}^2,
	\end{equation*}
	and
	\begin{equation*}
		\norm{g_2(c_{\eps_\ell})}^2_{\mathcal{L}_2^0}\leq 2\norm{g_2(c_{\eps_\ell})-g_2(0)}^2_{\mathcal{L}_2^0}+2\norm{g_2(0)}^2_{\mathcal{L}_2^0}\leq 2L_2^2\norm{c_{\eps_\ell}}_{L^2}^2+2\norm{g_2(0)}^2_{\mathcal{L}_2^0}.
	\end{equation*}
	By combining these inequalities, we infer that
	\begin{equation*}
		\begin{split}
			&	\be	\sup_{r\in[0,t]}\norm{ c_{\eps_\ell}(r)}_{L^2}^p+\be\left(\int_0^t\norm{ c_{\eps_\ell}(s)}_{W^{1,2}}^2ds\right)^\frac{p}{2}\\
			&\leq \bk\be\norm{v_0}^p_{{L^2}}+\bk\be\left(\int_0^t\norm{H_{\eps_\ell}(n_{\eps_\ell}(s))}^2_{L^2}ds\right)^\frac{p}{2}\\
			&\qquad +\bk\be\int_0^t\norm{c_{\eps_\ell}(s)}^p_{L^2}ds+\bk\be\sup_{r\in[0,t]}\abs{\int_0^r(g_2(c_{\eps_\ell}(s)),c_{\eps_\ell}(s))dW_2(s)}^\frac{p}{2}.
		\end{split}
	\end{equation*}
	Here, we have used the H\"older inequality to have $\Int_0^t\norm{c_{\eps_\ell}(s)}^2_{L^2}\leq T^{{\frac{p-2}{p}}}\left(\Int_0^t\norm{c_{\eps_\ell}(s)}^p_{L^2}{ds}\right)^\frac{2}{p}$, for $p>2$. Thanks to \re{equa2.4} and the H\"older inequality once more, we get
	\begin{equation*}
		\begin{split}
			&\bk\be\sup_{r\in[0,t]}\abs{\int_0^r(g_2(c_{\eps_\ell}(s)),c_{\eps_\ell}(s))dW_2(s)}^\frac{p}{2}\\
			&\leq \bk\be\left(\int_0^t\norm{g_2(c_{\eps_\ell}(s))-g_2(0)}^2_{\mathcal{L}^0_2}\norm{c_{\eps_\ell}(s)}_{L^2}^2ds\right)^\frac{p}{4}+ \bk\norm{g_2(0)}^\frac{p}{2}_{\mathcal{L}^0_2}\be\left(\int_0^t\norm{c_{\eps_\ell}(s)}_{L^2}^2ds\right)^\frac{p}{4}\\
			&\leq \frac{1}{2}\be\sup_{s\in[0,t]}\norm{c_{\eps_\ell}(s)}_{L^2}^p+\bk\be\left(\int_0^t\norm{c_{\eps_\ell}(s)}_{L^2}^2ds\right)^\frac{p}{2}+ \bk T^{{\frac{p-2}{4}}}\norm{g_2(0)}^\frac{p}{2}_{\mathcal{L}^0_2}\left(\be\int_0^t\norm{c_{\eps_\ell}(s)}_{L^2}^pds\right)^\frac{1}{2}\\
			&\leq \frac{1}{2}\be\sup_{s\in[0,t]}\norm{c_{\eps_\ell}(s)}_{L^2}^p+\bk\be\int_0^t\norm{c_{\eps_\ell}(s)}_{L^2}^pds+\bk.
		\end{split}
	\end{equation*}
	Note that this inequality holds also with $p=2$ but with different constants $\bk$. Combining these two inequalities and using \re{equa4.5}$_1$, we get
	\begin{equation*}
		\begin{split}
			&	\frac{1}{2}\be	\sup_{s\in[0,t]}\norm{ c_{\eps_\ell}(s)}_{L^2}^p+\be\left(\int_0^t\norm{ c_{\eps_\ell}(s)}_{W^{1,2}}^2ds\right)^\frac{p}{2}\\
			&\leq \bk\be\norm{v_0}^p_{L^{2}}+{\bk\mu^{-\frac p2}}\left(\be\norm{u_0}^{{\frac p2}}_{L^1}+\bk L_0^{{\frac p2}}\abs{\bo}^{{\frac p4}}\right)e^{\bk T}+\bk +\bk\int_0^t\be\norm{c_{\eps_\ell}(s)}^p_{L^2}ds, \quad \forall t\in [0,T],
		\end{split}
	\end{equation*}
	and by applying the Gronwall lemma, we get \re{equa4.5}$_2$. Let us now apply the generalised version of the It\^o formula from \cite[Proposition 3.3]{Scarpa} with  $\norm{\nabla c_{\eps_\ell}}_{L^2}^2$ to get
	\begin{equation*}
		\begin{split}
			&\norm{ \nabla c_{\eps_\ell}(t)}_{L^2}^2+2\int_0^t\left(\norm{\Delta c_{\eps_\ell}(s)}^2_{L^2}+\norm{\nabla c_{\eps_\ell}(s)}^2_{L^2}\right)ds\\
			&=\norm{\nabla v_0}^2_{L^2}+2\int_0^t(H_{\eps_\ell}(n_{\eps_\ell}(s)),-\Delta c_{\eps_\ell}(s))ds\\
			&\qquad +\int_0^t\norm{\nabla g_2(c_{\eps_\ell}(s))}^2_{\mathcal{L}_2^0}ds+2\int_0^t(g_2(c_{\eps_\ell}(s)),-\Delta c_{\eps_\ell}(s))dW_2(s).
		\end{split}
	\end{equation*}
	Since
	\begin{equation*}
		2(H_{\eps_\ell}(n_{\eps_\ell}),-\Delta c_{\eps_\ell})\leq2 \norm{H_{\eps_\ell}(n_{\eps_\ell})}_{L^2}\norm{\Delta c_{\eps_\ell}}_{L^2}\leq 2\norm{H_{\eps_\ell}(n_{\eps_\ell})}_{L^2}^2+\frac{1}{2}\norm{\Delta c_{\eps_\ell}}_{L^2}^2,
	\end{equation*}
	and
	\begin{equation*}
		\norm{\nabla g_2(c_{\eps_\ell})}^2_{\mathcal{L}_2^0}\leq 2\norm{g_2(c_{\eps_\ell})-g_2(0)}^2_{\mathcal{L}_2^1}+2\norm{g_2(0)}^2_{\mathcal{L}_2^1}\leq 2L_2^2\norm{c_{\eps_\ell}}_{W^{1,2}}^2+2\norm{g_2(0)}^2_{\mathcal{L}_2^1},
	\end{equation*}
	we derive by using \re{equa4.5}$_1$ and \re{equa4.5}$_2$ that
	\begin{equation*}
		\begin{split}
			&\be\sup_{s\in[0,T]}\norm{\nabla c_{\eps_\ell}(s)}_{L^2}^p+\be\left(\int_0^T\left(\norm{\Delta c_{\eps_\ell}(s)}^2_{L^2}+2\norm{\nabla c_{\eps_\ell}(s)}_{L^2}^2\right)ds\right)^\frac{p}{2}\\
			&\leq\bk\be\norm{v_0}^p_{W^{1,2}}+\bk\norm{g_2(0)}^p_{\mathcal{L}_2^1}+ \bk\be\left(\int_0^T\norm{H_{\eps_\ell}(n_{\eps_\ell}(s))}_{L^2}^2ds\right)^\frac{p}{2}+\bk\be\left(\int_0^T\norm{c_{\eps_\ell}(s)}_{W^{1,2}}^2ds\right)^\frac{p}{2}\\
			&\qquad +\bk\be\sup_{t\in[0,T]}\abs{\int_0^{{t}}(g_2(c_{\eps_\ell}(s)),-\Delta c_{\eps_\ell}(s))dW_2(s)}^\frac{p}{2}\\
			&\leq \bk\be\norm{v_0}^p_{W^{1,2}}+\bk +\bk\be\left(\int_0^T\norm{g_2(c_{\eps_\ell}(s))}^2_{\mathcal{L}^0_2}\norm{\Delta c_{\eps_\ell}(s)}_{L^2}^2ds\right)^\frac{p}{4}\\
			&\leq \bk\be\norm{v_0}^p_{W^{1,2}}+\bk +\bk\be\sup_{s\in[0,T]}\norm{c_{\eps_\ell}(s)}_{L^2}^\frac{p}{2}\left(\int_0^T\norm{\Delta c_{\eps_\ell}(s)}_{L^2}^2ds\right)^\frac{p}{4}+\bk\norm{g_2(0)}^\frac{p}{2}_{\mathcal{L}^0_2}\be\left(\int_0^T\norm{\Delta c_{\eps_\ell}(s)}_{L^2}^2ds\right)^\frac{p}{4}\\
			&\leq \bk\be\norm{v_0}^p_{W^{1,2}}+\bk +\bk\be\sup_{s\in[0,T]}\norm{c_{\eps_\ell}(s)}_{L^2}^p+\frac{1}{2}\be\left(\int_0^T\norm{\Delta c_{\eps_\ell}(s)}_{L^2}^2ds\right)^\frac{p}{2}\\
			&\leq \bk\be\norm{v_0}^p_{W^{1,2}}+\bk+\frac{1}{2}\be\left(\int_0^T\norm{\Delta c_{\eps_\ell}(s)}_{L^2}^2ds\right)^\frac{p}{2}.
		\end{split}
	\end{equation*}
	This implies that
	\begin{equation*}
		\be\sup_{s\in[0,T]}\norm{\nabla c_{\eps_\ell}(s)}_{L^2}^p+\frac{1}{2}\be\left(\int_0^T\norm{\Delta c_{\eps_\ell}(s)}^2_{L^2}ds\right)^\frac{p}{2}+\be\left(\int_0^T\norm{\nabla c_{\eps_\ell}(s)}^2_{L^2}ds\right)^\frac{p}{2}\leq \bk\be\norm{v_0}^p_{W^{1,2}}+\bk,
	\end{equation*}
	and \re{equa4.5}$_3$ follows {(using the elliptic estimate $\norm{\psi}_{W^{2,2}}\le\bk(\norm{\Delta \psi}_{L^2}+\norm{\psi}_{L^2})$ for $\psi\in D(A_2)$)}. Let us now consider the following $C^2(\br)$ function $\Psi$ defined by
	\begin{equation*}
		\Psi(x):=\begin{cases}
			\frac{1}{2}x^2+x,\quad\text{ if } x<0,\\
			(x+1)\ln(x+1),\quad\text{ if } x\geq 0,
		\end{cases}
	\end{equation*}
	and define the following functional $\Phi: L^2(\bo)\to \br$ as
	$
	\Phi(u):=\Int_\bo\Psi(u(x))dx, \quad u\in L^2(\bo).
	$
	The functional $\Phi$ is twice Fr{\'e}chet-differentiable and for any $u, h\in L^2(\bo),$
	\begin{equation*}
		\begin{split}
			&\Phi'(u):=\Psi'( u), \quad \Phi^{''}(u)(h):= \Psi^{''}(u)h, \quad \text{ in } L^2(\bo).
		\end{split}
	\end{equation*}
	Let us fix $R>0$ and $u\in L^2(\bo)$ such that $\norm{u}_{L^2}\leq R$.  Since $\ln(x+1)\leq x$, for any $x\geq 0$, it {holds}
	\begin{equation*}
		\begin{split}
			\abs{\Phi(u)}&\leq 
			\frac{1}{2}\int_{\{u(x)<0\}} u^2(x)dx+\int_{\{u(x)<0\}} \abs{u(x)}dx+\int_{\{u(x)\geq0\}} (u(x)+1)\ln (u(x)+1)dx\\
			&\leq \frac{1}{2}\norm{u}^2_{L^2}+\norm{u}_{L^1}+\int_{\{u(x)\geq0\}} (u(x)+1) u(x) dx\leq \frac{3}{2}\norm{u}^2_{L^2}+2\norm{u}_{L^1}\leq \frac{3}{2}R^2+2\abs{\bo}^\frac{1}{2}R.
		\end{split}
	\end{equation*}
	This implies that $\Phi$ is locally bounded. In a very similar way, we can show that $\Phi'$ and $\Phi^{''}$ are locally bounded. Moreover, for $u\in W^{1,2}(\bo)$, we have
	\begin{equation*}
		\begin{split}
			\norm{\Phi'(u)}^2_{W^{1,2}}
			&=\norm{\Psi'(u)}^2_{L^2}+\norm{\nabla\Psi'(u)}^2_{L^2}\\
			&=\int_{\{u(x)<0\}} (u(x)+1)^2dx+\int_{\{u(x)\geq0\}} (\ln(u(x)+1)+1)^2dx\\
			&\qquad+\int_{\{u(x)<0\}} \abs{\nabla u(x)}^2dx+\int_{\{u(x)\geq0\}}\frac{\abs{\nabla u(x)}^2}{(u(x)+1)^2}dx\\
			&\leq \int_{\bo} (u(x)+1)^2dx+\int_{\{u(x)\geq0\}} (u(x)+1)^2dx+\norm{\nabla u}^2_{L^2}\leq 4\norm{u}^2_{L^2}+4\abs{\bo}+\norm{\nabla u}^2_{L^2}.
		\end{split}
	\end{equation*}
	This implies that
	$
	\norm{\Phi'(u)}_{W^{1,2}}\leq \bk (1+\norm{u}_{W^{1,2}}).
	$
Using   this, and \re{equa2.3*}, 
we can	apply the It\^o formula in the version of \cite[Theorem 1.2]{Pardoux} with the functional $\Phi$ to get for any $r\in [0,t]$,
	\begin{equation}\label{equa4.6}
		\begin{split}
			&\norm{(n_{\eps_\ell}(r)+1)\ln(n_{\eps_\ell}(r)+1)}_{L^1}+4\Int_0^r \norm{\nabla\sqrt{n_{\eps_\ell}(s)+1}}^2_{L^2}ds\\
			&=\norm{(u_0+1)\ln(u_0+1)}_{L^1}+\chi\Int_0^r (F_{\eps_\ell}( n_{\eps_\ell}(s))\nabla  c_{\eps_\ell}(s),\nabla\ln(n_{\eps_\ell}(s)+1))ds\\
			&	\hspace{0.5cm}+\alpha\int_0^{r}(n_{\eps_\ell}(s),\ln(n_{\eps_\ell}(s)+1)+1)ds-\mu\int_0^{r} ([H_{\eps_\ell}(n_{\eps_{\ell}}(s))]^2,\ln(n_{\eps_\ell}(s)+1)+1)ds\\
			&	\hspace{0.5cm}+{\frac12}\int_0^{r}\sum_{k=1}^\infty\norm{\frac{g_1(n_{\eps_\ell}(s))e_k}{\sqrt{n_{\eps_\ell}+1}}}^2_{L^2}ds+\int_0^{r}(g_1(n_{\eps_\ell}(s)),\ln(n_{\eps_\ell}(s)+1)+1)dW_1(s).
		\end{split}
	\end{equation}
	It follows from the definition of $F_{\eps_\ell}$ that
	\begin{equation*}
		\begin{split}
			&\chi(F_{\eps_\ell}( n_{\eps_\ell})\nabla  c_{\eps_\ell},\nabla\ln(n_{\eps_\ell}+1))\\
			&=\chi\int_\bo\frac{n_{\eps_\ell}(x)\nabla  c_{\eps_\ell}(x)\cdot\nabla  n_{\eps_\ell}(x)}{(1+\eps_\ell n_{\eps_\ell}(x))(n_{\eps_\ell}(x)+1)}dx=2\chi\int_\bo\frac{n_{\eps_\ell}(x)\nabla  c_{\eps_\ell}(x)\cdot\nabla  \sqrt{n_{\eps_\ell}(x)+1}}{(1+\eps_\ell n_{\eps_\ell}(x))\sqrt{n_{\eps_\ell}(x)+1}}dx\\
			&=2\chi\int_\bo\frac{(n_{\eps_\ell}(x)+1)\nabla  c_{\eps_\ell}(x)\cdot\nabla  \sqrt{n_{\eps_\ell}(x)+1}}{(1+\eps_\ell n_{\eps_\ell}(x))\sqrt{n_{\eps_\ell}(x)+1}}dx-2\chi\int_\bo\frac{\nabla  c_{\eps_\ell}(x)\cdot\nabla  \sqrt{n_{\eps_\ell}(x)+1}}{(1+\eps_\ell n_{\eps_\ell}(x))\sqrt{n_{\eps_\ell}(x)+1}}dx\\
			&=2\chi\int_\bo\frac{\sqrt{n_{\eps_\ell}(x)+1}\nabla  c_{\eps_\ell}(x)\cdot\nabla  \sqrt{n_{\eps_\ell}(x)+1}}{1+\eps_\ell n_{\eps_\ell}(x)}dx-2\chi\int_\bo\frac{\nabla  c_{\eps_\ell}(x)\cdot\nabla  \sqrt{n_{\eps_\ell}(x)+1}}{(1+\eps_\ell n_{\eps_\ell}(x))\sqrt{n_{\eps_\ell}(x)+1}}dx.
		\end{split}
	\end{equation*}
	Since $\frac{1}{1+\eps_\ell n_{\eps_\ell}}\leq 1$ and $\frac{1}{(1+\eps_\ell n_{\eps_\ell})\sqrt{n_{\eps_\ell}+1}}\leq 1$, the inequality \re{equa2.3} and the Young inequality yield
	\begin{equation*}
		\begin{split}
			&\chi(F_{\eps_\ell}( n_{\eps_\ell})\nabla  c_{\eps_\ell},\nabla\ln(n_{\eps_\ell}+1))\\
			&\leq2\chi\norm{\sqrt{n_{\eps_\ell}+1}}_{L^4}\norm{\nabla  c_{\eps_\ell}}_{L^4}\norm{\nabla  \sqrt{n_{\eps_\ell}+1}}_{L^2}+2\chi\norm{\nabla  c_{\eps_\ell}}_{L^2}\norm{\nabla  \sqrt{n_{\eps_\ell}+1}}_{L^2}\\
			&\leq\bk\norm{\sqrt{n_{\eps_\ell}+1}}_{L^2}^\frac{1}{2}\norm{\nabla  c_{\eps_\ell}}_{L^4}\norm{\nabla  \sqrt{n_{\eps_\ell}+1}}^{{\frac32}}_{L^2}+\bk\norm{\sqrt{n_{\eps_\ell}+1}}_{L^2}\norm{\nabla  c_{\eps_\ell}}_{L^4}\norm{\nabla  \sqrt{n_{\eps_\ell}+1}}_{L^2}\\
			&\qquad+2\chi\norm{\nabla  c_{\eps_\ell}}_{L^2}\norm{\nabla  \sqrt{n_{\eps_\ell}+1}}_{L^2}\\
			&\leq\norm{\nabla  \sqrt{n_{\eps_\ell}+1}}^2_{L^2}+\bk\norm{\sqrt{n_{\eps_\ell}+1}}_{L^2}^2\norm{\nabla  c_{\eps_\ell}}^4_{L^4}+\bk\norm{\sqrt{n_{\eps_\ell}+1}}^2_{L^2}\norm{\nabla  c_{\eps_\ell}}^2_{L^4}+\bk\norm{\nabla  c_{\eps_\ell}}^2_{L^2}\\
			&\leq\norm{\nabla  \sqrt{n_{\eps_\ell}+1}}^2_{L^2}+\bk\norm{\sqrt{n_{\eps_\ell}+1}}_{L^2}^2(\norm{\nabla  c_{\eps_\ell}}^2_{L^2}\norm{  c_{\eps_\ell}}^2_{W^{2,2}}+1)+\bk\norm{\nabla  c_{\eps_\ell}}^2_{L^2}.
		\end{split}
	\end{equation*}
	Estimating $n_{\eps_\ell}\ln(n_{\eps_\ell}+1)\leq (n_{\eps_\ell}+1)\ln(n_{\eps_\ell}+1)$ and exploiting the non-negativity of $n_{\eps_\ell}$ shows
	\begin{equation*}
		\alpha(n_{\eps_\ell},\ln(n_{\eps_\ell}+1)+1)\leq {|}\alpha{|}\norm{(n_{\eps_\ell}+1)\ln(n_{\eps_\ell}+1)}_{L^1}+{|}\alpha{|}\norm{n_{\eps_\ell}}_{L^1},
	\end{equation*}
	and
	\begin{equation*}
		-([H_{\eps_\ell}(n_{\eps_{\ell}})]^2,\ln(n_{\eps_\ell}+1)+1)\leq 0.
	\end{equation*}
	Thanks to \re{equa2.3*} and the fact that $\frac{1}{n_{\eps_\ell}+1}\leq 1$,
	\begin{equation*}
		\sum_{k=1}^\infty\norm{\frac{g_1(n_{\eps_\ell})e_k}{\sqrt{n_{\eps_\ell}+1}}}^2_{L^2}\leq \norm{g_1(n_{\eps_\ell})}^2_{\mathcal{L}_2^0}\leq L_0^2.
	\end{equation*}
	In view of the Burkholder–Davis–Gundy inequality and \re{equa2.3*}, we obtain
	\begin{equation*}
		\begin{split}
			&	\be \sup_{{r}\in [0,t]}\abs{\int_0^{r}(g_1(n_{\eps_\ell}(s)),\ln(n_{\eps_\ell}(s)+1)+1)dW_1(s)}^4\\
			&\leq\bk \be \left(\int_0^t\norm{\ln(n_{\eps_\ell}(s)+1)+1}_{L^2}^2\norm{g_1(n_{\eps_\ell}(s))}_{\mathcal{L}_2^0}^2ds\right)^2\\
			&\leq \bk L_0^4+\bk L_0^4 \be \left(\int_0^t\norm{\ln(n_{\eps_\ell}(s)+1)}_{L^2}^2ds\right)^2.
		\end{split}
	\end{equation*}
	By using the estimate $\ln(n_{\eps_\ell}+1)\leq n_{\eps_\ell}\leq n_{\eps_\ell}+1$ we derive that $(\ln(n_{\eps_\ell}+1))^2\leq ( n_{\eps_\ell}+1)\ln(n_{\eps_\ell}+1)$ and hence, the H\"older inequality and the Young inequality yield
	\begin{equation*}
		\begin{split}
		&	\be \sup_{{r}\in [0,t]}\abs{\int_0^{r}(g_1(n_{\eps_\ell}(s)),\ln(n_{\eps_\ell}(s)+1)+1)dW_1(s)}^4\\
			&\leq \bk L_0^4+\bk L_0^4 \be \left(\int_0^t\norm{( n_{\eps_\ell}+1)\ln(n_{\eps_\ell}+1)}_{L^1}ds\right)^2\\
			&\leq \bk L_0^4+\bk L_0^4 \be \left(\int_0^t\norm{( n_{\eps_\ell}+1)\ln(n_{\eps_\ell}+1)}^4_{L^1}ds\right)^\frac{1}{2}\\
			&\leq \bk +\bk \be \int_0^t\norm{( n_{\eps_\ell}+1)\ln(n_{\eps_\ell}+1)}^4_{L^1}ds.
		\end{split}
	\end{equation*}
	By taking the supremum over $[0,{t}]$,  {raising to the fourth power} and taking the expectation in \re{equa4.6}, and using the above inequalities, we deduce that
	\begin{equation*}
		\begin{split}
			&\be\sup_{s\in[0,t]}\norm{(n_{\eps_\ell}(s)+1)\ln(n_{\eps_\ell}(s)+1)}^4_{L^1}+3\be\left(\Int_0^t \norm{\nabla\sqrt{n_{\eps_\ell}(s)+1}}^2_{L^2}ds\right)^4\\
			&\leq \be\norm{(u_0+1)\ln(u_0+1)}^4_{L^1}+\bk +\bk\be\left(\int_0^t\norm{\sqrt{n_{\eps_\ell}(s)+1}}_{L^2}^2\left(\norm{\nabla  c_{\eps_\ell}(s)}^2_{L^2}\norm{  c_{\eps_\ell}(s)}^2_{W^{2,2}}+1\right)ds\right)^4\\
			&\quad+\bk\be\left(\int_0^t\norm{\nabla  c_{\eps_\ell}(s)}^2_{L^2}ds\right)^4+\alpha\be\left(\int_0^t\norm{n_{\eps_\ell}(s)}_{L^1}ds\right)^4+\bk \be \int_0^t\norm{( n_{\eps_\ell}+1)\ln(n_{\eps_\ell}+1)}^4_{L^1}ds\\
			&= \be\norm{(u_0+1)\ln(u_0+1)}^4_{L^1}+\bk+J_1+J_2+J_3+J_4.
		\end{split}
	\end{equation*}
	Thanks to \re{equa4.5}$_1$, \re{equa4.5}$_2$, and \re{equa4.5}$_3$, we have
	\begin{equation*}
		\begin{split}
		J_1	&\leq \be\sup_{s\in[0,T]}\norm{\sqrt{n_{\eps_\ell}(s)+1}}_{L^2}^8\sup_{s\in[0,T]} \norm{\nabla  c_{\eps_\ell}(s)}^8_{L^2}\left(\int_0^T\norm{  c_{\eps_\ell}(s)}^2_{W^{2,2}}ds\right)^4+\bk\be\sup_{s\in[0,T]}\norm{\sqrt{n_{\eps_\ell}(s)+1}}_{L^2}^8\\
			&\leq \left[\be\sup_{s\in[0,T]}\norm{\sqrt{n_{\eps_\ell}(s)+1}}_{L^2}^{16}\right]^\frac{1}{2}\left[\be\sup_{s\in[0,T]} \norm{\nabla  c_{\eps_\ell}(s)}^{32}_{L^2}\right]^\frac{1}{4}\left[\be\left(\int_0^T\norm{  c_{\eps_\ell}(s)}^2_{W^{2,2}}ds\right)^{16}\right]^\frac{1}{4}\\
			&\qquad+\bk\be\sup_{s\in[0,T]}\norm{\sqrt{n_{\eps_\ell}(s)+1}}_{L^2}^8\\
			&\leq \bk,
		\end{split}
	\end{equation*}
	and
	\begin{equation*}
J_2+J_3 \leq T^{{4}} \be\sup_{s\in[0,T]} \norm{\nabla  c_{\eps_\ell}(s)}^{8}_{L^2}+{|\alpha|^4}T^{{4}}\be\sup_{s\in[0,T]} \norm{ n_{\eps_\ell}(s)}^{4}_{L^1}\leq \bk.
	\end{equation*}
	Combining these inequalities we obtain
	\begin{equation*}
		\begin{split}
			&\be\sup_{s\in[0,t]}\norm{(n_{\eps_\ell}(s)+1)\ln(n_{\eps_\ell}(s)+1)}^4_{L^1}+3\be\left(\Int_0^t \norm{\nabla\sqrt{n_{\eps_\ell}(s)+1}}^2_{L^2}ds\right)^4\\
			&\leq \be\norm{(u_0+1)\ln(u_0+1)}^4_{L^1}+\bk +\bk \int_0^t \be\norm{( n_{\eps_\ell}+1)\ln(n_{\eps_\ell}+1)}^4_{L^1}ds, \quad \forall t\in [0,T],
		\end{split}
	\end{equation*}
	and applying the Gronwall lemma to obtain \re{equa4.5}$_4$. By using \re{equa2.3}, we get
	\begin{equation*}
		\begin{split}
			\norm{n_{\eps_\ell}}_{L^2}^2&\leq \norm{n_{\eps_\ell}+1}_{L^2}^2=\norm{\sqrt{n_{\eps_\ell}+1}}_{L^4}^4\\
			&\leq {\bk}\norm{\sqrt{n_{\eps_\ell}+1}}_{L^2}^2\norm{\nabla\sqrt{n_{\eps_\ell}+1}}_{L^2}^2+{\bk}\norm{\sqrt{n_{\eps_\ell}+1}}_{L^2}^4\\
			&\leq {\bk}\norm{n_{\eps_\ell}+1}_{L^1}\norm{\nabla\sqrt{n_{\eps_\ell}+1}}_{L^2}^2+{\bk}\norm{n_{\eps_\ell}+1}_{L^1}^2,
		\end{split}
	\end{equation*}
	and therefore, by \re{equa4.5}$_1$ and \re{equa4.5}$_4$,
	\begin{equation*}
		\begin{split}
			\be\left(\int_0^t\norm{n_{\eps_\ell}(s)}_{L^2}^2ds\right)^2&\leq\bk\be\left(\int_0^T\norm{n_{\eps_\ell}(s)+1}_{L^1}\norm{\nabla\sqrt{n_{\eps_\ell}(s)+1}}_{L^2}^2ds\right)^2\\
			&\qquad+\bk\be \left(\int_0^T\norm{n_{\eps_\ell}(s)+1}_{L^1}^2ds\right)^2\\
			&\leq\bk\be \sup_{s\in[0,T]}\norm{n_{\eps_\ell}(s)+1}_{L^1}^2\left(\int_0^T\norm{\nabla\sqrt{n_{\eps_\ell}(s)+1}}_{L^2}^2ds\right)^2\\
			&\qquad+\bk\be \sup_{s\in[0,T]}\norm{n_{\eps_\ell}(s)+1}_{L^1}^4\\
			&\leq\bk\be \sup_{s\in[0,T]}\norm{n_{\eps_\ell}(s)+1}_{L^1}^4+\bk\be\left(\int_0^T\norm{\nabla\sqrt{n_{\eps_\ell}(s)+1}}_{L^2}^2ds\right)^4\leq \bk,
		\end{split}
	\end{equation*}
	and \re{equa4.5}$_5$ follows thanks to \re{equa4.5}$_1$ and \re{equa4.5}$_4$. Let us now apply the It\^o formula with $\norm{n_{\eps_\ell}}^2_{L^2}$ to obtain for any $r\in [0,t]$,
	\begin{equation}\label{equa4.7}
		\begin{split}
			&\norm{n_{\eps_\ell}(r)}_{L^2}^2+2\int_0^{r}\norm{\nabla n_{\eps_\ell}(s)}^2_{L^2}ds\\
			&=\norm{u_0}_{L^2}^2+ 2\chi \int_0^{r}(F_{\eps_\ell}(n_{\eps_\ell}(s))\nabla c_{\eps_\ell}(s),\nabla n_{\eps_\ell}(s) )ds-2\mu\int_0^{r} ([H_{\eps_\ell}(n_{\eps_{\ell}}(s))]^2,n_{\eps_\ell}(s))ds\\
			&	\hspace{0.1cm}+2\alpha\int_0^{r}\norm{n_{\eps_\ell}(s)}^2_{L^2}ds+\int_0^{r}\norm{g_1(n_{\eps_\ell}(s))}^2_{\mathcal{L}_2^0}ds+2\int_0^{r}(g_1(n_{\eps_\ell}(s)),n_{\eps_\ell}(s))dW_1(s).
		\end{split}
	\end{equation}
	The non-negativity of $n_{\eps_\ell}$ and \re{equa2.3} as well as  the fact that $1+\eps_\ell n_{\eps_\ell}{\geq}1$ yield
	\begin{equation*}
		\begin{split}
			2\chi(F_{\eps_\ell}(n_{\eps_\ell}(s))\nabla c_{\eps_\ell}(s),\nabla n_{\eps_\ell}(s) )&=2\chi\left(  \frac{n_{\eps_\ell}}{1+\eps_\ell n_{\eps_\ell}}\nabla  c_{\eps_\ell},\nabla n_{\eps_\ell}\right)\\
			&\leq 2\chi\norm{n_{\eps_\ell}}_{L^4}\norm{\nabla c_{\eps_\ell}}_{L^4}\norm{\nabla n_{\eps_\ell}}_{L^2}\\
			&\leq \frac{1}{2} \norm{\nabla n_{\eps_\ell}}_{L^2}^2+2\chi^2\norm{n_{\eps_\ell}}^2_{L^4}\norm{\nabla c_{\eps_\ell}}^2_{L^4}\\
			&\leq \frac{1}{2} \norm{\nabla n_{\eps_\ell}}_{L^2}^2+\bk\norm{\nabla n_{\eps_\ell}}_{L^2}\norm{n_{\eps_\ell}}_{L^2}\norm{\nabla c_{\eps_\ell}}^2_{L^4}+\bk\norm{n_{\eps_\ell}}^2_{L^2}\norm{\nabla c_{\eps_\ell}}^2_{L^4}\\
			&\leq  \norm{\nabla n_{\eps_\ell}}_{L^2}^2+\bk\norm{n_{\eps_\ell}}^2_{L^2}\norm{\nabla c_{\eps_\ell}}^4_{L^4}+\bk\norm{n_{\eps_\ell}}^2_{L^2},
		\end{split}
	\end{equation*}
	and $$-([H_{\eps_\ell}(n_{\eps_{\ell}})]^2,n_{\eps_\ell})\leq 0.$$ By taking the supremum over $[0,t]$ in \re{equa4.7} and using the above inequalities and \re{equa2.3*}, we derive that
	\begin{equation*}
		\begin{split}
			\sup_{r\in[0,t]}\norm{n_{\eps_\ell}(r)}_{L^2}^2+\int_0^{t}\norm{\nabla n_{\eps_\ell}(s)}^2_{L^2}ds
			&\leq \norm{u_0}_{L^2}^2+TL_0^2+ \bk \int_0^{t}\norm{n_{\eps_\ell}(s)}^2_{L^2}\left(\norm{\nabla c_{\eps_\ell}(s)}^4_{L^4}+1\right)ds\\
			&	\qquad+2\sup_{r\in[0,t]}\abs{\int_0^{r}(g_1(n_{\eps_\ell}(s)),n_{\eps_\ell}(s))dW_1(s)},\quad\forall t\in [0,T],
		\end{split}
	\end{equation*}
	and the Gronwall lemma yields \re{equa4.6*} and completes the proof.
\end{proof}

Let us fix $q>4$ and recall from \cite[Lemma 2.1]{Primio} the existence of a constant $\bk_0>0$ such that
\begin{equation}\label{equa4.9**}
	\norm{u}_{W^{\frac{2}{q},q}(0,T;(W^{3,2}(\bo))')}\leq \bk_0\norm{u}_{W^{1,2}(0,T;(W^{3,2}(\bo))')}\quad \forall u\in {W^{1,2}}(0,T;(W^{3,2}(\bo))').
\end{equation}

As a consequence of Lemma \ref{lem4.3*}, we have the following corollary.

\begin{cor}\label{cor4.6}
	There exists a constant $\bk>0$ such that for any $\ell\geq 1$ and $R> 2$,
	\begin{equation}\label{equa4.9*}
		\begin{split}
			&	\mathbb{P}\left\{\omega\in \Omega: \sup_{0\leq s\leq T}\norm{n_{\eps_\ell}(s)}_{L^2}^2> R\right\}\leq \frac{\bk}{\sqrt{R}}+\frac{\bk}{\ln (\sqrt{R})-\ln \sqrt{2}},\\
			&	\mathbb{P}\left\{\omega\in \Omega: \int_0^T\norm{\nabla n_{\eps_\ell}(s)}_{L^2}^2> R\right\}\leq \frac{\bk}{\sqrt{R}}+\frac{\bk}{\ln (\sqrt{R})-\ln \sqrt{2}}.
		\end{split}
	\end{equation}
	In addition, for any $q>4$, there exists $\bk>0$ such that for any $\ell\geq 1$ and $R>\max(2, 16\bk_0^2\bk_1^2\mu^2(T^2+1)T)$,
	\begin{equation}\label{equa4.10}
		\begin{split}
			&	\mathbb{P}\left\{ \omega\in \Omega: \norm{n_{\eps_\ell}}^2_{W^{\frac{2}{q},q}(0,T;(W^{3,2}(\bo))')}>R\right\}\\
			&\leq \frac{\bk}{R^\frac{1}{4}}+\frac{\bk}{\ln (R^\frac{1}{4})-\ln((4\bk^2_0\bk_1^2\mu^2(T^2+1) T)^\frac{1}{4})-\ln \sqrt{2}},
		\end{split}
	\end{equation}
	where $\bk_0>0$ and $\bk_1>0$ are constants coming respectively from \re{equa4.9**} and from the embedding of $W^{3,2}(\bo)$ into $L^\infty(\bo)$.
	Moreover, let $\xi>0$ and $\eta>0$. There exists $\zeta>0$ such that for every stopping times $\{\tau_\ell\}_{\ell\geq 1}$ satisfying for any $\ell\geq 1$, $0 \leq  \tau_\ell \leq  T$, one has for any $\ell\geq 1$,
	\begin{equation}\label{equa4.9}
		\sup_{0<\theta\leq \zeta}\mathbb{P}\left\{ \omega\in \Omega: \norm{c_{\eps_\ell}(\tau_\ell+\theta){-c_{\eps_\ell}(\tau_\ell)}}_{L^2}\geq \eta\right\}\leq \xi,
	\end{equation}
	where $c_{\eps_\ell}$ is extended by {$c_{\eps_\ell}(t):=c_{\eps_\ell}(T)$ for $t>T$}. 
\end{cor}
\begin{remark}
	Since we have not succeeded in obtaining a uniform estimate of $\be \Sup_{s\in[0,T]}\norm{n_{\eps_\ell}(s)}^2_{L^2}$, it is difficult to obtain a similar estimate like \re{equa4.9} for $n_{\eps_\ell}$ even in the dual space $(W^{3,2}(\bo))'$ due to the presence of the {quadratic term} in the logistic term.
\end{remark}

Let us now give a detailed proof of Corollary \ref{cor4.6}.

\begin{proof}
	Let $q>4$, $\ell\geq 1$, and $R>2$ be fixed. By using the fact that for any random variables $X$ and $Y$,
	\begin{equation*}
		\left\{\omega\in \Omega: X(\omega)\leq \frac{R}{2}\right\}\cap \left\{\omega\in \Omega: Y(\omega)\leq \frac{R}{2}\right\}\subset \left\{\omega\in \Omega: X(\omega)+Y(\omega)\leq R\right\},
	\end{equation*}
	the notion of the complementary, and \re{equa4.6*} as well, we deduce that
	\begin{equation}\label{equa4.11}
		\begin{split}
			&\mathbb{P}\left\{\omega\in \Omega: \sup_{0\leq s\leq T}\norm{n_{\eps_\ell}(s)}_{L^2}^2> R\right\}\\
			&\leq \mathbb{P}\left\{ \omega\in \Omega: \left(\norm{u_0(\omega)}_{L^2}^2+TL_0^2\right)\exp\left(\bk\Int_0^T\left(\norm{\nabla c_{\eps_\ell}(s)}^4_{L^4}+1\right)ds\right)>\frac{R}{2}\right\}\\
			& +\mathbb{P}\left\{\omega\in \Omega: 2\sup_{r\in[0,T]}\abs{\int_0^{r}(g_1(n_{\eps_\ell}),n_{\eps_\ell}(s))dW_1}\exp\left(\bk\Int_0^T\left(\norm{\nabla c_{\eps_\ell}(s)}^4_{L^4}+1\right)ds\right)> \frac{R}{2}\right\}\\
			&:= P_1+P_2.
		\end{split}
	\end{equation}
	To control $P_1$ and $P_2$, let us note that for any non-negative random variables $X$ and $Y$,
	\begin{equation}\label{equa4.12}
		\left\{\omega\in \Omega: X(\omega)\leq \sqrt{\frac{R}{2}}\right\}\cap\left\{\omega\in \Omega: Y(\omega)\leq \sqrt{\frac{R}{2}}\right\}\subset\left\{\omega\in \Omega: X(\omega)Y(\omega)\leq\frac{R}{2}\right\},
	\end{equation}
	and thanks to the notion of the complementary, the Markov inequality and \re{equa4.5}, the control of $P_1$ goes as follows
	\begin{equation*}
		\begin{split}
			P_1&\leq \mathbb{P}\left\{ \omega\in \Omega: \left(\norm{u_0}_{L^2}^2+TL_0^2\right)>\sqrt{\frac{R}{2}}\right\}
\\
&\qquad {} +\mathbb{P}\left\{ \omega\in \Omega: \exp\left(\bk\Int_0^T\left(\norm{\nabla c_{\eps_\ell}(s)}^4_{L^4}+1\right)ds\right)>\sqrt{\frac{R}{2}}\right\}
\\
			&= \mathbb{P}\left\{ \omega\in \Omega: \left(\norm{u_0}_{L^2}^2+TL_0^2\right)>\sqrt{\frac{R}{2}}\right\}+\mathbb{P}\left\{ \omega\in \Omega: \Int_0^T\left(\norm{\nabla c_{\eps_\ell}(s)}^4_{L^4}+1\right)ds>\frac{\ln\left(\sqrt{\frac{R}{2}}\right)}{\bk}\right\}\\
			&\leq \frac{\sqrt{2}\be \left(\norm{u_0}_{L^2}^2+TL_0^2\right)}{\sqrt{R}}+\frac{\bk\be  \Int_0^T\left(\norm{\nabla c_{\eps_\ell}(s)}^4_{L^4}+1\right)ds}{\ln (\sqrt{R})-\ln \sqrt{2}}\\
			&\leq \frac{\sqrt{2}\be \left(\norm{u_0}_{L^2}^2+TL_0^2\right)}{\sqrt{R}}+\frac{\bk}{\ln (\sqrt{R})-\ln \sqrt{2}},
		\end{split}
	\end{equation*}
	where we have used the fact that, by \re{equa2.3},
	\begin{equation*}
\begin{split}
	\be  \Int_0^T\norm{\nabla c_{\eps_\ell}(s)}^4_{L^4}ds&\leq \be  \sup_{s\in[0,T]}\norm{\nabla c_{\eps_\ell}(s)}^2_{L^2}\Int_0^T\norm{c_{\eps_\ell}}^2_{W^{2,2}}ds\\
	&\leq \be  \sup_{s\in[0,T]}\norm{\nabla c_{\eps_\ell}(s)}^4_{L^2}+\bk\be\left(\Int_0^T\norm{c_{\eps_\ell}}^2_{W^{2,2}}ds\right)^2.
\end{split}	
	\end{equation*}
	For the control of $P_2$, we are using \re{equa4.12}, the Markov inequality and \re{equa4.5} to obtain
	\begin{equation*}
		\begin{split}
			P_2&\leq  \mathbb{P}\left\{ \omega\in \Omega: 2\sup_{r\in[0,T]}\abs{\int_0^{r}(g_1(n_{\eps_\ell}(s)),n_{\eps_\ell}(s))dW_1(s)}>\sqrt{\frac{R}{2}}\right\}\\
			&\qquad+\mathbb{P}\left\{ \omega\in \Omega: \exp\left(\bk\Int_0^T\left(\norm{\nabla c_{\eps_\ell}(s)}^4_{L^4}+1\right)ds\right)>\sqrt{\frac{R}{2}}\right\}\\
			&\leq \frac{2\sqrt{2}\be \Sup_{r\in[0,T]}\abs{\Int_0^{r}(g_1(n_{\eps_\ell}(s)),n_{\eps_\ell}(s))dW_1(s)}}{\sqrt{R}}+\frac{\bk\be  \Int_0^T\left(\norm{\nabla c_{\eps_\ell}(s)}^4_{L^4}+1\right)ds}{\ln (\sqrt{R})-\ln \sqrt{2}}\\
			&\leq \frac{\bk}{\sqrt{R}}+\frac{\bk}{\ln (\sqrt{R})-\ln \sqrt{2}},
		\end{split}
	\end{equation*}
	where we have combined {\re{equa2.3*}}, the Burkholder–Davis–Gundy inequality, the Young inequality, with \re{equa4.5} to get
	\begin{equation*}
		\begin{split}
			\be \Sup_{r\in[0,T]}\abs{\Int_0^{r}(g_1(n_{\eps_\ell}(s)),n_{\eps_\ell}(s))dW_1(s)}& \leq\bk \be \left(\Int_0^{T}\norm{g_1(n_{\eps_\ell}(s))}^2_{\mathcal{L}_2^0}\norm{n_{\eps_\ell}(s)}^2_{L^2}ds\right)^\frac{1}{2}\\
			&\leq \bk \be \left(\Int_0^{T}\norm{n_{\eps_\ell}(s)}^2_{L^2}ds\right)^\frac{1}{2}\leq \bk\left( \be \Int_0^{T}\norm{n_{\eps_\ell}(s)}^2_{L^2}ds\right)^\frac{1}{2}\leq \bk.
		\end{split}
	\end{equation*}
	By using the controls of $P_1$ and $P_2$ we obtain \re{equa4.9*}$_1$ from \re{equa4.11}. The estimate \re{equa4.9*}$_2$ is obtained in a similar way due to \re{equa4.6*}. Let us prove now \re{equa4.10}. To start, let us fix $q>4$ and $R>{\max(2, 16\bk_0^2\bk_1^2\mu^2(T^2+1)T)}$. Defining the duality products $_{(W^{3,2}(\bo))'}\langle \nabla\cdot (F_{\eps_\ell}(n_{\eps_\ell})\nabla c_{\eps_\ell}), \cdot\rangle_{W^{3,2}(\bo)}$ and $_{(W^{3,2}(\bo))'}\langle An_{\eps_\ell}, \cdot\rangle_{W^{3,2}(\bo)}$ respectively by
	\begin{equation*}
		\begin{split}
			&_{(W^{3,2}(\bo))'}\langle \nabla\cdot (F_{\eps_\ell}(n_{\eps_\ell})\nabla c_{\eps_\ell}), \varphi\rangle_{W^{3,2}(\bo)}:=(F_{\eps_\ell}(n_{\eps_\ell})\nabla c_{\eps_\ell},\nabla\varphi) \\
			&_{(W^{3,2}(\bo))'}\langle An_{\eps_\ell}, \varphi\rangle_{W^{3,2}(\bo)}:=(\nabla n_{\eps_\ell},\nabla\varphi)
		\end{split}\qquad \forall\varphi\in W^{3,2}(\bo),
	\end{equation*}
	it follows from \re{equa3.1}, \re{equa4.9**} and the  inequality $(a+b)^2\leq 2a^2+2b^2$ for any $a,b>0$ that
	\begin{equation*}
		\begin{split}
			&\norm{n_{\eps_\ell}}^2_{W^{\frac{2}{q},q}(0,T;(W^{3,2}(\bo))')}\\
			&\leq\bk \norm{u_0}^2_{L^2} +\bk\norm{\int_0^\cdot A n_{\eps_\ell}(s)ds}^2_{W^{1,2}(0,T;(W^{3,2}(\bo))')}\\
			&\quad+\bk \norm{\int_0^\cdot\nabla\cdot(F_{\eps_\ell}(n_{\eps_\ell}(s))\nabla c_{\eps_\ell}(s))ds}^2_{W^{1,2}(0,T;(W^{3,2}(\bo))')}+\bk\norm{\int_0^\cdot n_{\eps_\ell}(s)ds}^2_{W^{1,2}(0,T;(W^{3,2}(\bo))')}\\
			&	\quad+2\bk^2_0\mu^2 \norm{\int_0^\cdot[H_{\eps_\ell}(n_\eps(s))]^2ds}^2_{W^{1,2}(0,T;(W^{3,2}(\bo))')}+\bk\norm{\int_0^\cdot g_1(n_\eps(s))dW_1(s)}^2_{W^{\frac{2}{q},q}(0,T;(W^{3,2}(\bo))')},
		\end{split}
	\end{equation*}
	and consequently
	\begin{equation}\label{equa4.14*}
		\begin{split}
			&	\mathbb{P}\left\{ \omega\in \Omega: \norm{n_{\eps_\ell}}^2_{W^{\frac{2}{q},q}(0,T;(W^{3,2}(\bo))')}>R\right\}\\
			&\leq \mathbb{P}\left\{ \omega\in \Omega:\bk \norm{u_0}^2_{L^2} +\bk\norm{\int_0^\cdot A n_{\eps_\ell}(s)ds}^2_{W^{1,2}(0,T;(W^{3,2}(\bo))')}\right.\\
			&\left.\qquad\qquad+\bk \norm{\int_0^\cdot\nabla\cdot(F_{\eps_\ell}(n_{\eps_\ell}(s))\nabla c_{\eps_\ell}(s))ds}^2_{W^{1,2}(0,T;(W^{3,2}(\bo))')}\right.\\
			&\left.\qquad\qquad+\bk\norm{\int_0^\cdot n_{\eps_\ell}(s)ds}^2_{W^{1,2}(0,T;(W^{3,2}(\bo))')}\right.\\
			&\left.\qquad\qquad+\bk\norm{\int_0^\cdot g_1(n_\eps(s))dW_1(s)}^2_{W^{\frac{2}{q},q}(0,T;(W^{3,2}(\bo))')}>\frac{R}{2}\right\}\\
			&\quad+\mathbb{P}\left\{ \omega\in \Omega: 2\bk^2_0\mu^2 \norm{\int_0^\cdot[H_{\eps_\ell}(n_\eps(s))]^2ds}^2_{W^{1,2}(0,T;(W^{3,2}(\bo))')}>\frac{R}{2}\right\}\\
			&:=P_3+P_4.
		\end{split}
	\end{equation}
	Let us control $P_3$. Thanks to the fact that $\nabla n_{\eps_\ell}=\nabla (n_{\eps_\ell}+1)=\nabla(\sqrt{n_{\eps_\ell}+1})^2$, and the embedding $H^2(\bo)\hookrightarrow L^\infty(\bo)$, we have
	\begin{equation*}
		\begin{split}
			( \nabla n_{\eps_\ell} ,\nabla\varphi)&=2( \sqrt{n_{\eps_\ell}+1} \nabla \sqrt{n_{\eps_\ell}+1},\nabla\varphi)\\
			&\leq2 \norm{\sqrt{n_{\eps_\ell}+1}}_{L^2}\norm{\nabla \sqrt{n_{\eps_\ell}+1}}_{L^2}\norm{\nabla\varphi}_{L^\infty}\\
			&\leq \bk \norm{n_{\eps_\ell}+1}_{L^1}^\frac{1}{2}\norm{\nabla \sqrt{n_{\eps_\ell}+1}}_{L^2}\norm{\nabla\varphi}_{W^{2,2}}\\
			&\leq \bk \norm{n_{\eps_\ell}+1}_{L^1}^\frac{1}{2}\norm{\nabla \sqrt{n_{\eps_\ell}+1}}_{L^2}\norm{\varphi}_{{W^{3,2}}} \text{ for any } \varphi\in  W^{3,2}(\bo).
		\end{split}
	\end{equation*}
	This, and \re{equa4.5} imply
	\begin{equation}\label{equa4.16}
		\begin{split}
			&	\be\norm{\int_0^\cdot A n_{\eps_\ell}(s)ds}^2_{W^{1,2}(0,T;(W^{3,2}(\bo))')}\\
			&\leq\bk\be \sup_{s\in[0,T]}\norm{n_{\eps_\ell}(s)+1}_{L^1}\int_0^T\norm{\nabla\sqrt{n_{\eps_\ell}(s)+1}}_{L^2}^2ds\\
			&\leq\bk\be \sup_{s\in[0,T]}\norm{n_{\eps_\ell}(s)+1}_{L^1}^2+\bk\be\left(\int_0^T\norm{\nabla\sqrt{n_{\eps_\ell}(s)+1}}_{L^2}^2ds\right)^2\leq \bk.
		\end{split}
	\end{equation}
	The H\"older inequality and the continuous embedding of $H^2(\bo)$ into $L^\infty(\bo)$ yield
	\begin{equation*}
	\begin{split}
		( F_{\eps_\ell}(n_{\eps_\ell })\nabla  c_{\eps_\ell },\nabla\varphi)&\leq \norm{ n_{\eps_\ell }}_{L^2}\norm{\nabla  c_{\eps_\ell }}_{L^2}\norm{\nabla\varphi}_{L^\infty}\\
		&\leq \bk \norm{ n_{\eps_\ell }}_{L^2}\norm{\nabla  c_{\eps_\ell }}_{L^2}\norm{\varphi}_{{W^{3,2}}}\text{ for any } \varphi\in  W^{3,2}(\bo),
	\end{split}
	\end{equation*}
	and \re{equa4.5} implies
	\begin{equation}\label{equa4.17}
		\begin{split}
			&\be\norm{\int_0^\cdot\nabla\cdot(F_{\eps_\ell}(n_{\eps_\ell}(s))\nabla c_{\eps_\ell}(s))ds}^2_{W^{1,2}(0,T;(W^{3,2}(\bo))')}\\&\leq  \bk	\be  \int_{0}^{T}  \norm{ n_{\eps_\ell }(s)}^2_{L^2}\norm{\nabla  c_{\eps_\ell }(s)}^2_{L^2}ds\\
			&\leq  \bk\left[	\be\sup_{0\leq s\leq T}\norm{ c_{\eps_\ell }(s)}^4_{W^{1,2}} +\be \left(\int_{0}^{T}  \norm{ n_{\eps_\ell }(s)}^2_{L^2}ds\right)^2\right]\leq  \bk,
		\end{split}
	\end{equation}
	and
	\begin{equation}\label{equa4.18}
		\be\norm{\int_0^\cdot n_{\eps_\ell}(s)ds}^2_{W^{1,2}(0,T;(W^{3,2}(\bo))')}\leq \be \int_{0}^{T}  \norm{ n_{\eps_\ell }(s)}^2_{L^2}ds\leq \bk.
	\end{equation}
	For the stochastic term, we apply the {Jensen} inequality, \cite[Lemma 2.1]{Flandoli} {(which requires $\frac2q<\frac12$, i.e.\ $q>4$)}, and \re{equa2.3*} to obtain
	\begin{equation}\label{equa4.19*}
		\begin{split}
			\be \norm{\int_0^\cdot g_1(n_{\eps_\ell}(s))dW_1(s)}^2_{W^{\frac{2}{q},q}(0,T;(W^{3,2}(\bo))')}&\leq \left(\be \norm{\int_0^\cdot g_1(n_{\eps_\ell}(s))dW_1(s)}^q_{W^{\frac{2}{q},q}(0,T;(W^{3,2}(\bo))')}\right)^\frac{2}{q}\\
			&\leq \bk\left(\be \int_0^T \norm{g_1(n_{\eps_\ell}(s))}^q_{\mathcal{L}^0_2}ds\right)^\frac{2}{q}\leq \bk L_0^2.
		\end{split}
	\end{equation}
	By using the Markov inequality, \re{equa4.16}, \re{equa4.17}, \re{equa4.18} and \re{equa4.19*} as well as the fact that $R>2$, we get
	\begin{equation*}
		P_3\leq \frac{2(\be \norm{u_0}^2_{L^2}+\bk)}{R}\leq \frac{2(\be \norm{u_0}^2_{L^2}+\bk)}{\sqrt{R}}\leq \frac{2(\be \norm{u_0}^2_{L^2}+\bk)}{R^\frac{1}{4}}.
	\end{equation*}
	Note that by  \eqref{equa3.3},  
	\begin{equation*}
	([H_{\eps_\ell}(n_{\eps_\ell})]^2,\varphi)\leq \norm{n_{\eps_\ell}}^2_{L^2}\norm{\varphi}_{L^\infty}\leq \bk_1\norm{n_{\eps_\ell}}^2_{L^2}\norm{\varphi}_{{W^{3,2}}}, \text{ for any } \varphi\in W^{3,2}(\bo). 
	\end{equation*}
	Hence,
	\begin{equation*}
		\norm{\int_0^\cdot[H_{\eps_\ell}(n_\eps(s))]^2ds}^2_{W^{1,2}(0,T;(W^{3,2}(\bo))')}\leq  (T^2+1)\bk_1^2  \int_{0}^{T}  \norm{ n_{\eps_\ell }(s)}^4_{L^2}ds\leq (T^2+1) T\bk_1^2 \sup_{s\in[0,T]}\norm{ n_{\eps_\ell }(s)}^4_{L^2},
	\end{equation*}
	and it follows from \re{equa4.9*} (since $R>\max(2, 16\bk_0^2\bk_1^2\mu^2(T^2+1)T)\geq 2$) that
	\begin{equation*}
		\begin{split}
			P_4&\leq \mathbb{P}\left\{ \omega\in \Omega: 2\bk^2_0\bk_1^2\mu^2(T^2+1) T\sup_{s\in[0,T]}\norm{ n_{\eps_\ell }(s)}^4_{L^2}>\frac{R}{2}\right\}\\
			&= \mathbb{P}\left\{ \omega\in \Omega: \sup_{s\in[0,T]}\norm{ n_{\eps_\ell }(s)}^2_{L^2}>\sqrt{\frac{R}{4\bk^2_0\bk_1^2\mu^2(T^2+1) T}}\right\}\\
			&\leq \frac{\bk}{R^\frac{1}{4}}+\frac{\bk}{\ln (R^\frac{1}{4})-\ln((4\bk^2_0\bk_1^2\mu^2(T^2+1) T)^\frac{1}{4})-\ln \sqrt{2}}.
		\end{split}
	\end{equation*}
	Combining the estimates of $P_3$ and $P_4$, we obtain \re{equa4.10} from \re{equa4.14*}. To prove \re{equa4.9}, let us fix stopping times $\{\tau_\ell\}_{\ell\geq 1}$ such that for any $\ell\geq 1$, $0 \leq  \tau_\ell \leq  T$. Let $\theta>0$ be fixed. It follows from \re{equa3.1}$_2$ and Remark \ref{remark4.2} that
	\begin{equation*}
		\begin{split}
			\be  \norm{c_{\eps_\ell }(\tau_{\ell}+\theta)-c_{\eps_\ell }(\tau_{\ell})}^2_{L^2}
			&\leq
			\bk	\be \abs{\int_{\tau_{\ell}}^{\tau_{\ell}+\theta}   \norm{A_2 c_{\eps_\ell}(s)}_{L^2}ds}^2
			+	\be\abs{\int_{\tau_{\ell}}^{\tau_{\ell}+\theta}  \norm{H_{\eps_\ell}( n_{\eps_\ell }(s))}_{L^2}ds}^2 \\
			&\qquad+	\be \norm{  \int_{\tau_{\ell}}^{\tau_{\ell}+\theta}  g_2( c_{\eps_\ell}(s))d  W_{2 }(s) }^2_{L^2}\\
			&:=II_1+II_2+II_3.
		\end{split}
	\end{equation*}
	By using the {Cauchy--Schwarz} inequality along with \re{equa4.5} and \re{equa3.2}, we get 
	\begin{equation*}
		II_1+II_2\leq \bk\theta\be\int_0^T \norm{ c_{\eps_\ell}(s)}_{{W^{2,2}}}^2ds+\frac{16}{9}\theta\be\int_0^T \norm{ n_{\eps_\ell}(s)}_{L^2}^2ds\leq \bk \theta,
	\end{equation*}
	and the It\^o isometry yields
	\begin{equation*}
		\begin{split}
			II_3
			\leq \be \int_{\tau_{\ell}}^{\tau_{\ell}+\theta} \norm{g_2(c_{\eps_\ell}(s))}^2_{\mathcal{L}_2^0}ds
			&\leq 2\be \int_{\tau_{\ell}}^{\tau_{\ell}+\theta} \norm{g_2(c_{\eps_\ell}(s))-g_2(0)}^2_{\mathcal{L}_2^0}ds+2\theta\norm{g_2(0)}^2_{\mathcal{L}_2^0}\\
			&\leq 2L_2^2\be \int_{\tau_{\ell}}^{\tau_{\ell}+\theta} \norm{c_{\eps_\ell}(s)}^2_{L^2}ds+2\theta\norm{g_2(0)}^2_{\mathcal{L}_2^0}\\
			&\leq 2L^2_2\theta\be \sup_{s\in[0,T]} \norm{c_{\eps_\ell}(s)}^2_{L^2}+2\theta\norm{g_2(0)}^2_{\mathcal{L}_2^0}\leq \bk\theta.
		\end{split}
	\end{equation*}
	Combining these calculations, we derive that
	\begin{equation*}
		\be  \norm{c_{\eps_\ell }(\tau_{\ell}+\theta)-c_{\eps_\ell }(\tau_{\ell})}^2_{L^2}\leq \bk\theta.
	\end{equation*}
	Let us fix now $\xi>0$ and $\eta>0$. By the Markov inequality, we infer that
	\begin{equation*}
		\mathbb{P}\left\{ \omega\in \Omega: \norm{c_{\eps_\ell}(\tau_\ell+\theta){-}c_{\eps_\ell}(\tau_\ell)}_{L^2}\geq \eta\right\}\leq \frac{1}{\eta^2}\be  \norm{c_{\eps_\ell }(\tau_{\ell}+\theta)-c_{\eps_\ell }(\tau_{\ell})}^2_{L^2}\leq \frac{\bk}{\eta^2}\theta.
	\end{equation*}
	Setting $\zeta:= \frac{\eta^2}{\bk}\xi$, the inequality \re{equa4.9} follows from the above inequality and completes the proof of Corollary \ref{cor4.6}.
\end{proof}

The purpose now is to establish the tightness property of the probability law of $n_{\eps_\ell}$ and $c_{\eps_\ell}$. We start by giving the following compactness result, where the details of the proof can be done similarly to \cite[Lemma 3.1]{Brze} because of $D(A_2)\hookrightarrow\hookrightarrow W^{1,2}(\bo)\hookrightarrow L^2(\bo)$, where $\hookrightarrow\hookrightarrow$ means the compact embedding. As a main ingredient, we define the sets $\mathcal{W}^T_n$, $\mathcal{W}^T_c$ and $\mathcal{W}^T_W$ as well as $\mathcal{W}_{u_0}$ and $\mathcal{W}_{v_0}$ respectively by
\begin{equation*}
	\begin{split}
		&\mathcal{W}_n^T:=C([0,T];L_{\text{weak}}^2(\bo))\cap C([0,T];(W^{3,2}(\bo))')\cap  L^2(0,T;L^2(\bo))\cap L^2_{\text{weak}}(0,T;W^{1,2}(\bo)),\\
		&\mathcal{W}_c^T:=C([0,T];{W^{1,2}_{\text{weak}}}(\bo))\cap C([0,T];L^2(\bo))\cap L^2(0,T;W^{1,2}(\bo))\cap L^2_{\text{weak}}(0,T;D(A_2)),\\
		&\mathcal{W}_W^T:=C([0,T];{U_1}), \ \mathcal{W}_{u_0}:=(L^2(\bo))_+,  \ \text{ and }\ \mathcal{W}_{v_0}:=(B^\frac{3}{2}_{4,4,\mathbf{n}}(\bo))_+.
	\end{split}
\end{equation*}

\begin{lemma}\label{lem4.7}
	Let $\Lambda$ be the supremum of the topologies corresponding to $\mathcal{W}_c^T$. Then a set $\mathcal{Q} \subset\mathcal{W}_c^T$ is $\Lambda$-relatively compact if the following {three} conditions hold:
	\begin{equation*}
	\sup_{u\in \mathcal{Q}}\sup_{s\in[0,T]}\norm{u(s)}^2_{W^{1,2}}<\infty,\ 	\sup_{u\in \mathcal{Q}}\int_0^T\norm{u(s)}^2_{W^{2,2}}<\infty,
\end{equation*}
and 
\begin{equation*}
\ \sup_{u\in \mathcal{Q}}\left[\sup_{s,t\in [0,T],\ \abs{t-s}\leq \delta}\norm{u(t)-u(s)}_{L^2}\right]\xrightarrow{\delta{\downarrow} 0}0.
\end{equation*}
\end{lemma}

On the space
\begin{equation*}
	\mathcal{Z}^T:=	\mathcal{W}_n^T\times \mathcal{W}_c^T\times \mathcal{W}_W^T\times  \mathcal{W}_W^T\times\mathcal{W}_{u_0}\times \mathcal{W}_{v_0},
\end{equation*}
we consider the family of probability measures
\begin{equation*}
	\Gamma_{\eps_\ell}:=\Gamma_{\eps_\ell}^n\otimes \Gamma_{\eps_\ell}^c\otimes \Gamma^{W_1}\otimes \Gamma^{W_2}\otimes \Gamma^{u_0}\otimes \Gamma^{v_0},\qquad\ell\geq 1,
\end{equation*}
where $\Gamma_{\eps_\ell}^n$, $\Gamma_{\eps_\ell}^c$, $\Gamma^{W_1}$, $\Gamma^{W_2}$, $\Gamma^{u_0}$, and $\Gamma^{v_0}$ denote respectively the {marginal} probability measure{s} of $n_{\eps_\ell}$, $c_{\eps_\ell}$, $W_1$, $W_2$, $u_0$, and $v_0$ on $(\mathcal{W}_n^T, \mathcal{B}(\mathcal{W}_n^T))$, $(\mathcal{W}_c^T, \mathcal{B}(\mathcal{W}_c^T))$, $(\mathcal{W}_W^T, \mathcal{B}(\mathcal{W}_W^T))$, $(\mathcal{W}_W^T, \mathcal{B}(\mathcal{W}_W^T))$, $(\mathcal{W}_{u_0}, \mathcal{B}(\mathcal{W}_{u_0}))$ and $(\mathcal{W}_{v_0}, \mathcal{B}(\mathcal{W}_{v_0}))$ respectively (where by $\mathcal{B}(X)$, we mean the $\sigma$-algebra induced by the topology of a given topological space $X$). That is, for any $\ell\geq 1$,
\begin{equation*}
	\begin{split}
		&\Gamma_{\eps_\ell}^n(A):=\mathbb{P}\left\{ \omega\in \Omega: n_{\eps_\ell}(\omega)\in A \right\}, \ \forall A\in \mathcal{B}(\mathcal{W}_n^T),\\ &\Gamma_{\eps_\ell}^c(A):=\mathbb{P}\left\{ \omega\in \Omega: c_{\eps_\ell}(\omega)\in A \right\}, \ \forall A\in \mathcal{B}(\mathcal{W}_c^T),\\
		&\Gamma^{W_1}(A):=\mathbb{P}\left\{ \omega\in \Omega: W_1(\omega)\in A \right\}, \ \forall A\in \mathcal{B}(\mathcal{W}_W^T),\\
	&	 \Gamma^{W_2}(A):=\mathbb{P}\left\{ \omega\in \Omega: W_2(\omega)\in A \right\}, \ \forall A\in \mathcal{B}(\mathcal{W}_W^T),\\
		&\Gamma^{u_0}(A):=\mathbb{P}\left\{ \omega\in \Omega: u_0(\omega)\in A \right\}, \ \forall A\in \mathcal{B}(\mathcal{W}_{u_0}),\\ 
	&\Gamma^{v_0}(A):=\mathbb{P}\left\{ \omega\in \Omega: v_0(\omega)\in A \right\}, \ \forall A\in \mathcal{B}(\mathcal{W}_{v_0}),\\
		&\Gamma_{\eps_\ell}(\mathbf{A}):=\mathbb{P}\left\{ \omega\in \Omega: (n_{\eps_\ell}(\omega),c_{\eps_\ell}(\omega),W_{1}(\omega),W_2(\omega){,u_0(\omega),v_0(\omega)})\in \mathbf{A} \right\}, \ \forall \mathbf{A}\in \mathcal{B}(	\mathcal{Z}^T).
	\end{split}
\end{equation*}

\begin{lemma}\label{lem4.8}
	The sequence $(\Gamma_{\eps_\ell})_{\ell\geq 1}$ is tight as probability measures defined on $(\mathcal{Z}^T, \mathcal{B}(\mathcal{Z}^T))$.
\end{lemma}

\begin{proof}
	For any given $\varsigma>0$, our goal is to construct a compact subset $\mathcal{A}_\varsigma$ of $\mathcal{Z}^T$ such that
	\begin{equation*}
		\Gamma_{\eps_\ell}(\mathcal{A}_\varsigma)\geq 1-\varsigma\qquad
		\forall\ell\geq 1.
	\end{equation*}
	Let $\varsigma>0$ be fixed. Note that, since the laws of $W_1$, $W_2$, $u_0$, and $v_0$ are tight as Radon measures on the (Polish) spaces $\mathcal{W}_W^T$, $\mathcal{W}_W^T$, $\mathcal{W}_{u_0}$, and $\mathcal{W}_{v_0}$ respectively, there are compact subsets $\mathcal{A}_{\varsigma,W_1}\subset \mathcal{W}_W^T$, $\mathcal{A}_{\varsigma,W_2}\subset \mathcal{W}_W^T$, $\mathcal{A}_{\varsigma,u_0}\subset\mathcal{W}_{u_0}$ and $\mathcal{A}_{\varsigma,v_0}\subset\mathcal{W}_{v_0}$ such that
	\begin{equation*}
		\Gamma^{W_1}(\mathcal{A}_{\varsigma,W_1})\geq 1-\frac{\varsigma}{8}, \  \Gamma^{W_2}(\mathcal{A}_{\varsigma,W_2})\geq 1-\frac{\varsigma}{8},\  \Gamma^{u_0}(\mathcal{A}_{\varsigma,u_0})\geq 1-\frac{\varsigma}{8} \ \text{ and }\  \Gamma^{v_0}(\mathcal{A}_{\varsigma,v_0})\geq 1-\frac{\varsigma}{8}.
	\end{equation*}
	We point out that $(L^2(\bo))_+$ and $(B^\frac{3}{2}_{4,4,\mathbf{n}}(\bo))_+$ are Polish spaces, because they are complete metric spaces and also separable (see, e.g. \cite[Proposition 3.25]{Brezis} for the separability). Next, we arbitrarily fix $q>4$ and claim that
	\begin{align}
			&	L^\infty(0,T;L^2(\bo))\cap W^{\frac{2}{q},q}(0,T;(W^{3,2}(\bo))')\hookrightarrow\hookrightarrow C([0,T];L_{\text{weak}}^2(\bo)),\label{equa4.20}\\
			&	L^2(0,T;W^{1,2}(\bo))\cap W^{\frac{2}{q},q}(0,T;(W^{3,2}(\bo))')\hookrightarrow\hookrightarrow L^2([0,T];L^2(\bo))\cap C([0,T];(W^{3,2}(\bo))').\notag
	\end{align}
	Indeed, let us arbitrarily fix $\delta\in (0,\frac{1}{q})$ and note that thanks to the continuous embedding
	\begin{equation*}
		W^{\frac{2}{q},q}(0,T;(W^{3,2}(\bo))')\hookrightarrow C^\delta([0,T];(W^{3,2}(\bo))'),
	\end{equation*}
	any bounded subset of $W^{\frac{2}{q},q}(0,T;(W^{3,2}(\bo))')$ is equi-continuous in $C([0,T];(W^{3,2}(\bo))')$. Now, let $\mathcal{A}\subset L^\infty(0,T;L^2(\bo))\cap W^{\frac{2}{q},q}(0,T;(W^{3,2}(\bo))')$ be bounded. Then, it is equi-continuous in $C([0,T];(W^{3,2}(\bo))')$, and since it is also bounded in $L^\infty(0,T;L^2(\bo))$, we can follow \cite[Proposition 4.2]{Brze1} and derive that $\mathcal{A}$ is relatively compact in $C([0,T];L_{\text{weak}}^2(\bo))$. Let us now fix a bounded set $\mathcal{A}_1$ of $L^2(0,T;W^{1,2}(\bo))\cap W^{\frac{2}{q},q}(0,T;(W^{3,2}(\bo))')$. Then, it is equi-continuous in $C([0,T];(W^{3,2}(\bo))')$, and due to its boundedness in $L^2(0,T;W^{1,2}(\bo))$, we can apply \cite[Chapter IV, Theorem 4.1]{Vishik} to conclude that $\mathcal{A}_1$ is relatively compact in $L^2([0,T];L^2(\bo))\cap C([0,T];(W^{3,2}(\bo))')$.
	
{Since the right-hand sides of \eqref{equa4.9*}--\eqref{equa4.10} tend to $0$ as $R\to\infty$ uniformly in $\ell$, one can choose $R_1$, $R_2$, $R_3$ large enough
such that }
	\begin{equation}\label{equa4.21}
		\begin{split}
			&	\mathbb{P}\left\{ \omega\in \Omega: \norm{n_{\eps_\ell}}^2_{L^\infty(0,T;L^2(\bo))}+\norm{n_{\eps_\ell}}^2_{W^{\frac{2}{q},q}(0,T;(W^{3,2}(\bo))')}>2R_1\right\}\\
			&\leq \mathbb{P}\left\{ \omega\in \Omega: \norm{n_{\eps_\ell}}^2_{L^\infty(0,T;L^2(\bo))}>R_1\right\}+\mathbb{P}\left\{ \omega\in \Omega: \norm{n_{\eps_\ell}}^2_{W^{\frac{2}{q},q}(0,T;(W^{3,2}(\bo))')}>R_1\right\}\\
			&\leq\frac{\bk}{\ln (\sqrt{R_1})-\ln \sqrt{2}}+\frac{\bk}{R_1^\frac{1}{4}}+\frac{\bk}{\ln (R_1^\frac{1}{4})-\ln((4\bk^2_0\bk_1^2\mu^2(T^2+1) T)^\frac{1}{4})-\ln \sqrt{2}}\\
			&\leq \frac{\varsigma}{12},
		\end{split}
	\end{equation}
	\begin{align}\label{equa4.22}
			&	\mathbb{P}\left\{ \omega\in \Omega: \norm{n_{\eps_\ell}}^2_{L^2(0,T;W^{1,2}(\bo))}+\norm{n_{\eps_\ell}}^2_{W^{\frac{2}{q},q}(0,T;(W^{3,2}(\bo))')}>3R_2\right\}\notag\\
			&{\leq}\mathbb{P}\left\{ \omega\in \Omega: T\norm{n_{\eps_\ell}}^2_{L^\infty(0,T;L^2(\bo))}+\norm{\nabla n_{\eps_\ell}}^2_{L^2(0,T;L^2(\bo))}+\norm{n_{\eps_\ell}}^2_{W^{\frac{2}{q},q}(0,T;(W^{3,2}(\bo))')}>3R_2\right\}\notag\\
			&\leq \mathbb{P}\left\{ \omega\in \Omega: \norm{n_{\eps_\ell}}^2_{L^\infty(0,T;L^2(\bo))}>\frac{R_2}{T}\right\}+\mathbb{P}\left\{ \omega\in \Omega: \norm{\nabla n_{\eps_\ell}}^2_{L^2(0,T;L^2(\bo))}>R_2\right\}\\
			&\qquad+\mathbb{P}\left\{ \omega\in \Omega: \norm{n_{\eps_\ell}}^2_{W^{\frac{2}{q},q}(0,T;(W^{3,2}(\bo))')}>R_2\right\}\notag\\
			&\leq\frac{\bk}{\ln (\sqrt{R_2})-\ln \sqrt{2T}}+\frac{\bk}{\ln (\sqrt{R_2})-\ln \sqrt{2}}+\frac{\bk}{R_2^\frac{1}{4}}+\frac{\bk}{\ln (R_2^\frac{1}{4})-\ln((4\bk^2_0\bk_1^2\mu^2(T^2+1) T)^\frac{1}{4})-\ln \sqrt{2}}\notag\\
			&\leq \frac{\varsigma}{12},\notag
	\end{align}
and
	\begin{equation}\label{equa4.23}
		\begin{split}
			&\mathbb{P}\left\{ \omega\in \Omega: \norm{ n_{\eps_\ell}}^2_{L^2(0,T;W^{1,2}(\bo))}>2R_3\right\}\\
			&\leq \mathbb{P}\left\{ \omega\in \Omega: \norm{n_{\eps_\ell}}^2_{L^\infty(0,T;L^2(\bo))}>\frac{R_3}{T}\right\}+\mathbb{P}\left\{ \omega\in \Omega: \norm{\nabla n_{\eps_\ell}}^2_{L^2(0,T;L^2(\bo))}>R_3\right\}\\
			&\leq\frac{\bk}{\ln (\sqrt{R_3})-\ln \sqrt{2T}}+\frac{\bk}{\ln (\sqrt{R_3})-\ln \sqrt{2}}+\frac{\bk}{\sqrt{R_3}}\\
			&\leq \frac{\varsigma}{12},
		\end{split}
	\end{equation}
Let us now define $\mathcal{A}_{\varsigma,n}$ as $\mathcal{A}^{1}_{\varsigma,n}\cap \mathcal{A}^{2}_{\varsigma,n}\cap \mathcal{A}^{3}_{\varsigma,n}$, where
	\begin{equation*}
		\begin{split}
			&\mathcal{A}^{1}_{\varsigma,n}:=	\left\{ u\in L^\infty(0,T;L^2(\bo))\cap W^{\frac{2}{q},q}(0,T;(W^{3,2}(\bo))'): \right. \\
			&\left.\hspace{7cm} \norm{{u}}^2_{L^\infty(0,T;L^2(\bo))}+\norm{{u}}^2_{W^{\frac{2}{q},q}(0,T;(W^{3,2}(\bo))')}\leq2R_1 \right\},\\
			&\mathcal{A}^{2}_{\varsigma,n}:=	\left\{ u\in L^2(0,T;W^{1,2}(\bo))\cap W^{\frac{2}{q},q}(0,T;(W^{3,2}(\bo))'):  \right. \\
			&\left.\hspace{7cm} \norm{{u}}^2_{L^2(0,T;W^{1,2}(\bo))}+\norm{{u}}^2_{W^{\frac{2}{q},q}(0,T;(W^{3,2}(\bo))')}\leq3R_2 \right\},\\
			&\mathcal{A}^{3}_{\varsigma,n}:=	\left\{ u\in L^2(0,T;W^{1,2}(\bo)):   \norm{{u}}^2_{L^2(0,T;W^{1,2}(\bo))}\leq2R_3 \right\}.
		\end{split}
	\end{equation*}
	By the compact embedding \re{equa4.20}, the sets $\mathcal{A}^{1}_{\varsigma,n}$ and $\mathcal{A}^{2}_{\varsigma,n}$ are {relatively} compact sets in $C([0,T];L_{\text{weak}}^2(\bo))$ and $L^2([0,T];L^2(\bo))\cap C([0,T];(W^{3,2}(\bo))')$ respectively. In addition, by the Banach-Alaoglu theorem, $\mathcal{A}^{3}_{\varsigma,n}$ is a compact set of $L^2_{\text{weak}}(0,T;W^{1,2}(\bo))$. This justifies that {the closure of} $\mathcal{A}_{\varsigma,n}$ is a compact set of $\mathcal{W}_n^T$. By using \re{equa4.21}, \re{equa4.22} and \re{equa4.23}, we get
	\begin{equation}\label{equa4.13}
		\Gamma^n_{\eps_\ell}(\mathcal{A}_{\varsigma,n})\geq 1-\frac{\varsigma}{4}\quad\text{ for any } \ell\geq 1.
	\end{equation}
	Owing to \re{equa4.9}, we can apply \cite[Lemma 3.8]{Brze} and \cite[Lemma 3.6]{Brze} to derive the existence of a subset $A_{\frac{\varsigma}{12}} \subset C([0,T];L^2(\bo))$ such that
	\begin{equation*}
		\Gamma^c_{\eps_\ell}(A_{\frac{\varsigma}{12}})\geq 1-\frac{\varsigma}{12}\text{ and } \sup_{u\in A_{\frac{\varsigma}{12}}}\left[\sup_{s,t\in [0,T],\ \abs{t-s}\leq \delta}\norm{u(t)-u(s)}_{L^2}\right]\xrightarrow{\delta{\downarrow} 0}0\text{ for any } \ell\geq 1.
	\end{equation*}
	Let $R_4>0$ to be chosen later. It follows from \re{equa4.5} and the Markov inequality that
	\begin{equation*}
		\mathbb{P}\left\{ \omega\in \Omega: \norm{c_{\eps_\ell}}^2_{L^\infty(0,T;W^{1,2}(\bo))}>R_4\right\}
		\leq \frac{\be \Sup_{s\in[0,T]}\norm{ c_{\eps_\ell}(s)}^2_{W^{1,2}}}{R_4}\leq \frac{\bk}{R_4},
	\end{equation*}
	and
	\begin{equation*}
		\mathbb{P}\left\{ \omega\in \Omega: \int_0^T\norm{c_{\eps_\ell}(s)}^2_{W^{2,2}}>R_4\right\}
		\leq \frac{\be \Int_0^T\norm{c_{\eps_\ell}(s)}^2_{W^{2,2}}}{R_4}\leq \frac{\bk}{R_4},
	\end{equation*}
	Choosing $R_4\geq {\frac{12\bk}{\varsigma}}$, we get
	\begin{equation*}
		\mathbb{P}\left\{ \omega\in \Omega: \norm{c_{\eps_\ell}}^2_{L^\infty(0,T;W^{1,2}(\bo))}>R_4\right\}\leq \frac{\varsigma}{12}\quad\text{ for any } \ell\geq 1,
	\end{equation*}
	and
	\begin{equation*}
		\mathbb{P}\left\{ \omega\in \Omega: \int_0^T\norm{c_{\eps_\ell}(s)}^2_{W^{2,2}}>R_4\right\}\leq \frac{\varsigma}{12}\quad\text{ for any } \ell\geq 1.
	\end{equation*}
	Let us define $\mathcal{A}_{\varsigma,c}$ as the closure of the set
	\begin{equation*}
		\left\{ u\in \mathcal{W}_c^T:   \norm{u}^2_{L^\infty(0,T;{W^{1,2}}(\bo))}\leq R_4\text{ and } \norm{u}^2_{L^2(0,T;D(A_2))}\leq R_4 \right\}\cap A_{\frac{\varsigma}{12}},
	\end{equation*}
	in $\mathcal{W}_c^T$. By Lemma \ref{lem4.7}, ${\mathcal{A}_{\varsigma,c}}$ is a compact set in $\mathcal{W}_c^T$. In addition, based on our construction, it is not hard to show that
	\begin{equation*}
		\Gamma^c_{\eps_\ell}(\mathcal{A}_{\varsigma,c})\geq 1-\frac{\varsigma}{4}\quad\text{ for any } \ell\geq 1.
	\end{equation*}
	Setting $\mathcal{A}_\varsigma:=\mathcal{A}_{\varsigma,n}\times \mathcal{A}_{\varsigma,c}\times \mathcal{A}_{\varsigma,W_1}\times \mathcal{A}_{\varsigma,W_2}{\times\mathcal{A}_{\varsigma,u_0}\times\mathcal{A}_{\varsigma,v_0}}$, we deduce from this inequality and \re{equa4.13} {(by a union bound over the complements)} that
	\begin{equation*}
		\begin{split}
			\Gamma_{\eps_\ell}(\mathcal{A}_\varsigma)&
			{\geq 1-\bigl(1-\Gamma^n_{\eps_\ell}(\mathcal{A}_{\varsigma,n})\bigr)-\bigl(1-\Gamma^c_{\eps_\ell}(\mathcal{A}_{\varsigma,c})\bigr)-\sum_{i=1}^2\bigl(1-\Gamma^{W_i}(\mathcal{A}_{\varsigma,W_i})\bigr)}\\
			&\qquad{-\bigl(1-\Gamma^{u_0}(\mathcal{A}_{\varsigma,u_0})\bigr)-\bigl(1-\Gamma^{v_0}(\mathcal{A}_{\varsigma,v_0})\bigr)}\\
			&\geq 1-\tfrac{\varsigma}{4}-\tfrac{\varsigma}{4}-4\cdot\tfrac{\varsigma}{8}=1-\varsigma\quad\text{ for any } \ell\geq 1,
		\end{split}
	\end{equation*}
	and completes the proof of Lemma \ref{lem4.8}.
\end{proof}

\subsection{Passage to the limit and proof of Theorem \ref{theo1.5}}
\label{subsec:limit}
Since {all the spaces comprising $\mathcal Z^T$ admit a countable family of continuous real-valued functions separating points, it follows from  Lemma~\ref{lem4.8} together with the Skorokhod representation theorem for non-metric spaces from \cite[Theorem~2]{Jakubowski} that, }
 there exist $\mathcal{Z}^T$-valued random variables $(\tilde n_{\eps_\ell},\tilde c_{\eps_\ell},\tilde{W}^{\eps_\ell}_1,\tilde{W}^{\eps_\ell}_2,\tilde{u}^{\eps_\ell}_0,\tilde{v}^{\eps_\ell}_0)$, for $\ell\geq1$, and $(n,c, \tilde{W}_1, \tilde{W}_2,\tilde{u}_0,\tilde{v}_0)$, on the probability space
\begin{equation*}
	(\tilde \Omega,\tilde{\mathcal{F}}, \tilde{\mathbb{P}}):=([0,1],\mathcal{B}([0,1]), \text{Leb}),
\end{equation*}
(where Leb is the Lebesgue measure), such that along a subsequence (notationally not relabelled), the laws of $(\tilde n_{{\eps_\ell}},\tilde c_{\eps_\ell},\tilde{W}_1^{\eps_\ell},\tilde{W}_2^{\eps_\ell},\tilde{u}^{\eps_\ell}_0,\tilde{v}^{\eps_\ell}_0)$ and $( n_{\eps_\ell},  c_{\eps_\ell}, W_1,W_2,{u_0,v_0})$ coincide for all $\ell\geq 1$, and the following $\tilde{\mathbb{P}}\text{-a.s.}$ convergences hold
\begin{equation}\label{equa4.38}
	\begin{split}
		&\tilde n_{\eps_\ell}\xrightarrow{\ell\uparrow\infty}n \text{ in } \mathcal{W}_n^T,\qquad\tilde c_{\eps_\ell}\xrightarrow{\ell\uparrow\infty}c \text{ in } \mathcal{W}_c^T,\\
		&\tilde W_1^{\eps_\ell}\xrightarrow{\ell\uparrow\infty}\tilde{W}_1 \text{ in } C([0,T];{U_1}),  \qquad \tilde W_2^{\eps_\ell}\xrightarrow{\ell\uparrow\infty}\tilde{W}_2 \text{ in } C([0,T];{U_1}),\\
		&	\tilde{u}^{\eps_\ell}_0\xrightarrow{\ell\uparrow\infty}\tilde{u}_0 \text{ in } (L^2(\bo))_+,\qquad \text{ and }\qquad \tilde{v}^{\eps_\ell}_0\xrightarrow{\ell\uparrow\infty}\tilde{v}_0 \text{ in } (B^{{\frac32}}_{4,4,\mathbf{n}}(\bo))_+.
	\end{split}
\end{equation}
Due to the equality of laws of $(\tilde{u}^{\eps_\ell}_0,\tilde{v}^{\eps_\ell}_0)$ and $(u_0,v_0)$ for any $\ell\geq 1$, combined with \re{equa2.2}, we get
\begin{equation}\label{equa4.26*}
	\text{Law}_{\tilde{\mathbb{P}}}(\tilde{u}^{\eps_\ell}_0)=\vartheta_1 \text{ and }\text{Law}_{\tilde{\mathbb{P}}}({\tilde{v}^{\eps_\ell}_0})=\vartheta_2 \qquad \forall \ell\geq 1.
\end{equation}
Since almost sure convergence implies convergence in law, we derive from \re{equa4.38} that $\text{Law}_{\tilde{\mathbb{P}}}(\tilde{u}^{\eps_\ell}_0)\xrightarrow{\ell\uparrow\infty} \text{Law}_{\tilde{\mathbb{P}}}(\tilde{u}_0)$ and $\text{Law}_{\tilde{\mathbb{P}}}(\tilde{v}^{\eps_\ell}_0)\xrightarrow{\ell\uparrow\infty} \text{Law}_{\tilde{\mathbb{P}}}(\tilde{v}_0)$. This, and \re{equa4.26*} imply that
\begin{equation}\label{equa4.27}
	\text{Law}_{\tilde{\mathbb{P}}}(\tilde{u}_0)=\vartheta_1 \text{ and }\text{Law}_{\tilde{\mathbb{P}}}({\tilde{v}_0})=\vartheta_2 .
\end{equation}
Thanks to the equality of laws of $(\tilde n_{\eps_\ell},\tilde c_{\eps_\ell})$ and $( n_{\eps_\ell},  c_{\eps_\ell})$ for all $\ell\geq 1$, it follows from \re{equa4.5} that for any ${r}\in[2,4]$,
\begin{equation*}
	\tilde{\be}\norm{\tilde n_{\eps_\ell}}^{{r}}_{L^2(0,T;L^2(\bo))}+\tilde{\be}\norm{\tilde c_{\eps_\ell}}^{{r}}_{L^2(0,T;D(A_2))}\leq \bk.
\end{equation*}
and applying the Fatou lemma, we get for any ${r}\in[2,4]$,
\begin{equation}\label{equa4.39}
	\tilde{\be}\norm{ n}^{{r}}_{L^2(0,T;L^2(\bo))}+\tilde{\be}\norm{ c}^{{r}}_{L^2(0,T;D(A_2))}<\infty.
\end{equation}
The following result concerns the non-negativity property of the new random variables $n$ and $c$.

\begin{lemma}\label{lem4.10*}
	Under the same assumption as in Lemma \ref{lem4.3}, there exists a subsequence of $(\tilde{n}_{\eps_\ell} ,\tilde{c}_{\eps_\ell} )$ still labelled by $\ell$ such that for any $\ell\geq 1$, the following inequalities hold $\tilde{\mathbb{P}}$-a.s.,
	\begin{equation*}
		\tilde{n}_{\eps_\ell} (t,x)\geq 0,\  \tilde{c}_{\eps_\ell} (t,x)\geq 0,\	n(t,x)\geq 0,\  c(t,x)\geq 0,
		\text{ for almost every } (t,x)\in  [0,T]\times \bo.
	\end{equation*}
\end{lemma}

\begin{proof}
	For $\mathcal{C}\in \{L^2(\bo),W^{1,2}(\bo)\}$, we define
	\begin{equation}\label{equa2.67}
		L^2(0,T;\mathcal{C}_+):=\left\{u\in L^2(0,T;\mathcal{C}): u(t,x)\geq 0\ \text{ for  almost every } (t,x)\in (0,T)\times\bo  \right\}.
	\end{equation}
	It is not too difficult to prove that $L^2(0,T;\mathcal{C}_+)$ is a closed set of $L^2(0,T;\mathcal{C})$ with respect to the strong topology and therefore, $\mathcal{C}_+^n\times \mathcal{C}_+^c$ where
	\begin{equation*}
		\mathcal{C}^n_+:=C([0,T];L_{\text{weak}}^2(\bo))\cap C([0,T];(W^{3,2}(\bo))')\cap  L^2(0,T;(L^2(\bo))_+)\cap L^2_{\text{weak}}(0,T;W^{1,2}(\bo)),
	\end{equation*}
	and
	$
	\mathcal{C}_+^c:=C([0,T];{W^{1,2}_{\text{weak}}}(\bo))\cap C([0,T];L^2(\bo))\cap L^2(0,T;(W^{1,2}(\bo))_+)\cap L^2_{\text{weak}}(0,T;D(A_2))
	$
	is a measurable set of $(\mathcal{W}_n^T\times\mathcal{W}_c^T, \mathcal{B}(\mathcal{W}_n^T\times\mathcal{W}_c^T))$. Let us fix $\ell\geq 1$ and denote by $\tilde{\Gamma}_{\eps_\ell}$ the probability law of $(\tilde{n}_{\eps_\ell},\tilde{c}_{\eps_\ell})$ on $\mathcal{W}_n^T\times\mathcal{W}_c^T$, i.e.
	$
	\tilde{\Gamma}_{\eps_\ell}(\mathbf{A}):=\tilde{\mathbb{P}}\left\{ \omega\in \tilde{\Omega}: (\tilde{n}_{\eps_\ell}(\omega,\cdot),\tilde{c}_{\eps_\ell}(\omega,\cdot))\in \mathbf{A} \right\}, \ \forall \mathbf{A}\in \mathcal{B}(\mathcal{W}_n^T\times\mathcal{W}_c^T).
	$
	It follows from the non-negativity property of $n_{\eps_\ell}$ and $c_{\eps_\ell}$ proved in Lemma \ref{lem4.3} that
	\begin{equation*}
		\mathbb{P}\left\{ \omega\in \Omega: (n_{\eps_\ell}(\omega),{c_{\eps_\ell}}(\omega))\in \mathcal{C}_+^n\times \mathcal{C}_+^c \right\}=1.
	\end{equation*}
Using the equality of laws of $(  n_{\eps_\ell},  c_{\eps_\ell})$ and $(\tilde n_{\eps_\ell},\tilde c_{\eps_\ell})$ on $\mathcal{W}_n^T\times\mathcal{W}_c^T$, we get
	\begin{equation*}
		\begin{split}
			\tilde{\mathbb{P}}\left\{ \omega\in \tilde{\Omega}: (\tilde{n}_{\eps_\ell}(\omega),\tilde{c}_{\eps_\ell}(\omega))\in \mathcal{C}_+^n\times \mathcal{C}_+^c \right\}&=\tilde{\Gamma}_{\eps_\ell}(\mathcal{C}_+^n\times \mathcal{C}_+^c)\\
			&=\Gamma_{\eps_\ell}(\mathcal{C}_+^n\times \mathcal{C}_+^c)\\
			&=\mathbb{P}\left\{ \omega\in \Omega: (n_{\eps_\ell}(\omega),c_{\eps_\ell}(\omega))\in \mathcal{C}_+^n\times \mathcal{C}_+^c \right\}
			=1.
		\end{split}
	\end{equation*}
	This implies that $(\tilde{n}_{\eps_\ell} ,\tilde{c}_{\eps_\ell} )\in \mathcal{C}_+^n\times \mathcal{C}_+^c$, $\tilde{\mathbb{P}}$-almost surely and in particular,
	\begin{equation}\label{equa4.40*}
		(\tilde{n}_{\eps_\ell} ,\tilde{c}_{\eps_\ell} )\in L^2(0,T;(L^2(\bo))_+)\times L^2(0,T;(W^{1,2}(\bo))_+)\qquad\tilde{\mathbb{P}}\text{-almost surely}.
	\end{equation}
	From \re{equa4.38} we infer the existence of a subsequence still labelled by $\ell$ of $(\tilde{n}_{\eps_\ell} ,\tilde{c}_{\eps_\ell} )$ such that $\tilde{\mathbb{P}}$-almost surely,
	$
	\tilde{n}_{\eps_\ell}(t,x)\xrightarrow{\ell\uparrow\infty}n(t,x)\text{ and } \tilde{c}_{\eps_\ell}(t,x)\xrightarrow{\ell\uparrow\infty}c(t,x)\text{ for almost every }(t,x)\in (0,T)\times \bo.
	$
	Thanks to \re{equa4.40*}, we get
	$
	0\leq\tilde{n}_{\eps_\ell}(t,x)\xrightarrow{\ell\uparrow\infty}n(t,x)$ and $0\leq \tilde{c}_{\eps_\ell}(t,x)\xrightarrow{\ell\uparrow\infty}c(t,x)\text{ for almost every }(t,x)\in (0,T)\times \bo,
	$
	which implies that
	$
	0\leq n(t,x)$ and  $0\leq c(t,x)$  for almost every $(t,x)\in (0,T)\times \bo,
	$
	and completes the proof of Lemma \ref{lem4.10*}.
\end{proof}

Our goal now is to show that $\tilde W_i^{\eps_\ell}$ and ${\tilde W_i}$, $i=1,2$, are in fact cylindrical Wiener processes. To this end, we need to construct a new filtration. Hereafter, for any $t\in [0,T]$ and any real-valued function $f$ defined on $[0,T]$, $f\mapsto f|_{[0,t]}$ denotes the restriction to the interval $[0,t]$ of $f$. Moreover, we will denote by $\sigma(X)$ the smallest $\sigma$-algebra containing a collection $X$ of subsets of $\tilde{\Omega}$. Let us now consider the following filtration
\begin{equation}\label{equa4.12*}
	\begin{split}
		&\tilde{\mathcal{F}}_t:=\sigma\left(\sigma((n,c,\tilde{W}_1,\tilde{W}_2)|_{[0,t]})\cup\{ N\in \tilde{\mathcal{F}}: \tilde{\mathbb{P}}(N)=0 \}\right),\qquad t\in [0,T],\\
		&\tilde{\mathcal{F}}^{\eps_\ell}_t:=\sigma\left(\sigma((\tilde{n}_{\eps_\ell},\tilde{c}_{\eps_\ell},\tilde{W}_1^{\eps_\ell},\tilde{W}_2^{\eps_\ell})|_{[0,t]})\cup\{ N\in \tilde{\mathcal{F}}: \tilde{\mathbb{P}}(N)=0 \}\right),\quad t\in [0,T], \quad \ell\geq 1.
	\end{split}
\end{equation}
By using the equality of laws of $\tilde W_i^{\eps_\ell}$ and ${W_i}$, $i=1,2$, and the L\'evy martingale characterisation of a Wiener process (see e.g. \cite[Theorem IV.3.6]{Revuz}), we can show that for any $\ell\geq 1$, $\tilde W_i^{\eps_\ell}$ is a cylindrical Wiener process with respect to the filtered probability space $(\tilde \Omega,\tilde{\mathcal{F}}, \{\tilde{\mathcal{F}}^{\eps_\ell}_t\}_{t\in [0,T]}, \tilde{\mathbb{P}})$. A standard argument (see e.g., \cite[Lemma 4.8]{Debussche}) reveals that the almost sure limit $\tilde{W}_i$ of $\tilde W_i^{\eps_\ell}$ is a cylindrical Wiener process relative to the filtered probability space $(\tilde \Omega,\tilde{\mathcal{F}}, \{\tilde{\mathcal{F}}_t\}_{t\in [0,T]} ,\tilde{\mathbb{P}})$. Moreover, by a similar argument as in the proof of Lemma 6.2 in \cite{Galimberti} we can use the equality of laws, Lemma \ref{lem4.10*}, Lemma \ref{lem4.3}, and \re{equa3.1} to show that for any $\ell\geq 1$, the process $(\tilde n_{\eps_\ell},\tilde c_{\eps_\ell}, \tilde W_1^{\eps_\ell},\tilde W_2^{\eps_\ell})$ satisfies the following differential stochastic equations for any $t\in [0,T]$,
\begin{equation}\label{equa4.30}
	\begin{split}
		&	(\tilde n_{\eps_\ell}(t),\varphi)+\int_0^t( \nabla \tilde n_{\eps_\ell}(s),\nabla\varphi)ds-\chi \int_0^t(F_{\eps_\ell}(\tilde n_{\eps_\ell}(s))\nabla \tilde c_{\eps_\ell}(s),\nabla \varphi )ds-(\tilde{u}^{\eps_\ell}_0,\varphi)\\
		&=\int_0^t(\alpha\tilde  n_{\eps_\ell}(s)-\mu [H_{\eps_\ell}(\tilde n_{\eps_\ell}(s))]^2,\varphi)ds+\int_0^t(g_1(\tilde n_{\eps_\ell}(s)),\varphi)d\tilde{W}_1^{\eps_\ell}(s),\ \forall \varphi\in W^{1,2}(\bo), \\
		&	(\tilde c_{\eps_\ell}(t),\psi)+\int_0^t(A_2 \tilde c_{\eps_\ell}(s),\psi)ds\\
		&=(\tilde{v}^{\eps_\ell}_0,\psi)+\int_0^t(H_{\eps_\ell}(\tilde n_{\eps_\ell}(s)),\psi)ds+\int_0^t(g_2(\tilde c_{\eps_\ell}(s)),\psi)d\tilde{W}_2^{\eps_\ell}(s),\quad \forall \psi\in L^{2}(\bo),
	\end{split}
\end{equation}
$\tilde{\mathbb{P}}$-almost surely. With \re{equa4.38} and Lemma \ref{lem4.10*} in hand, it is not too difficult to prove that for any $t\in [0,T]$, the following convergences hold $\tilde{\mathbb{P}}$-almost surely,
\begin{equation}\label{equa4.40}
	\begin{split}
		&( \tilde n_{\eps_\ell }(t),\varphi)\xrightarrow{\ell\uparrow\infty}(  n(t),\varphi),\quad( {\tilde u^{\eps_\ell}_0},\varphi)\xrightarrow{\ell\uparrow\infty}(  \tilde{u}_0,\varphi)\quad\text{ for all } \varphi\in W^{3,2}(\bo),\\
		& ( \tilde c_{\eps_\ell }(t),\psi)\xrightarrow{\ell\uparrow\infty}(  c(t),\psi),\quad\hspace{0.2cm} ( \tilde v_0^{\eps_\ell },\psi)\xrightarrow{\ell\uparrow\infty}(  \tilde{v}_0,\psi)\quad\text{ for all } \psi\in L^2(\bo),\\
		&\int_0^t( \nabla\tilde n_{\eps_\ell }(s),\nabla\varphi)ds\xrightarrow{\ell\uparrow\infty}\int_0^t(  \nabla n(s),\nabla\varphi)ds\quad\text{ for all } \varphi\in W^{1,2}(\bo),\\
		&\int_0^t( \nabla\tilde c_{\eps_\ell }(s),\nabla\psi)ds\xrightarrow{\ell\uparrow\infty}\int_0^t(  \nabla c(s),\nabla\psi)ds\quad\text{ for all } \psi\in W^{1,2}(\bo),\\
		&\int_0^t( \tilde c_{\eps_\ell }(s),\psi)ds\xrightarrow{\ell\uparrow\infty}\int_0^t(   c(s),\psi)ds\quad\text{ for all } \psi\in L^2(\bo),\\
		&\int_0^t( H_{\eps_\ell}(\tilde n_{\eps_\ell }(s)),\psi)ds\xrightarrow{\ell\uparrow\infty}\int_0^t(   n(s),\psi)ds\quad\text{ for all } \psi\in L^2(\bo),\\
		&\int_0^t( [H_{\eps_\ell}(\tilde n_{\eps_\ell }(s))]^2,\varphi)ds\xrightarrow{\ell\uparrow\infty}\int_0^t(   n^2(s),\varphi)ds\quad\text{ for all } \varphi\in W^{2,2}(\bo),\\
		&\int_0^t( \tilde n_{\eps_\ell }(s),\varphi)ds\xrightarrow{\ell\uparrow\infty}\int_0^t(   n(s),\varphi)ds\quad\text{ for all } \varphi\in L^2(\bo),\\
		&\int_0^t( F_{\eps_\ell}(\tilde n_{\eps_\ell }(s))\nabla\tilde c_{\eps_\ell }(s),\nabla\varphi)ds\xrightarrow{\ell\uparrow\infty}\int_0^t(  n(s)\nabla c(s),\nabla\varphi)ds\quad\text{ for all } \varphi\in W^{2,2}(\bo).
	\end{split}
\end{equation}
Note that the last convergence of \re{equa4.40} uses, in addition to \re{equa4.38} and Lemma \ref{lem4.10*}, the fact that
\begin{equation*}
	\int_0^t\norm{F_{\eps_\ell}(\tilde n_{\eps_\ell }(s))-n(s)}^2_{L^2}ds\xrightarrow{\ell\uparrow\infty}0\quad \tilde{\mathbb{P}}\text{-almost surely.}
\end{equation*}
Indeed, we have $F_{\eps_\ell}(\tilde n_{\eps_\ell })-n=\frac{\tilde n_{\eps_\ell }}{1+\eps_\ell\tilde n_{\eps_\ell }}-\frac{n}{1+\eps_\ell n}+\frac{n}{1+\eps_\ell n}-n$ and therefore,
\begin{equation*}
	\begin{split}
		&	\int_0^t\norm{F_{\eps_\ell}(\tilde n_{\eps_\ell }(s))-n(s)}^2_{L^2}ds\\
		&\leq 2\int_0^T\norm{\frac{\tilde n_{\eps_\ell }(s)-n(s)}{(1+\eps_\ell\tilde n_{\eps_\ell }(s))(1+\eps_\ell n (s))}}^2_{L^2}ds+2\eps_\ell^2\int_0^T\norm{\frac{n^2(s)}{1+\eps_\ell n (s)}}^2_{L^2}ds\\
		&\leq  2\int_0^T\norm{\tilde n_{\eps_\ell }(s)-n(s)}^2_{L^2}ds+2\eps_\ell^2\int_0^t\norm{n(s)}^4_{L^4}ds\\
		&\leq  2\int_0^T\norm{\tilde n_{\eps_\ell }(s)-n(s)}^2_{L^2}ds+2\eps_\ell^2\sup_{s\in[0,T]}\norm{n(s)}^2_{L^2}\int_0^T\norm{n(s)}^2_{W^{1,2}}ds\xrightarrow{\ell\uparrow\infty}0\quad \tilde{\mathbb{P}}-\text{a.s.},
	\end{split}
\end{equation*}
where in the last line, we have used \re{equa2.3}, the embedding
\begin{equation*}
	C([0,T];L_{\text{weak}}^2(\bo))\cap C([0,T];(W^{3,2}(\bo))')\hookrightarrow L^\infty(0,T;L^2(\bo)),
\end{equation*}
the convergence \re{equa4.38} and the fact that $\eps_\ell\xrightarrow{\ell\uparrow\infty}0$. As far as the stochastic terms are concerned, we also have the following convergence result up to a subsequence.
\begin{equation}\label{equa4.35}
	\begin{split}
		&\int_0^t(g_1(\tilde n_{\eps_\ell}(s)),\varphi)d\tilde{W}_1^{\eps_\ell}(s)\xrightarrow{\ell\uparrow\infty}\int_0^t(g_1( n(s)),\varphi)d\tilde{W}_1(s)\text{ for any } \varphi\in W^{3,2}(\bo),\\
		&\int_0^t(g_2(\tilde c_{\eps_\ell}(s)),\psi)d\tilde{W}_2^{\eps_\ell}(s)\xrightarrow{\ell\uparrow\infty}\int_0^t(g_2( c(s)),\psi)d\tilde{W}_2(s)\text{ for any } \psi\in L^2(\bo).
	\end{split}
\end{equation}
In fact, using \re{equa2.3*}, we get
\begin{equation*}
	\int_0^t\norm{g_1(\tilde n_{\eps_\ell}(s))-g_1( n(s))}^2_{\mathcal{L}_{2}^0}ds\leq L^2_1 \int_0^T\norm{ \tilde n_{\eps_\ell}(s)-n(s)}_{L^2}^2ds\xrightarrow{\ell\uparrow\infty}0,\quad \tilde{\mathbb{P}}\text{-a.s.},
\end{equation*}
which implies that
$
g_1(\tilde n_{\eps_\ell})\xrightarrow{\ell\uparrow\infty} g_1(n)  \text{ in } L^2(0,t;\mathcal{L}_{2}^0) , \text{ in probability}.
$
Next, from \re{equa4.38}, we infer that $\tilde{W}_1^{\eps_\ell}\xrightarrow{\ell\uparrow\infty}\tilde{W}_1\text{ in } C([0,t];{U_1})$ in probability and realize that the assumptions of \cite[Lemma 2.1]{Debussche1} are therefore fulfilled, with the result that
\begin{equation*}
	\int_0^\cdot g_1(\tilde n_{\eps_\ell}(s))d\tilde{W}_1^{\eps_\ell}(s) \xrightarrow{\ell\uparrow\infty}\int_0^\cdot g_1( n(s))d\tilde{W}_1(s) \text{ in } L^2(0,t; L^2(\bo)) ,
\end{equation*}
in probability. This implies the existence of a subsequence, still denoted by $\ell$, such that \re{equa4.35}$_1$ holds. By a similar argument, we obtain \re{equa4.35}$_2$. Combining \re{equa4.35} and \re{equa4.40}, together with the density of $C^\infty_c(\mathcal{O})$ in $L^2(\mathcal{O})$ and \re{equa4.38}, we can pass to the limit as $\ell\to\infty$ in \re{equa4.30}, for $\varphi\in D(A_2)\hookrightarrow W^{1,2}(\mathcal{O})$. Consequently,
$
((\tilde{\Omega},\tilde{\mathcal{F}},\{\tilde{\mathcal{F}}_t\}_{t\in[0,T]},\tilde{\mathbb{P}}),(\tilde{W}_1,\tilde{W}_2),n,c)
$
is a martingale solution to \eqref{equa1.1}--\eqref{equa1.2}, in the sense of Definition~\ref{defi1.3}, by \re{equa4.27} and \re{equa4.39}. Finally, Lemma~\ref{lem4.10*} yields that $(n,c)$ satisfies \re{equa1.23}, thereby completing the proof of Theorem~\ref{theo1.5}. 

\section{Proof of Theorem \ref{theo1.6}}{\label{sec:proof2}}

We start by proving item (i). Let $(({\Omega},{\mathcal{F}},{\mathbb{P}},{\mathbb{F}}_T), (W_1,W_2),  n_i, c_i)$ , $i=1,2$, be given two martingale solutions to \re{equa1.1}-\re{equa1.2} with  the same   initial {data, i.e.\ $n_1(0)=n_2(0)$ and $c_1(0)=c_2(0)$ $\mathbb P$-a.s.} The differences  $u:=n_1-n_2$ and $v:=c_1-c_2$  satisfy for any $t\in [0,T]$, 
\begin{equation}
	\begin{split}
		&	(u(t),\varphi)+\int_0^t( \nabla u(s),\nabla\varphi)ds-\chi \int_0^t(n_1(s)\nabla v(s),\nabla \varphi )ds\\
		&=\chi \int_0^t(u(s)\nabla c_2(s),\nabla \varphi )ds+\alpha\int_0^t( u(s),\varphi)ds-\mu\int_0^t( (n_1(s)+n_2(s))u(s),\varphi)ds\\
		&	\hspace{0.5cm}+\int_0^t(g_1(n_1(s))-g_1(n_2(s)),\varphi)dW_1(s),\  \varphi\in W^{3,2}(\bo), \\
		&	(v(t),\psi)+\int_0^t(A_2 v(s)-u(s),\psi)ds=\int_0^t(g_2(c_1(s))-g_2(c_2(s)),\psi)dW_2(s),\ \psi\in L^2(\bo).
	\end{split} 
\end{equation}
For every integer $m>0$, we define the stopping time $\tau_m^i$, $i\in \{1,2,3\}$ as 
\begin{equation*}
\begin{split}
&\tau_m^i:=\inf\left\{t\in [0,T]:\sup_{s\in[0,t]}\norm{n_i(s)}^2_{L^2}+\int_0^t\norm{n_i(s)}^2_{W^{1,2}}ds\geq m^2 \right\},\quad i\in \{1,2\}, \\
&\tau_m^3:=\inf\left\{t\in [0,T]:\sup_{s\in[0,t]}\norm{c_2(s)}^2_{W^{1,2}}+\int_0^t\norm{c_2(s)}^2_{W^{2,2}}ds\geq m^2 \right\},
\end{split}
\end{equation*}
with the usual  convention  $\inf\emptyset=T$,  and set
$
\tau_m:=\tau_m^1\wedge\tau_m^2\wedge\tau_m^3.
$

Let $t\in [0,T]$ be arbitrary fixed. Applying the It\^o formula to $\norm{u}_{L^2}^2$ {in the version of \cite[Theorem 1.2]{Pardoux} with the Gelfand triple $W^{3,2}(\bo)\hookrightarrow L^2(\bo)\cong L^2(\bo)\hookrightarrow(W^{3,2}(\bo))'$,} and stopping it at time $\tau_m$, we get for any $r\in [0,t\wedge\tau_m]$,
\begin{equation}\label{equa4.36}
\begin{split}
&\norm{u(r)}^2_{L^2}+2\int_0^{r}\norm{\nabla u}^2_{L^2}ds\\
&=2\chi \int_0^{r}(n_1(s)\nabla v(s),\nabla u(s) )ds+2\chi \int_0^{r\wedge\tau_m}(u(s)\nabla c_2(s),\nabla u(s))ds\\
&	\hspace{0.5cm}+2\alpha\int_0^{r} \norm{u(s)}_{L^2}^2ds-2\mu\int_0^{r}( (n_1(s)+n_2(s))u^2(s),1)ds\\
&	\hspace{0.5cm}+\int_0^{r}\norm{g_1(n_1(s))-g_1(n_2(s))}^2_{\mathcal{L}^0_2}ds+2\int_0^{r}(g_1(n_1(s))-g_1(n_2(s)),u(s))dW_1(s).
\end{split}
\end{equation}
The goal now is to address the various terms in \re{equa4.36} separately.  The H\"older and Young inequality, together with  \re{equa2.3} and the  inequality $\norm{v}_{W^{2,2}}\leq \bk(\norm{\Delta v}_{L^2}+\norm{v}_{W^{1,2}})$ yield
\begin{equation}\label{equa4.37}
\begin{split}
&2\chi(n_1\nabla v,\nabla u )\\
&\leq2\chi \norm{n_1}_{L^4}\norm{\nabla v}_{L^4}\norm{\nabla u}_{L^2}\\
&\leq\frac{1}{4}\norm{\nabla u}_{L^2}^2+\bk \norm{n_1}_{L^2}\norm{n_1}_{W^{1,2}}\norm{\nabla v}_{L^2}\norm{ v}_{W^{2,2}}\\
&\leq\frac{1}{4}\norm{\nabla u}_{L^2}^2+\bk \norm{n_1}_{L^2}\norm{n_1}_{W^{1,2}}\norm{\nabla v}_{L^2}\norm{ \Delta v}_{L^2}+\bk \norm{n_1}_{L^2}\norm{n_1}_{W^{1,2}}\norm{v}_{W^{1,2}}^2\\
&\leq\frac{1}{4}\norm{\nabla u}_{L^2}^2+\frac{1}{2}\norm{ \Delta v}_{L^2}^2+\bk \norm{n_1}^2_{L^2}\norm{n_1}^2_{W^{1,2}}\norm{\nabla v}^2_{L^2}+\bk \norm{n_1}^2_{W^{1,2}}\norm{ v}_{W^{1,2}}^2\\
&\leq\frac{1}{4}\norm{\nabla u}_{L^2}^2+\frac{1}{2}\norm{ \Delta v}_{L^2}^2+\bk (\norm{n_1}^2_{L^2}+1)\norm{n_1}^2_{W^{1,2}}\norm{v}^2_{W^{1,2}},
\end{split}
\end{equation}
and
\begin{equation}\label{equa4.38*}
	\begin{split}
	2\chi	(u\nabla c_2,\nabla u )&\leq 2\chi\norm{u}_{L^4}\norm{\nabla c_2}_{L^4}\norm{\nabla u}_{L^2}\\
		&\leq\frac{1}{8}\norm{\nabla u}_{L^2}^2+\bk \norm{u}_{L^2}\norm{u}_{W^{1,2}}\norm{c_2}_{W^{2,2}}\norm{\nabla c_2}_{L^2}\\
		&\leq\frac{1}{8}\norm{\nabla u}_{L^2}^2+\bk \norm{u}_{L^2}^2\norm{c_2}_{W^{2,2}}\norm{\nabla c_2}_{L^2}+\bk \norm{u}_{L^2}\norm{\nabla u}_{L^2}\norm{c_2}_{W^{2,2}}\norm{\nabla c_2}_{L^2}\\
		&\leq\frac{1}{4}\norm{\nabla u}_{L^2}^2+\bk \norm{u}_{L^2}^2\norm{c_2}^2_{W^{2,2}}+\bk \norm{u}_{L^2}^2\norm{\nabla c_2}_{L^2}^2\norm{c_2}_{W^{2,2}}^2\\
	&\leq\frac{1}{4}\norm{\nabla u}_{L^2}^2+\bk(\norm{ c_2}_{W^{1,2}}^2+1) \norm{c_2}^2_{W^{2,2}}\norm{u}_{L^2}^2.
	\end{split}
\end{equation}
It follows from the H\"older inequality,  the Gagliardo-Nirenberg interpolation inequality \re{equa2.3}, and the Young inequality that 
\begin{equation}\label{equa4.39*}
\begin{split}
-2\mu((n_1+n_2)u^2,1)&\leq 2\mu\norm{n_1}_{L^2}\norm{u}^2_{L^4}+ 2\mu\norm{n_2}_{L^2}\norm{u}^2_{L^4}\\
&\leq \bk\norm{n_1}_{L^2}\norm{u}_{L^2}\norm{ u}_{W^{1,2}}+\bk\norm{n_2}_{L^2}\norm{u}_{L^2}\norm{ u}_{W^{1,2}}\\
&\leq \frac{1}{4}\norm{ \nabla u}^2_{L^2}+\bk(1+\norm{n_1}^2_{L^2}+\norm{n_2}^2_{L^2})\norm{u}^2_{L^2},
\end{split}
\end{equation}
and from \re{equa2.3*} we get 
\begin{equation}\label{equa4.40**}
\norm{g_1(n_1)-g_1(n_2)}^2_{\mathcal{L}^0_2}\leq L^2_1\norm{u}^2_{L^2}.
\end{equation}
Collecting \re{equa4.37}-\re{equa4.40**}, we infer from \re{equa4.36} that  for any $r\in [0,t\wedge\tau_m]$, 
\begin{equation}\label{equa4.41}
	\begin{split}
		&\norm{u(r)}^2_{L^2}+\int_0^{r}\norm{\nabla u(s)}^2_{L^2}ds\\
		&\leq \frac{1}{2} \int_0^{r}\norm{\Delta v(s)}^2_{L^2} ds+\bk \int_0^{r}(\norm{n_1}^2_{L^2}+1)\norm{n_1(s)}^2_{W^{1,2}}\norm{\nabla v(s)}^2_{L^2}ds\\
		&	\quad+(2{|}\alpha{|}+L^2_1)\int_0^{r} \norm{u(s)}_{L^2}^2ds+\bk\int_0^{r}(\norm{ c_2}_{W^{1,2}}^2+1)\norm{c_2(s)}^2_{W^{2,2}}\norm{u(s)}_{L^2}^2ds\\
		&	\quad +\bk\int_0^r(1+\norm{n_1(s)}^2_{L^2}+\norm{n_2(s)}^2_{L^2})\norm{u(s)}^2_{L^2}ds\\
		&\quad+2\abs{\int_0^{r}(g_1(n_1(s))-g_1(n_2(s)),u(s))dW_1(s)}.
	\end{split}
\end{equation}
Before controlling the stochastic term in 	\re{equa4.41}, let us apply the It\^o formula to the process $r\mapsto\norm{v(r)}^2_{L^2}$ (in version of \cite[Theorem 4.32]{DaPrato}) and to the process $r\longmapsto\norm{\nabla v(r)}^2_{L^2}$ (in version of  \cite[Proposition 3.3]{Scarpa} ), and combine  the results  to obtain for any $r\in [0,t\wedge\tau_m]$,
\begin{equation}\label{equa4.42}
	\begin{split}
		&\norm{ v(r)}_{W^{1,2}}^2+2\int_0^{r}\left(\norm{\Delta v(s)}^2_{L^2}+\norm{\nabla v(s)}_{L^2}^2+\norm{ v(s)}_{W^{1,2}}^2\right)ds\\
		&=2\int_0^{r }(u(s),-\Delta v(s)+v(s))ds +\int_0^{r }\norm{g_2(c_1(s))-g_2(c_2(s))}^2_{\mathcal{L}_2^1}ds\\
		&\qquad+2\int_0^{r }(g_2(c_1(s))-g_2(c_2(s)),v(s))_{W^{1,2}}dW_2(s)\\
		&\leq  \frac{1}{2}\int_0^{r }\norm{\Delta v(s)}_{L^2}^2ds+ \int_0^{r }\norm{v(s)}^2_{L^2}ds+\bk\int_0^{r }\norm{u(s)}^2_{L^2}ds\\
		&\qquad +L^2_2\int_0^{r }\norm{v(s)}^2_{W^{1,2}}ds+2\abs{\int_0^{r}(g_2(c_1(s))-g_2(c_2(s)),v(s))_{W^{1,2}}dW_2(s)}.
	\end{split}
\end{equation}
Here, we have used the H\"older inequality, the Young inequality and \re{equa2.4}. 

Adding \re{equa4.41} and \re{equa4.42} together, taking the supremum over $[0,t\wedge\tau_m]$ and take into account the definition of the stopping time $\tau_m$, we obtain 
\begin{equation*}
	\begin{split}
		&\sup_{r\in[0,t\wedge\tau_m]}\norm{u({r})}^2_{L^2}+\sup_{r\in[0,t\wedge\tau_m]}\norm{ v({r})}_{W^{1,2}}^2\\
		&\qquad+\int_0^{t\wedge\tau_m}\left(\norm{\nabla u(s)}^2_{L^2}+\norm{\Delta v(s)}^2_{L^2}+2\norm{\nabla v(s)}_{L^2}^2+\norm{ v(s)}_{W^{1,2}}^2\right)ds\\
		&\leq \bk (m^2+1)\int_0^{t\wedge\tau_m}(\norm{n_1(s)}^2_{W^{1,2}}+1)\norm{ v(s)}^2_{W^{1,2}}ds+\bk(m^2+1)\int_0^{t\wedge\tau_m}(\norm{c_2(s)}^2_{W^{2,2}}+1)\norm{u(s)}_{L^2}^2ds\\
		&	\hspace{0.5cm}+(\bk+2m{^2}\bk +2{|}\alpha{|}+L^2_1)\int_0^{t\wedge\tau_m}\norm{u(s)}_{L^2}^2ds+2\sup_{r\in[0,t\wedge\tau_m]}\abs{\int_0^{r}(g_1(n_1(s))-g_1(n_2(s)),u(s))dW_1(s)}\\
		&	\hspace{0.5cm}+(1+L^2_2)\int_0^{t\wedge\tau_m}\norm{v(s)}_{W^{1,2}}^2ds+2\sup_{r\in[0,t\wedge\tau_m]}\abs{\int_0^{r}(g_2(c_1(s))-g_2(c_2(s)),v(s))_{W^{1,2}}dW_2(s)},
	\end{split}
\end{equation*}
and the Gronwall Lemma and the definition of $\tau_m$ yield
\begin{equation}\label{equa4.43}
	\begin{split}
		&\sup_{r\in[0,t\wedge\tau_m]}\norm{u(r)}^2_{L^2}+\sup_{r\in[0,t\wedge\tau_m]}\norm{ v(r)}_{W^{1,2}}^2\\
		&\qquad+\int_0^{t\wedge\tau_m}\left(\norm{\nabla u(s)}^2_{L^2}+\norm{\Delta v(s)}^2_{L^2}+2\norm{\nabla v(s)}_{L^2}^2+\norm{ v(s)}_{W^{1,2}}^2\right)ds\\
		&\leq2e^{\bk_mT(m^2+1)}\sup_{r\in[0,t\wedge\tau_m]}\abs{\int_0^{r}(g_1(n_1(s))-g_1(n_2(s)),u(s))dW_1(s)}\\
		&	\hspace{0.5cm}+2e^{\bk_mT(m^2+1)}\sup_{r\in[0,t\wedge\tau_m]}\abs{\int_0^{r}(g_2(c_1(s))-g_2(c_2(s)),v(s))_{W^{1,2}}dW_2(s)}\\
		&:=J_8+J_9,
	\end{split}
\end{equation}
{Taking expectations,} the Burkholder–Davis–Gundy inequality combined with the Young inequality and \re{equa2.3*} entail 
\begin{equation}\label{equa4.44}
	\begin{split}
{\be}J_8 &\leq 2\bk e^{\bk_mT(m^2+1)}\be\left(\int_0^{t\wedge\tau_{m}}\norm{g_1(n_1(s))-g_1(n_2(s))}^2_{\mathcal{L}_2^0}\norm{u(s)}^2_{L^2}ds\right)^\frac{1}{2}\\
		&\leq \frac{1}{2}\be\sup_{r\in[0,t\wedge\tau_{m}]}\norm{u(r)}^2_{L^2}+2\bk e^{2\bk_mT(m^2+1)}\be\int_0^{t\wedge\tau_{m}}\norm{u(s)}^2_{L^2}ds,
	\end{split}
\end{equation}
and the same inequalities along with \re{equa2.4} also yield
\begin{equation}\label{equa4.45}
	\begin{split}
	{\be}J_9 &\leq 2\bk e^{\bk_mT(m^2+1)}\be\left(\int_0^{t\wedge\tau_{m}}\norm{g_2(c_1(s))-g_2(c_2(s))}^2_{\mathcal{L}_2^1}\norm{v(s)}^2_{W^{1,2}}ds\right)^\frac{1}{2}\\
		&\leq \frac{1}{2}\be\sup_{r\in[0,t]}\norm{v(r\wedge\tau_{m})}^2_{W^{1,2}}+2\bk e^{\bk_mT(m^2+1)}\be\int_0^{t\wedge\tau_{m}}\norm{v(s)}^2_{W^{1,2}}ds.
	\end{split}
\end{equation}
Taking the mathematical expectation in \re{equa4.43} and inserting \re{equa4.44} and \re{equa4.45} into the result, we infer that
\begin{equation*}
\frac{1}{2}\be\sup_{r\in[0,t\wedge\tau_m]}\norm{u(r)}^2_{L^2}+\frac{1}{2}\be\sup_{r\in[0,t\wedge\tau_m]}\norm{ v(r)}_{W^{1,2}}^2\leq \bk_m\int_0^{{t}}{\be\sup_{r\in[0,s\wedge\tau_m]}}(\norm{u({r})}^2_{L^2}+\norm{v({r})}^2_{W^{1,2}})ds,\ \forall t\in [0,T],
\end{equation*}  
and the Gronwall inequality yield 
\begin{equation*}
	\be\sup_{r\in[0,t\wedge\tau_m]}\norm{u(r)}^2_{L^2}+\be\sup_{r\in[0,t\wedge\tau_m]}\norm{ v(r)}_{W^{1,2}}^2=0,\quad \forall t\in [0,T].
\end{equation*}
 Since for any $t\in [0,T]$, we have  $t\wedge\tau_m\nearrow t$,  $\mathbb{P}$-almost surely, the Fatou lemma implies that $\be\norm{u(t)}^2_{L^2}+\be\norm{ v(t)}_{W^{1,2}}^2=0$ for all $t\in [0,T]$. This shows that $u(t)=v(t)=0$ for every $t\in[0,T]$, almost surely {(first for each $t$, then for all $t$ simultaneously by the continuity of the paths in $(W^{3,2}(\bo))'$ and $L^2(\bo)$)}, and thus completes the proof of item~(i) of Theorem~\ref{theo1.6}. Given item~(i), item~(ii) follows from \cite[Theorem~5.3]{Kunze}. The first part of item~(iii) follows from item~(i) and \cite[Corollary~5.4]{Kunze}. It remains to prove \re{equa1.28}. Since $(n,c)$ is the weak solution to the system \re{equa1.1}-\re{equa1.2} with the initial data $(n_0,c_0)$ satisfying \re{equa1.8}, then $((\Omega,\mathcal{F},\mathbb{P},\mathbb{F}_T), (W_1,W_2),  n, c)$   is a martingale solution to the system \re{equa1.1}-\re{equa1.2} with initial data $(\text{Law}_\mathbb{P}(n_0),\text{Law}_\mathbb{P}({c_0}))$.   Let $((\bar{\Omega},\bar{\mathcal{F}},\bar{\mathbb{P}},\bar{\mathbb{F}}_T), (\bar{W}_1,\bar{W}_2), \bar n,\bar c)$ be a martingale solution to the system \re{equa1.1}-\re{equa1.2}  obtained from Theorem \ref{theo1.5} with initial data $(\text{Law}_\mathbb{P}(n_0),\text{Law}_\mathbb{P}({c_0}))$. By item (ii) of Theorem \ref{theo1.6}, 
 \begin{equation}\label{equa3.12*}
 	\text{Law}_{\mathbb{P}}(n, c)=\text{Law}_{\bar{\mathbb{P}}}(\bar n,\bar c) \text{ on } \bz_1\times \bz_2, 
 \end{equation}
 where $\bz_1:=C([0,T];(W^{3,2}(\bo))')\cap C([0,T];L_{\text{weak}}^2(\bo))\cap L^2(0,T;W^{1,2}(\bo)) $ and \\
 $\bz_2:=C([0,T];L^2(\bo))\cap C([0,T];W_{\text{weak}}^{1,2}(\bo))\cap L^2(0,T;D(A_2))$. For $\mathcal{C}\in \{W^{1,2}(\bo), D(A_2)\}$, we consider the space $	L^2(0,T;\mathcal{C}_+)$  defined by \re{equa2.67} and 
 the measurable sets 
 \begin{equation*}
 	\mathcal{A}^n_+:=C([0,T];L_{\text{weak}}^2(\bo))\cap C([0,T];(W^{3,2}(\bo))')\cap  L^2(0,T;(W^{1,2}(\bo))_+),
 \end{equation*}
 and 
 $
 \mathcal{A}_+^c:=C([0,T];{W^{1,2}_{\text{weak}}}(\bo))\cap C([0,T];L^2(\bo))\cap L^2(0,T;(D(A_2))_+)
 $ on $(\mathcal{Z}_1, \mathcal{B}(\mathcal{Z}_1))$ and \\ $(\mathcal{Z}_2, \mathcal{B}(\mathcal{Z}_2))$ respectively.
 It follows from \re{equa1.23} and {\re{equa3.12*}} that  
 \begin{equation*}
 	\begin{split}
 		\mathbb{P}\left\{ \omega\in  \Omega: (n(\omega),c(\omega))\in \mathcal{A}_+^n\times \mathcal{A}_+^c \right\}&=	\text{Law}_{\mathbb{P}}(n, c)[\mathcal{A}_+^n\times \mathcal{A}_+^c]\\
 		&=\text{Law}_{\bar{\mathbb{P}}}(\bar n, \bar c)[\mathcal{A}_+^n\times \mathcal{A}_+^c]\\
 		&=\bar{\mathbb{P}}\left\{ \omega\in \bar{\Omega}: (\bar n(\omega),\bar c(\omega))\in \mathcal{A}_+^n\times \mathcal{A}_+^c \right\}
 		=1.
 	\end{split}
 \end{equation*}
 This implies that $(n,c )\in \mathcal{A}_+^n\times \mathcal{A}_+^c$, $\mathbb{P}$-almost surely and therefore,
 \begin{equation*}
 	(n,c )\in L^2(0,T;(W^{1,2}(\bo))_+)\times L^2(0,T;(D(A_2))_+)\quad\mathbb{P}\text{-a.s.}. 
 \end{equation*}
This implies \re{equa1.28} and completes the proof of Theorem \ref{theo1.6}.
\appendix
\section{Construction of the regularised solution}
\label{app:reg}

This appendix contains the proof of Proposition \ref{propo3.5}: the construction of a global solution to the regularised stochastic variational system \re{equa3.1}. The proof proceeds in four steps. In the first step, we construct a local solution by combining a localization procedure with the Banach fixed-point theorem. In the second step, we prove the pathwise uniqueness of local solutions. In the third step, we construct a unique maximal local solution and establish a blow-up criterion. In the fourth step, we use the a priori bound to upgrade the maximal local solution to a global one, thereby completing the proof of Proposition \ref{propo3.5}.
We begin by recalling the notions of local, maximal local, and global solutions to the system \re{equa3.1}.

\begin{definition}\label{defi3.1}
	A local solution to the system \re{equa3.1} is given by a triple $(n_\eps,c_\eps,\tau_\eps)$ consisting of a strictly positive stopping time $\tau_\eps$, and a $\mathbb{F}_T$-progressively measurable, $L^2(\bo)\times B_{4,4}^\frac{3}{2}(\bo)$-valued process $\{(n_\eps(t),c_\eps(t))\}_{t\in [0,T]}$ satisfying
	\begin{equation}\label{equa2.6*}
		\begin{split}
			&(n_\eps(\cdot\wedge\tau_\eps),c_\eps(\cdot\wedge\tau_\eps))\in C([0,T]; L^2(\bo)\times B_{4,4,\mathbf{n}}^\frac{3}{2}(\bo)),\\
			&	(\mathbf{1}_{[0,\tau_\eps]}(\cdot)n_\eps,\mathbf{1}_{[0,\tau_\eps]}(\cdot)c_\eps)\in L^2(0,T; W^{1,2}(\bo)\times D(A_2)),
		\end{split} \qquad\mathbb{P}\text{-a.s.},
	\end{equation}
	and
	\begin{equation}\label{equa3.20}
		\begin{split}
			&\be \sup_{s\in[0,T\wedge\tau_\eps]}\norm{n_\eps(s)}^4_{L^2}+\be\left(\int_0^{T\wedge\tau_\eps}\norm{n_\eps(s)}^2_{W^{1,2}}ds\right)^2<\infty,\\
			&\be \sup_{s\in[0,T\wedge\tau_\eps]}\norm{c_\eps(s)}^4_{B_{4,4}^\frac{3}{2}}+\be\left(\int_0^{T\wedge\tau_\eps}\norm{A_2c_\eps(s)}^2_{L^2}ds\right)^2<\infty,
		\end{split}
	\end{equation}
	as well as the following identities for all  $t\in[0,T]$,
	\begin{equation}\label{equa3.5*}
		\begin{split}
			&	(n_\eps(t\wedge\tau_\eps),\varphi)+\int_0^{t\wedge\tau_\eps}( \nabla n_\eps(s),\nabla\varphi)ds-\chi \int_0^{t\wedge\tau_\eps}(F_\eps(n^+_\eps(s))\nabla c_\eps(s),\nabla \varphi )ds-(u_0,\varphi)\\
			&=\int_0^{t\wedge\tau_\eps}(\alpha n_\eps(s)-\mu [H_\eps(n^+_\eps(s))]^2,\varphi)ds+\int_0^{t\wedge\tau_\eps}(g_1(n^+_\eps),\varphi)dW_1(s)\  \forall \varphi\in W^{1,2}(\bo),  \   \mathbb{P}\text{-a.s.}, \\
			&	(c_\eps(t\wedge\tau_\eps),\psi)+\int_0^{t\wedge\tau_\eps}(A_2 c_\eps(s),\psi)ds\\
			&=(v_0,\psi)+\int_0^{t\wedge\tau_\eps}(H_\eps(n^+_\eps(s)),\psi)ds+\int_0^{t\wedge\tau_\eps}(g_2(c_\eps(s)),\psi)dW_2(s)\ \quad\forall \psi\in L^2(\bo),  \quad    \mathbb{P}\text{-a.s.}.
		\end{split}
	\end{equation}
\end{definition}

\begin{definition}\label{defi3.3}
	A maximal local solution to the system \re{equa3.1} is given by a triple $(n_\eps,c_\eps,  \tau_\eps^\text{max})$ consisting of a strictly positive stopping time $\tau^\text{max}_\eps$ and a $\mathbb{F}_T$-progressively measurable, $L^2(\bo)\times B_{4,4}^\frac{3}{2}(\bo)$-valued process $\{(n_\eps(t),c_\eps(t))\}_{t\in [0,T]}$ such that,
	\begin{itemize}
		\item[(i)] there exists an increasing sequence $\{\tau_{\eps,\kappa}\}_{\kappa\geq 2\chi}$ 
		 of $\mathbb{F}_T$-adapted strictly positive stopping times such that $\tau_{\eps,\kappa}\nearrow\tau_\eps^\text{max}$ almost surely and for any integer $\kappa\geq 2\chi$, the triple $(n_\eps,c_\eps,  \tau_{\eps,\kappa})$ is a local solution to the system \re{equa3.1};
		\item[(ii)] for all local solutions $(\tilde{n}_\eps,\tilde{c}_\eps,  \tilde{\tau}_\eps)$ to the system \re{equa3.1}, we have $\tilde{\tau}_{\eps}\leq \tau_{\eps}^\text{max}$ almost surely and $(n_\eps(t\wedge\tilde{\tau}_{\eps}),c_\eps(t\wedge\tilde{\tau}_{\eps}))=(\tilde{n}_\eps(t\wedge\tilde{\tau}_{\eps}),\tilde{c}_\eps(t\wedge\tilde{\tau}_{\eps}))$ for all $t\in [0,T]$, almost surely.
	\end{itemize}
\end{definition}

\begin{definition}\label{defi3.4}
	A global solution to the system \re{equa3.1} is a maximal local solution $(n_\eps,c_\eps, \tau_\eps^\text{max})$ such that $\tau_\eps^\text{max}={T}$ almost surely. 
\end{definition}

\begin{remark}\label{remark2.4}
	If $(n_\eps,c_\eps,  \tau_\eps^\text{max})$ is a global solution to the system \re{equa3.1}, then the couple $(n_\eps,c_\eps)$ is a $\mathbb{F}_T$-progressively measurable, $L^2(\bo)\times B_{4,4}^\frac{3}{2}(\bo)$-valued process satisfying \re{equa3.1}.
		Indeed, note that by Definition \ref{defi3.4}, $\tau_\eps^\text{max}={T}$ almost surely and by Definition \ref{defi3.3}, there exists a family of stopping times $\{\tau_{\eps,\kappa}\}_{\kappa\geq 2\chi}$ increasing to $\tau_\eps^\text{max}={T}$ such that for any integer $\kappa\geq 2\chi$, the triple $(n_\eps,c_\eps,  \tau_{\eps,\kappa})$ is a local solution to \re{equa3.1}. In particular, it satisfies \re{equa2.6*} and \re{equa3.5*} with $\tau_{\eps,\kappa}$ in place of $\tau_{\eps}$ and therefore, we can pass to the limit almost surely in the corresponding equation of \re{equa3.5*} and obtain \re{equa3.1}.
\end{remark}
We now provide the detailed proof of Proposition~\ref{propo3.5}.

$\bullet$ \textbf{Step 1: Existence of local solutions.} To construct local solutions of problem \re{equa3.1}, we have to modify the cross-diffusion coefficient by introducing a family of smooth cut-off functions $\theta_R:=\theta(\frac{\cdot}{R})$, $R>0$, where $\theta\in C^2([0,{\infty}),[0,1])$ is a fixed function satisfying
\begin{enumerate}
	\item[(i):] $\theta(z)=1$,\ $z\in [0,1]$, and $\theta(z)=0$, \ $z>2$,
	\item[(ii):] $\Sup_{z{\geq 0}}{|}\theta'(z){|}\leq \bk<\infty$.
\end{enumerate}
By using the Banach fixed-point theorem, we are going to construct for any arbitrary fixed
{$
R>0,
$ }
a unique $\mathbb{F}_T$-adapted $((n^R_\eps(t),c^R_\eps(t)))_{t\in [0,T]}$ with $\mathbb{P}$-a.s. continuous sample paths in $L^2(\bo)\times B_{4,4,\mathbf{n}}^\frac{3}{2}(\bo)$, satisfying for any $t\in[0,T]$,
\begin{equation}\label{equa3.5}
	\begin{split}
		&	(n^R_\eps(t),\varphi)+\int_0^t( \nabla n^R_\eps(s),\nabla\varphi)ds-\chi \int_0^t\theta_R(\norm{c^R_\eps(s)}_{W^{1,4}})(F_\eps((n^R_\eps(s))^+)\nabla c^R_\eps(s),\nabla \varphi )ds\\
		&=(u_0,\varphi)+\int_0^t(\alpha n^R_\eps(s)-\mu [H_\eps((n^R_\eps(s))^+)]^2,\varphi)ds\\
		&\hspace{2cm}+\int_0^t(g_1((n^R_\eps(s))^+),\varphi)dW_1(s), \ \forall \varphi\in W^{1,2}(\bo), \quad    \mathbb{P}\text{-a.s.}, \\
		&	(c^R_\eps(t),\psi)+\int_0^t(A_2 c^R_\eps(s),\psi)ds\\
		&=(v_0,\psi)+\int_0^t(H_\eps((n^R_\eps(s))^+),\psi)ds+\int_0^t(g_2(c^R_\eps(s)),\psi)dW_2(s),\quad\forall \psi\in L^2(\bo),  \quad    \mathbb{P}\text{-a.s.},
	\end{split}
\end{equation}
and later on, construct a local solution of system \re{equa3.1} from the process $(n^R_\eps(t),c^R_\eps(t))_{t\in [0,T]}$.

\begin{claim}\label{lem3.1}
	For any $R>2\chi$, there exists a unique $\mathbb{F}_T$-adapted $((n^R_\eps(t),c^R_\eps(t)))_{t\in [0,T]}$ with $\mathbb{P}$-a.s. continuous sample paths in $L^2(\bo)\times B_{4,4,\mathbf{n}}^\frac{3}{2}(\bo)$, satisfying \re{equa3.5} and such that
	\begin{equation}\label{equa3.6*}
		\begin{split}
			&n_\eps^R\in L^4_{\mathbb{F}_T}(\Omega ; C([0,T]; L^2(\bo)))\cap L^4_{\mathbb{F}_T}(\Omega ; L^2(0,T; W^{1,2}(\bo))),\\
			&c_\eps^R\in L^4_{\mathbb{F}_T}(\Omega; C([0,T]; B_{4,4,\mathbf{n}}^\frac{3}{2}(\bo)))\cap L^4_{\mathbb{F}_T}(\Omega ; L^4(0,T; D(A_4))).
		\end{split}
	\end{equation}
\end{claim}

 \textbf{Proof of Claim \ref{lem3.1}.}
	Let us fix $R>2\chi$. We will first use the Banach fixed-point theorem to construct a unique $\mathbb{F}_T$-adapted $((n^R_\eps(t),c^R_\eps(t)))_{t\in [0,T]}$ with $\mathbb{P}$-a.s. continuous sample paths in $L^2(\bo)\times W^{1,2}(\bo)$ satisfying \re{equa3.5} and later on, prove that the process $\{c^R_\eps(t)\}_{t\in [0,T]}$ is in fact a $B_{4,4,\mathbf{n}}^\frac{3}{2}(\bo)$-valued continuous process. Let $T_0>0$ to be fixed later. Let us define the Banach spaces
	\begin{equation*}
		\begin{split}
			&\bx_{T_0}:=L^4_{\mathbb{F}_{T_0}}(\Omega ; C([0,T_0]; L^2(\bo)))\cap L^4_{\mathbb{F}_{T_0}}(\Omega ; L^2(0,T_0; W^{1,2}(\bo))),\\
			&\by_{T_0}:=L^4_{\mathbb{F}_{T_0}}(\Omega ; C([0,T_0]; W^{1,2}(\bo)))\cap L^4_{\mathbb{F}_{T_0}}(\Omega ; L^2(0,T_0; D(A_2))),
		\end{split}
	\end{equation*}
	endowed respectively with the norm
	\begin{equation}
		\begin{split}
			\norm{u}^4_{\bx_{T_0}}:=\be\sup_{t\in [0,T_0]}\abs{u(t)}_{L^2}^4+\be\left(\int_0^{T_0}\norm{u(s)}^2_{W^{1,2}}ds\right)^{2}\qquad u\in \bx_{T_0}, \\
			\norm{v}^4_{\by_{T_0}}:=\be\sup_{t\in [0,T_0]}\abs{v(t)}^4_{W^{1,2}}+\be\left(\int_0^{T_0}\norm{v(s)}_{W^{2,2}}^2ds\right)^{2}\qquad v\in \by_{T_0}. \\
		\end{split}
	\end{equation}
	For each $(u,v)\in \bx_{T_0}\times \by_{T_0}$, let us consider the system of stochastic differential equations with additive noise that consists of finding $\mathbb{F}_{T_0}$-adapted processes $(m_\eps^R(t))_{t\in [0,T_0]}$ and $(d_\eps^R(t))_{t\in [0,T_0] }$ such that for all $t\in[0,T_0]$,
	\begin{equation}\label{equa3.6}
		\begin{split}
			&	(m^R_\eps(t),\varphi)+\int_0^t( \nabla m^R_\eps(s),\nabla\varphi)ds-\chi \int_0^t\theta_R(\norm{d^R_\eps(s)}_{W^{1,4}})(F_\eps((u(s))^+)\nabla d^R_\eps(s),\nabla \varphi )ds\\
			&=(u_0,\varphi)+\int_0^t(\alpha u(s)-\mu [H_\eps((u(s))^+)]^2,\varphi)ds\\
			&\quad+\int_0^t(g_1(u^+(s)),\varphi)dW_1(s),  \ \forall \varphi\in W^{1,2}(\bo),\quad    \mathbb{P}\text{-a.s.}, \\
			&	(d^R_\eps(t),\psi)+\int_0^t(A_2 d^R_\eps(s),\psi)ds\\
			&=(v_0,\psi)+\int_0^t(H_\eps((u(s))^+),\psi)ds+\int_0^t(g_2(v(s)),\psi)dW_2(s),\ \forall \psi\in L^2(\bo),  \quad    \mathbb{P}\text{-a.s.}.
		\end{split}
	\end{equation}
	For any $t\in[0,T_0]$, define
	$
	f_\eps(t):=H_\eps\big((u(t))^+\big)$ and $
	B_\eps(t):=g_2(v(t)).
$
	Then \re{equa3.6}$_2$ is equivalently written as
	\begin{equation*}
		d^R_\eps(t)+\int_0^tA_2 d^R_\eps(s)ds=v_0+\int_0^tf_\eps(s)ds+\int_0^tB_\eps(s)dW_2(s) \text{ in } L^2(\bo), \ \forall t\in [0,T_0], \quad    \mathbb{P}\text{-a.s.},
	\end{equation*}
	that is a linear equation in $L^2(\bo)$ driven by additive noise. By \re{equa3.4}, we have $\norm{f_\eps(t)}_{L^2}\leq \bk_\eps$ for any $t\in [0,T]$. In addition, we can use \re{equa2.1} and \re{equa1.5} to observe that
	\begin{equation} \label{equa2.15}
		\be \norm{v_0}^4_{W^{1,2}}\leq \bk \be\norm{v_0}^4_{B^\frac{3}{2}_{4,4}}<\infty.
	\end{equation}
	Hence, using the Galerkin method, we can adapt the idea of the proof of \cite[Proposition 3.3]{Tachim} to construct a unique $\mathbb{F}_{T_0}$-adapted process $(d_\eps^R (t))_{t\in [0,T_0]}$ satisfying \re{equa3.6}$_2$ and such that
	\begin{equation}\label{equa3.7}
		d_\eps^R\in \by_{T_0}.
	\end{equation}
	With the process $(d_\eps^R (t))_{t\in [0,T_0]}$ in hand, we can replace it in \re{equa3.6}$_1$ and obtain the following stochastic variational equation with the process $(m_\eps^R (t))_{t\in [0,T_0]}$ as unknown,
	\begin{equation*}
		(m_\eps^R(t),\varphi)=({u_0},\varphi)+\int_0^t{}_{(W^{1,2}(\bo))'}\langle \mathcal{A}(t,m^R_\eps),\varphi\rangle_{W^{1,2}(\bo)}+ \int_0^t(B(s,m_\eps^R),\varphi)dW_1(s), \quad t\in [0,T_0],
	\end{equation*}
	where
	\begin{equation*}
		\mathcal{A}:\Omega\times [0,T_0]\times W^{1,2}(\bo)\longrightarrow (W^{1,2}(\bo))'\qquad B:\Omega\times[0,T_0]\times  W^{1,2}(\bo)\longrightarrow \mathcal{L}_2^0(\bo),
	\end{equation*}
	are mappings defined respectively by for any $(\omega,t ,z)\in\Omega\times  [0,T_0]\times W^{1,2}(\bo)$,
	\begin{equation*}
		\begin{split}
			_{(W^{1,2}(\bo))'}\langle \mathcal{A}(\omega,t ,z),\varphi\rangle_{W^{1,2}(\bo)}:=&-( \nabla z,\nabla\varphi)+\chi \theta_R(\norm{d^R_\eps(\omega,t )}_{W^{1,4}})(F_\eps((u(\omega,t ))^+)\nabla d^R_\eps(\omega,t ),\nabla \varphi )\\
			&\qquad +(\alpha u(\omega,t )-\mu [H_\eps((u(\omega,t ))^+)]^2,\varphi),\qquad \varphi\in W^{1,2}(\bo),
		\end{split}
	\end{equation*}
	and
	\begin{equation*}
		(B(\omega,t ,z))=g_1(u^+(\omega,t )).
	\end{equation*}
	With \re{equa2.3*} in hand, it is straightforward to show that $\mathcal{A}$ and $B$ satisfy (H1)-(H4) of \cite{Prevot}. So, due to the regularity of $u_0$, we can apply \cite[Theorem 4.2.4]{Prevot} and infer the existence of a unique continuous $L^2(\bo)$-valued $\mathbb{F}_{T_0}$-adapted process $m_\eps^R$ satisfying \re{equa3.6}$_1$ and
	\begin{equation*}
		m_\eps^R\in L^2(\Omega ; C([0,T_0]; L^2(\bo)))\cap L^2(\Omega ; L^2(0,T_0; W^{1,2}(\bo))).
	\end{equation*}
	Thanks to the It\^o formula in the version of \cite[Theorem 1.2]{Pardoux}, we can show that
	\begin{equation}\label{equa3.8}
		m_\eps^R\in \bx_{T_0}.
	\end{equation}
	Combining \re{equa3.7} and \re{equa3.8} we derive that the map
	\begin{equation*}
		\varUpsilon:(u,v)\in\bx_{T_0}\times \by_{T_0}\longmapsto \varUpsilon(u,v):=(m_\eps^R, d_\eps^R)\in \bx_{T_0}\times \by_{T_0},
	\end{equation*}
	is well defined. Our aim now is to prove that $\varUpsilon$ has a unique fixed-point by showing that it satisfies the condition for applying the contraction mapping principle. Let $(u_1,v_1)$ and $(u_2,v_2)$ be two elements of $\bx_{T_0}\times\by_{T_0}$. Setting
	\begin{equation*}
		(m_{i,\eps}^R, d_{i,\eps}^R)=\varUpsilon(u_i,v_i)\qquad \forall i\in \{1,2\},
	\end{equation*}
	the differences $M_\eps^R:=m_{1,\eps}^R-m_{2,\eps}^R$ and $D_\eps^R:=d_{1,\eps}^R-d_{2,\eps}^R$ satisfy for any $t\in [0,T_0]$,
	\begin{equation*}
		\begin{split}
			&	(M^R_\eps(t),\varphi)+\int_0^t( \nabla M^R_\eps(s),\nabla\varphi)ds\\
			&=\chi \int_0^t\left(\theta_R\left(\norm{d^R_{1,\eps}}_{W^{1,4}}\right)F_\eps((u_1)^+)\nabla d^R_{1,\eps}(s)-\theta_R\left(\norm{d^R_{2,\eps}}_{W^{1,4}}\right)F_\eps((u_2(s))^+)\nabla d^R_{2,\eps}(s),\nabla \varphi \right)ds\\
			&\qquad+\alpha\int_0^t( u_1(s)-u_2(s),\varphi)ds-\mu\int_0^t( [H_\eps((u_1(s))^+)]^2-[H_\eps((u_2(s))^+)]^2,\varphi)ds\\
			&\qquad+\int_0^t(g_1(u_1{^+}(s))-g_1(u_2{^+}(s)),\varphi)dW_1(s),\qquad\forall \varphi\in W^{1,2}(\bo), \quad \mathbb{P}\text{-a.s.}, \\
			&	(D^R_\eps(t),\psi)+\int_0^t(A{_2} D^R_\eps(s),\psi)ds\\
			&=\int_0^t(H_\eps((u_1(s))^+)-H_\eps((u_2(s))^+),\psi)ds+\int_0^t(g_2(v_1(s))-g_2(v_2(s)),\psi)dW_2(s),\ \forall \psi\in L^2(\bo),\  \mathbb{P}\text{-a.s.}.
		\end{split}
	\end{equation*}
	Let $t\in [0,T_0]$ be fixed. By applying the It\^o formula in the version of \cite[Theorem 4.32]{DaPrato} to the process $D_\eps^R$ with functional $\varphi\mapsto \norm{\varphi}_{L^2}^2$ defined over $L^2(\bo)$, we obtain
	\begin{equation*}
		\begin{split}
			&\norm{D_\eps^R(t)}_{L^2}^2+2\int_0^t\norm{D_\eps^R(s)}^2_{W^{1,2}}ds\\
			&=2\int_0^t(H_\eps((u_1(s))^+)-H_\eps((u_2(s))^+),D_\eps^R(s))ds\\
			&\qquad +\int_0^t\norm{g_2(v_1(s))-g_2(v_2(s))}^2_{\mathcal{L}_2^0}+2\int_0^t(g_2(v_1(s))-g_2(v_2(s)),D_\eps^R(s))dW_2(s).
		\end{split}
	\end{equation*}
	Since the functional $\varphi\mapsto\norm{\nabla\varphi}^2_{L^2}$ is twice-differentiable only in $D(A_2)$ and is not even well-defined on the whole space $L^2(\bo)$, instead of using the It\^o formula in the version of \cite[Theorem 4.32]{DaPrato}, we use the generalised version of the It\^o formula from \cite[Proposition 3.3]{Scarpa} (see also \cite[Proof of Lemma A.2]{Scarpa1}) with the process $(D_\eps^R (t))_{t\in [0,T_0]}$ to get
	\begin{equation*}
		\begin{split}
			&\norm{\nabla D_\eps^R(t)}_{L^2}^2+2\int_0^t\left(\norm{\Delta D_\eps^R(s)}^2_{L^2}+\norm{\nabla D_\eps^R(s)}_{L^2}^2\right)ds\\
			&= 2\int_0^t(H_\eps((u_1(s))^+)-H_\eps((u_2(s))^+),-\Delta D_\eps^R(s))ds\\
			&\qquad +\int_0^t\norm{\nabla\left(g_2(v_1(s))-g_2(v_2(s))\right)}^2_{\mathcal{L}_2^0}+2\int_0^t(g_2(v_1(s))-g_2(v_2(s)),-\Delta D_\eps^R(s))dW_2(s).
		\end{split}
	\end{equation*}
	Adding up these two identities, we get
	\begin{equation}\label{equa3.10}
		\begin{split}
			&\norm{ D_\eps^R(t)}_{W^{1,2}}^2+2\int_0^t\left(\norm{\Delta D_\eps^R(s)}^2_{L^2}+\norm{\nabla D_\eps^R(s)}_{L^2}^2+\norm{ D_\eps^R(s)}_{W^{1,2}}^2\right)ds\\
			&=2\int_0^t(H_\eps((u_1(s))^+)-H_\eps((u_2(s))^+),-\Delta D_\eps^R(s)+D_\eps^R(s))ds\\
			&\qquad +\int_0^t\norm{g_2(v_1(s))-g_2(v_2(s))}^2_{\mathcal{L}_2^1}+2\int_0^t(g_2(v_1(s))-g_2(v_2(s)),D_\eps^R(s))dW_2(s)\\
			&\qquad+2\int_0^t(g_2(v_1(s))-g_2(v_2(s)),-\Delta D_\eps^R(s))dW_2(s).
		\end{split}
	\end{equation}
	Let us now control the terms appearing on the right-hand side of \re{equa3.10}. By using 
	\re{equa3.2} and the Lipschitz property of the function $\max(0,\cdot)$, we get
	\begin{equation*}
		\begin{split}
			2(H_\eps((u_1)^+)-H_\eps((u_2)^+),-\Delta D_\eps^R)&\leq 2\norm{H_\eps((u_1)^+)-H_\eps((u_2)^+)}_{L^2}\norm{\Delta D_\eps^R}_{L^2}\\
			&\leq \norm{H_\eps((u_1)^+)-H_\eps((u_2)^+)}^2_{L^2}+\norm{\Delta D_\eps^R}_{L^2}^2\\
			&\leq \frac{16}{9}\norm{(u_1)^+-(u_2)^+}^2_{L^2}+\norm{\Delta D_\eps^R}_{L^2}^2\\
			&\leq \frac{16}{9}\norm{u_1-u_2}^2_{L^2}+\norm{\Delta D_\eps^R}_{L^2}^2,
		\end{split}
	\end{equation*}
	and
	\begin{equation*}
		\begin{split}
			2(H_\eps((u_1)^+)-H_\eps((u_2)^+),D_\eps^R)&\leq 2\norm{H_\eps((u_1)^+)-H_\eps((u_2)^+)}_{L^2}\norm{D_\eps^R}_{L^2}\\
			&\leq \frac{16}{9}\norm{u_1-u_2}^2_{L^2}+\norm{D_\eps^R}_{W^{1,2}}^2.
		\end{split}
	\end{equation*}
	Thanks to \re{equa2.4}, we have
	\begin{equation*}
		\norm{g_2(v_1)-g_2(v_2)}^2_{\mathcal{L}_2^1}\leq L^2_2 \norm{v_1-v_2}^2_{W^{1,2}}.
	\end{equation*}
	Replacing these inequalities in \re{equa3.10} after using an integration-by-parts on the stochastic term (note that the boundary integral term will vanish since $D^R_\eps\in D(A_2)$), we obtain
	\begin{equation*}
		\begin{split}
			&\norm{ D_\eps^R(t)}_{W^{1,2}}^2+\int_0^t\left(\norm{\Delta D_\eps^R(s)}^2_{L^2}+2\norm{\nabla D_\eps^R(s)}_{L^2}^2+\norm{ D_\eps^R(s)}_{W^{1,2}}^2\right)ds\\
			&\leq {\frac{32}{18}}\int_0^t\norm{u_1(s)-u_2(s)}^2_{L^2}ds +L^2_2\int_0^t\norm{v_1(s)-v_2(s)}^2_{W^{1,2}}ds\\
			&\quad+2\int_0^t(g_2(v_1(s))-g_2(v_2(s)),D_\eps^R(s))dW_2(s)+2\int_0^t(\nabla\left(g_2(v_1(s))-g_2(v_2(s))\right),\nabla D_\eps^R(s))dW_2(s),
		\end{split}
	\end{equation*}
	where we have used an integration-by-parts on the last term (note that since $D_\eps^R\in D(A_2)$, the boundary integral vanishes). Taking the square, the supremum in time over $[0,T_0]$ and the expectation, we deduce that
	\begin{equation*}
		\begin{split}
			\norm{ D_\eps^R}_{\by_{T_0}}^4
			&\leq \bk T_0^2\norm{u_1-u_2}^{{4}}_{\bx_{T_0}} +\bk T_0^2\norm{v_1-v_2}^{{4}}_{\by_{T_0}}\\
			&\quad+\bk\be\sup_{t\in[0,T_0]}\abs{\int_0^t(g_2(v_1(s))-g_2(v_2(s)),D_\eps^R(s))dW_2(s)}^{2}\\
			&\qquad+\bk\be\sup_{t\in[0,T_0]}\abs{\int_0^t(\nabla\left(g_2(v_1(s))-g_2(v_2(s))\right),\nabla D_\eps^R(s))dW_2(s)}^{2}\\
			&:= \bk T_0^2\norm{u_1-u_2}^{{4}}_{\bx_{T_0}} +\bk T_0^2\norm{v_1-v_2}^{{4}}_{\by_{T_0}}+I_1+I_2.
		\end{split}
	\end{equation*}
	The Burkholder–Davis–Gundy inequality and \re{equa2.4} imply
	\begin{equation*}
		\begin{split}
			I_1
			&\leq \bk\be {\int_0^{T_0}\norm{g_2(v_1(s))-g_2(v_2(s))}^2_{\mathcal{L}_2^0}\norm{D_\eps^R(s)}^2_{L^2}ds} \\
			&\leq \bk\be\sup_{s\in[0,T_0]}\norm{D_\eps^R(s)}^2_{W^{1,2}} {\int_0^{T_0}\norm{g_2(v_1(s))-g_2(v_2(s))}^2_{\mathcal{L}_2^1}ds} \\
			&\leq \frac{1}{4}\be\sup_{s\in[0,T_0]}\norm{D_\eps^R(s)}^4_{W^{1,2}}+\bk T_0^2\norm{v_1-v_2}^4_{\by_{T_0}},
		\end{split}
	\end{equation*}
	and
	\begin{equation*}
		\begin{split}
			I_2
			&\leq \bk\be\int_0^{T_0}\norm{\nabla\left(g_2(v_1(s))-g_2(v_2(s))\right)}^2_{\mathcal{L}_2^0}\norm{\nabla D_\eps^R(s)}^2_{L^2}ds\\
			&\leq \frac{1}{4}\be\sup_{s\in[0,T_0]}\norm{D_\eps^R(s)}^4_{W^{1,2}}+\bk T_0^2\norm{v_1-v_2}^4_{\by_{T_0}}.
		\end{split}
	\end{equation*}
	Combining these three inequalities, we get
	\begin{equation}\label{equa3.12}
		\begin{split}
			\norm{ D_\eps^R}_{\by_{T_0}}^4
			&\leq \bk T_0^2\norm{u_1-u_2}^4_{\bx_{T_0}} +\bk T_0^2\norm{v_1-v_2}^4_{\by_{T_0}}.
		\end{split}
	\end{equation}
	Let us now estimate $\norm{M_\eps^R}^4_{\bx_{T_0}}$. We start by pointing out that {the equation for $M^R_\eps$} is a stochastic variational equation with the Gelfand triple $W^{1,2}(\bo)\hookrightarrow L^2(\bo)\hookrightarrow(W^{1,2}(\bo))'$ instead of a stochastic differential equation in $L^2(\bo)$. Then, by applying the It\^o formula from \cite[Theorem 1.2]{Pardoux} with the functional $\psi\mapsto\norm{\psi}^2_{L^2}$, we obtain
	\begin{equation*}
		\begin{split}
			&\norm{M_\eps^R(t)}^2_{L^2}+2\int_0^t\norm{\nabla M_\eps^R(s)}_{L^2}^2ds\\
			&=2\chi \int_0^t\left(\theta_R\left(\norm{d^R_{1,\eps}(s)}_{W^{1,4}}\right)F_\eps((u_1(s))^+)\nabla d^R_{1,\eps}(s)-\theta_R\left(\norm{d^R_{2,\eps}(s)}_{W^{1,4}}\right)F_\eps((u_2(s))^+)\nabla d^R_{2,\eps}(s),\nabla M_\eps^R(s) \right)ds\\
			&\qquad+2\alpha\int_0^t( u_1(s)-u_2(s),M_\eps^R(s))ds-2\mu\int_0^t\left( [H_\eps((u_1(s))^+)]^2-[H_\eps((u_2(s))^+)]^2,M_\eps^R(s)\right)ds\\
			&\qquad+\int_0^t\norm{g_1(u^+_1(s))-g_1(u^+_2(s))}^2_{\mathcal{L}_2^0}+2\int_0^t(g_1(u^+_1(s))-g_1(u^+_2(s)),M_\eps^R(s))dW_1(s).
		\end{split}
	\end{equation*}
	In order to control the right-hand side terms of this identity, let us define
	\begin{equation*}
	{\Gamma:=\theta_R\bigl(\norm{d^R_{1,\eps}}_{W^{1,4}}\bigr)F_\eps((u_1)^+)\nabla d^R_{1,\eps}-\theta_R\bigl(\norm{d^R_{2,\eps}}_{W^{1,4}}\bigr)F_\eps((u_2)^+)\nabla d^R_{2,\eps}.}
	\end{equation*}
	Due to the definition of $\theta_R$, to be able to estimate $\Gamma$, we have to distinguish four cases.
	
	\textbf{Case 1: $\max\left(\norm{d^R_{1,\eps}}_{W^{1,4}}, \norm{d^R_{2,\eps}}_{W^{1,4}} \right)\leq 2R$.} Using the fact that $\Sup_{x\geq 0}\abs{\theta'_R(x)}\leq \frac{\bk}{R}$, \re{equa3.2}, and \re{equa3.3}, we get
	\begin{equation*}
		\begin{split}
			\norm{\Gamma}_{L^2}&\leq \norm{\left(\theta_R\left(\norm{d^R_{1,\eps}}_{W^{1,4}}\right)-\theta_R\left(\norm{d^R_{2,\eps}}_{W^{1,4}}\right)\right)F_\eps((u_1)_+)\nabla d_{1,\eps}^R(s)}_{L^2}\\
			&\qquad+\norm{\theta_R\left(\norm{d^R_{2,\eps}}_{W^{1,4}}\right)F_\eps((u_1)_+)\nabla (d_{1,\eps}^R-d_{2,\eps}^R)}_{L^2}\\
			&\qquad+\norm{\theta_R\left(\norm{d^R_{2,\eps}}_{W^{1,4}}\right)\left(F_\eps((u_1)_+)-F_\eps((u_2)_+)\right)\nabla d_{2,\eps}^R}_{L^2}\\
			&\leq \frac{\bk}{R\eps}\norm{d_{1,\eps}^R-d_{2,\eps}^R}_{W^{1,4}}\norm{\nabla d_{1,\eps}^R}_{L^2}+\frac{1}{\eps}\norm{\nabla(d_{1,\eps}^R-d^R_{2,\eps})}_{L^2}+\norm{u_1-u_2}_{L^4}\norm{\nabla d^R_{2,\eps}}_{L^4}\\
			&\leq \left(\frac{2\bk}{\eps}+\frac{1}{\eps}\right)\norm{d_{1,\eps}^R-d_{2,\eps}^R}_{W^{1,4}}+2R\norm{u_1-u_2}_{L^4}.
		\end{split}
	\end{equation*}
	\textbf{Case 2: ($\norm{d_{1,\eps}^R}_{W^{1,4}}\leq 2R$ and $\norm{d_{2,\eps}^R}_{W^{1,4}}> 2R$).} Here, we have
	\begin{equation*}
		\begin{split}
			\norm{\Gamma}_{L^2}&=\norm{\theta_R\left(\norm{d^R_{1,\eps}}_{W^{1,4}}\right)F_\eps((u_1)_+)\nabla d_{1,\eps}}_{L^2}\\
			&=\abs{\theta_R\left(\norm{d^R_{1,\eps}}_{W^{1,4}}\right)-\theta_R\left(\norm{d^R_{2,\eps}}_{W^{1,4}}\right)}\norm{F_\eps((u_1)_+)\nabla d_{1,\eps}^R}_{L^2}\\
			&\leq \frac{\bk}{\eps}\norm{d_{1,\eps}^R-d_{2,\eps}^R}_{W^{1,4}}.
		\end{split}
	\end{equation*}
	\textbf{Case 3: ($\norm{d_{1,\eps}^R}_{W^{1,4}}> 2R$ and $\norm{d_{2,\eps}^R}_{W^{1,4}}\leq 2R$).} We have
	\begin{equation*}
		\norm{\Gamma}_{L^2}=\norm{\theta_R\left(\norm{d^R_{2,\eps}}_{W^{1,4}}\right)F_\eps((u_2)_+)\nabla d_{2,\eps}}_{L^2}
		\leq \frac{\bk}{\eps}\norm{d_{1,\eps}^R-d_{2,\eps}^R}_{W^{1,4}}.
	\end{equation*}
	\textbf{Case 4: $\min\left(\norm{d_{1,\eps}^R}_{W^{1,4}}, \norm{d_{2,\eps}^R}_{W^{1,4}} \right)> 2R$.} Here we obviously have $\Gamma=0$.
	
	It then follows from these cases that
	\begin{equation*}
		\norm{\Gamma}_{L^2}\leq \bk_{\eps,R}\norm{d_{1,\eps}^R-d_{2,\eps}^R}_{W^{1,4}}+2R\norm{u_1-u_2}_{L^4}.
	\end{equation*}
	The H\"older inequality and the Sobolev embedding $H^2(\bo)\hookrightarrow  W^{1,4}(\bo)$, together with \re{equa2.3} and the Young inequality yield
	\begin{equation*}
		\begin{split}
			&2\chi(\theta_R\left(\norm{d^R_{1,\eps}}_{W^{1,4}}\right)F_\eps((u_1)^+)\nabla d^R_{1,\eps}-\theta_R\left(\norm{d^R_{2,\eps}}_{W^{1,4}}\right)F_\eps((u_2)^+)\nabla d^R_{2,\eps},\nabla M_\eps^R )\\
			&\leq 2\chi \norm{\Gamma}_{L^2}\norm{\nabla M_\eps^R}_{L^2}\\
			&\leq \bk_{\eps,R}\norm{d_{1,\eps}^R-d_{2,\eps}^R}_{W^{2,2}}^2+16R^2\chi^2\bk\norm{u_1-u_2}_{L^2}\norm{u_1-u_2}_{W^{1,2}}+\norm{\nabla M_\eps^R}_{L^2}^2\\
			&\leq \bk_{\eps,R}\norm{d_{1,\eps}^R-d_{2,\eps}^R}_{W^{2,2}}^2+R^6\chi^2\bk\norm{u_1-u_2}_{L^2}^2+\frac{\chi^2}{R^2}\norm{u_1-u_2}^2_{W^{1,2}}+\norm{\nabla M_\eps^R}_{L^2}^2,
		\end{split}
	\end{equation*}
	and
	\begin{equation*}
		2\alpha( u_1-u_2,M_\eps^R)\leq 2\alpha\norm{u_1-u_2}_{L^2}\norm{M_\eps^R}_{L^2}.
	\end{equation*}
	Thanks to \re{equa3.4}, \re{equa3.2}, and the Lipschitz property of the function $(\cdot)^+$, we get
	\begin{equation*}
		\begin{split}
			&-2\mu( [H_\eps((u_1)^+)]^2-[H_\eps((u_2)^+)]^2,M_\eps^R)\\
			&=-2\mu( \left[[H_\eps((u_1)^+)]+[H_\eps((u_2)^+)]\right]\left[[H_\eps((u_1)^+)]-[H_\eps((u_2)^+)]\right],M_\eps^R)\\
			&\leq \bk_\eps\abs{( [H_\eps((u_1)^+)]-[H_\eps((u_2)^+)],M_\eps^R)}\\
			&\leq 7\bk_\eps\norm{u_1-u_2}_{L^2}\norm{M_\eps^R}_{L^2}.
		\end{split}
	\end{equation*}
	Using \re{equa2.4}, we get
	\begin{equation*}
		\norm{g_1(u_1^+)-g_1(u_2^+)}^2_{\mathcal{L}_2^0}\leq L_1^2 \norm{u^+_1-u^+_2}_{L^2}^2\leq  L_1^2 \norm{u_1-u_2}_{L^2}^2.
	\end{equation*}
	Combining these inequalities together with the fact that $(a+b)^2\leq 2(a^2+b^2)$ for any real numbers $a, b>0$, we obtain
	\begin{equation}\label{equa3.13}
		\begin{split}
			&\be\sup_{t\in[0,T_0]}\norm{M_\eps^R(t)}^4_{L^2}+\be\left(\int_0^{T_0}\norm{\nabla M_\eps^R(s)}_{L^2}^2ds\right)^{2}\\
			&\leq\bk_{\eps,R} \be\left(\int_0^{T_0}\norm{D_\eps^R(s)}_{W^{2,2}}^2ds\right)^2+\bk_{R}\be\left(\int_0^{T_0}\norm{u_1(s)-u_2(s)}_{L^2}^2ds\right)^{2}\\
			&\quad+\frac{2\chi^4}{R^{4}}\be\left(\int_0^{T_0}\norm{u_1(s)-u_2(s)}^2_{W^{1,2}}ds\right)^2+\bk\be\left(\int_0^{T_0} \norm{u_1(s)-u_2(s)}_{L^2}^2ds\right)^2\\
			&\quad+\bk_{\eps}\be\left(\int_0^{T_0}\norm{u_1(s)-u_2(s)}_{L^2}\norm{M_\eps^R(s)}_{L^2}ds\right)^2\\
			&\quad+\bk\be\sup_{t\in[0,T_0]}\abs{\int_0^t(g_1(u^+_1(s))-g_1(u^+_2(s)),M_\eps^R(s))dW_1(s)}^2\\
			&:=J_1+J_2+J_3+J_4+J_5+J_6.
		\end{split}
	\end{equation}
	Let us note that
	\begin{equation*}
		\begin{split}
			J_5
			&\leq\bk_{\eps}\be \sup_{s\in[0,T_0]}\norm{M_\eps^R(s)}^2_{L^2}\left(\int_0^{T_0}\norm{u_1(s)-u_2(s)}_{L^2}ds\right)^2\\
			&\leq \frac{1}{4}\be\sup_{s\in[0,T_0]}\norm{M_\eps^R(s)}_{L^2}^4+\bk_{\eps}\be\left(\int_0^{T_0}\norm{u_1(s)-u_2(s)}_{L^2}ds \right)^4\\
			&\leq \frac{1}{4}\be\sup_{s\in[0,T{_0}]}\norm{M_\eps^R(s)}_{L^2}^{{4}}+ T_0^2\bk_{\eps}\be\left(\int_0^{T_0}\norm{u_1(s)-u_2(s)}_{L^2}^2ds\right)^2.
		\end{split}
	\end{equation*}
	The Burkholder–Davis–Gundy inequality along with \re{equa2.4} give
	\begin{equation*}
		\begin{split}
			J_6
			&\leq \bk\be{\int_0^{T_0}\norm{g_1(u^+_1(s))-g_1(u^+_2(s))}^2_{\mathcal{L}_2^0}\norm{M_\eps^R(s)}^2_{L^2}ds}\\
			&\leq \frac{1}{4}\be\sup_{s\in[0,T_0]}\norm{M_\eps^R(s)}^4_{L^2}+\bk\be\left(\int_0^{T_0}\norm{u_1(s)-u_2(s)}^2_{L^2}ds\right)^{2}.
		\end{split}
	\end{equation*}
	It then follows from these two inequalities, together with \re{equa3.13} and \re{equa3.12}, that
	\begin{equation*}
		\begin{split}
			&\frac{1}{2}\be\sup_{t\in[0,T_0]}\norm{M_\eps^R(t)}^4_{L^2}+\be\left(\int_0^{T_0}\norm{\nabla M_\eps^R(s)}_{L^2}^2ds\right)^2\\
			&\leq\bk_{\eps,R} \norm{D_\eps^R}^4_{\by_{T_0}}+\bk_{R} T_0^2\norm{u_1-u_2}^4_{\bx_{T_0}}+\frac{2\chi^4}{R^{4}}\norm{u_1-u_2}^4_{\bx_{T_0}}\\
			&\quad+\bk_{\eps}T_0^4\norm{u_1-u_2}^4_{\bx_{T_0}}\\
			&\leq \left[\bk_{\eps,R,\chi} (T_0^2+T_0^4)+\frac{2\chi^4}{R^{4}}\right]\norm{u_1-u_2}^4_{\bx_{T_0}}+\bk T_0^2\norm{v_1-v_2}^4_{\by_{T_0}}\\
			&\leq \left[\bk_{\eps,R,\chi} (T_0^2+T_0^4)+\frac{2\chi^4}{R^{4}}+\bk T_0^2\right]\left(\norm{u_1-u_2}^4_{\bx_{T_0}}+\norm{v_1-v_2}^4_{\by_{T_0}}\right).
		\end{split}
	\end{equation*}
	By remarking that
	\begin{equation*}
		\norm{M_\eps^R}^4_{\bx_{T_0}}\leq (1+2T_0^2)\be\sup_{t\in[0,T_0]}\norm{M_\eps^R(t)}^4_{L^2}+2\be\left(\int_0^{T_0}\norm{\nabla M_\eps^R(s)}_{L^2}^2ds\right)^{{2}},
	\end{equation*}
	we infer from these two inequalities that
	\begin{equation*}
		\begin{split}
			\norm{M_\eps^R}^4_{\bx_{T_0}}&\leq4(1+T^2_0)\left(\bk_{\eps,R,\chi} (T^2_0+T^4_0)+\frac{2\chi^4}{R^{4}}+\bk T_0^2\right)\left(\norm{u_1-u_2}^4_{\bx_{T_0}}+\norm{v_1-v_2}^4_{\by_{T_0}}\right)\\
			&\leq\left(\bk_{\eps,R,\chi} (T_0^2+T_0^4+T_0^6)+\frac{8\chi^4}{R^{4}}\right)\left(\norm{u_1-u_2}^4_{\bx_{T_0}}+\norm{v_1-v_2}^4_{\by_{T_0}}\right).
		\end{split}
	\end{equation*}
	Combining this, with \re{equa3.12}, we get
	\begin{equation*}
		\norm{M_\eps^R}^4_{\bx_{T_0}}+	\norm{D_\eps^R}^4_{\by_{T_0}}\leq\left(\bk_{\eps,R,\chi} (T_0^2+T_0^4+T_0^6)+\frac{8\chi^4}{R^{4}}\right)\left(\norm{u_1-u_2}^4_{\bx_{T_0}}+\norm{v_1-v_2}^4_{\by_{T_0}}\right),
	\end{equation*}
	which implies that
	\begin{equation*}
		\begin{split}
			&\norm{\varUpsilon(u_1,v_1)-\varUpsilon(u_2,v_2)}^4_{\bx_{T_0}\times\by_{T_0}}\\
			&\leq\left(\bk_{\eps,R,\chi} T^2_0(1+T_0^2+T_0^4)+\frac{8\chi^4}{R^{4}}\right)\norm{(u_1,v_1)-(u_2,v_2)}^4_{\bx_{T_0}\times\by_{T_0}}.
		\end{split}
	\end{equation*}
	Let us choose
	\begin{equation*}
		T_0=T_{\eps,R}:=\min\left(\frac{1}{4},\frac{\sqrt{2}}{\sqrt{9\bk_{\eps,R,\chi}}}\right).
	\end{equation*}
	With this choice, we have $T_0\in (0,1)$ and consequently, $1+T_0^2+T_0^4\leq 1+2T_0^2\leq \frac{9}{8}$. This implies that $\bk_{\eps,R,\chi} T_0^2(1+T_0^2+T_0^4)\leq \frac{9\bk_{\eps,R,\chi}}{8}T_0^2\leq \frac{1}{4}$.
	On the other hand, choosing arbitrary $R\geq 2\chi$, it holds that $\frac{8\chi^4}{R^4}\leq \frac{1}{2}$ and therefore
	\begin{equation*}
		\norm{\varUpsilon(u_1,v_1)-\varUpsilon(u_2,v_2)}_{\bx_{T_0}\times\by_{T_0}}
		\leq\left(\frac{3}{4}\right)^\frac{1}{4}\norm{(u_1,v_1)-(u_2,v_2)}_{\bx_{T_0}\times\by_{T_0}}.
	\end{equation*}
	Then $\varUpsilon$ is a contraction mapping on the space $\bx_{T_{\eps,R}}\times\by_{T_{\eps,R}}$. By applying the Banach fixed-point theorem, there exists a unique $\mathbb{F}_{T_{\eps,R}}$-adapted process $(n_{\eps,1}^R,c_{\eps,1}^R)\in\bx_{T_{\eps,R}}\times\by_{T_{\eps,R}}$ such that $(n_{\eps,1}^R,c_{\eps,1}^R)$ satisfies \re{equa3.5} for all $t\in [0,T_{\eps,R}]$. Thanks to the fact that $(n_{\eps,1}^R,c_{\eps,1}^R)\in\bx_{T_{\eps,R}}\times\by_{T_{\eps,R}}$, we get
	\begin{equation}\label{equa3.14}
		\begin{split}
			&\be\norm{n_{\eps,1}^R(T_{\eps,R})}_{L^2}^4\leq\be\sup_{t\in[0,T_{\eps,R}]}\norm{n_{\eps,1}^R(t)}_{L^2}^4<\infty,\\
			&\be\norm{c_{\eps,1}^R(T_{\eps,R})}_{W^{1,2}}^4\leq\be\sup_{t\in[0,T_{\eps,R}]}\norm{c_{\eps,1}^R(t)}_{W^{1,2}}^4<\infty.
		\end{split}
	\end{equation}
	It is important to note that the condition $\be \norm{v_0}^4_{W^{1,2}}<\infty$ is sufficient to construct the process $(n_{\eps,1}^R,c_{\eps,1}^R)\in\bx_{T_{\eps,R}}\times\by_{T_{\eps,R}}$ (see, the inequality \re{equa2.15}).
	
	Let us now define the spaces
	\begin{equation*}
		\begin{split}
			&\bx_{T_{\eps,R}}^2:=L^4_{\mathbb{F}^2_T}(\Omega; C([T_{\eps,R},2T_{\eps,R}]; L^2(\bo)))\cap L^4_{\mathbb{F}^2_T}(\Omega; L^2(T_{\eps,R},2T_{\eps,R}; W^{1,2}(\bo))),\\
			&\by^2_{T_{\eps,R}}:=L^4_{\mathbb{F}^2_T}(\Omega; C([T_{\eps,R},2T_{\eps,R}]; W^{1,2}(\bo)))\cap L^4_{\mathbb{F}^2_T}(\Omega; L^2(T_{\eps,R},2T_{\eps,R}; D(A_2))),
		\end{split}
	\end{equation*}
	with $\mathbb{F}^2_T=(\mathcal{F}_t)_{t\in[T_{\eps,R},2T_{\eps,R}]}$.
	By using \re{equa3.14}, we can follow the same methodology of construction of $(n_{\eps,1}^R,c_{\eps,1}^R)$ and construct a unique $\mathbb{F}^2_T$-adapted process $(n_{\eps,2}^R,c_{\eps,2}^R)\in\bx^2_{T_{\eps,R}}\times\by^2_{T_{\eps,R}}$ satisfying for all $t\in [T_{\eps,R},2T_{\eps,R}]$,
	\begin{equation*}
		\begin{split}
			&	(n^R_{\eps,2}(t),\varphi)+\int_{T_{\eps,R}}^t( \nabla n^R_{\eps,2}(s),\nabla\varphi)ds-\chi \int_{T_{\eps,R}}^t\theta_R(\norm{c^R_{\eps,2}(s)}_{W^{1,4}})(F_{\eps}((n^R_{\eps,2}(s))^+)\nabla c^R_{\eps,2}(s),\nabla \varphi )ds\\
			&=({n^R_{\eps,1}}(T_{\eps,R}),\varphi)+\int_{T_{\eps,R}}^t(\alpha n^R_{\eps,2}(s)-\mu [H_{\eps}((n^R_{\eps,2}(s))^+)]^2,\varphi)ds\\
			&	\hspace{0.9cm}+\int_{T_{\eps,R}}^t(g_1((n^R_{\eps,2}(s))^+),\varphi)dW_1(s),\quad  \varphi\in W^{1,2}(\bo),\quad  \mathbb{P}\text{-a.s.}, \\
			&	(c^R_{\eps,2}(t),\psi)+\int_{T_{\eps,R}}^t(A_2 c^R_{\eps,2}(s),\psi)ds\\
			&=({c^R_{\eps,1}}(T_{\eps,R}),\psi)+\int_{T_{\eps,R}}^t(H_{\eps}((n^R_{\eps,2}(s))^+),\psi)ds+\int_{T_{\eps,R}}^t(g_2(c^R_{\eps,2}(s)),\psi)dW_2(s),\quad \psi\in L^2(\bo)\quad\mathbb{P}\text{-a.s.},
		\end{split}
	\end{equation*}
	and
	\begin{equation*}
		\begin{split}
			&\be\norm{n_{\eps,2}^R(2T_{\eps,R})}_{L^2}^4\leq\be\sup_{t\in[T_{\eps,R},2T_{\eps,R}]}\norm{n_{\eps,2}^R(t)}_{L^2}^4<\infty,\\
			&\be\norm{c_{\eps,2}^R(2T_{\eps,R})}_{W^{1,2}}^4\leq\be\sup_{t\in[T_{\eps,R},2T_{\eps,R}]}\norm{c_{\eps,2}^R(t)}_{W^{1,2}}^4<\infty.
		\end{split}
	\end{equation*}
	Let $1\leq \ell\leq {N_{\eps,R}:=\lceil T/T_{\eps,R}\rceil}$ be an integer, where by $\lceil T/T_{\eps,R}\rceil$ we mean the integer part of $T/T_{\eps,R}$. 
Let us define
	\begin{equation*}
		\begin{split}
			&\bx_{T_{\eps,R}}^\ell:=L^4_{\mathbb{F}_{T}^\ell}(\Omega; C([(\ell-1)T_{\eps,R},\ell T_{\eps,R}]; L^2(\bo)))\cap L^4_{\mathbb{F}^\ell_T}(\Omega ; L^2((\ell-1)T_{\eps,R},\ell T_{\eps,R}; W^{1,2}(\bo))),\\
			&\by^\ell_{T_{\eps,R}}:=L^4_{\mathbb{F}^\ell_T}(\Omega ; C([(\ell-1)T_{\eps,R},\ell T_{\eps,R}]; W^{1,2}(\bo)))\cap L^4_{\mathbb{F}^\ell_T}(\Omega; L^2((\ell-1)T_{\eps,R},\ell T_{\eps,R}; D(A_2))),
		\end{split}
	\end{equation*}
	with $\mathbb{F}^\ell_T=(\mathcal{F}_t)_{t\in[(\ell-1)T_{\eps,R},(\ell)T_{\eps,R}]}$. Repeating the above arguments, we can recursively construct a unique $\mathbb{F}^\ell_T$-adapted process $(n_{\eps,\ell}^R,c_{\eps,\ell}^R)\in\bx^\ell_{T_{\eps,R}}\times\by^\ell_{T_{\eps,R}}$ satisfying for any $t\in [(\ell-1)T_{\eps,R},\ell T_{\eps,R}]$,
	\begin{equation*}
		\begin{split}
			&	(n^R_{\eps,\ell}(t),\varphi)+\int_{(\ell-1)T_{\eps,R}}^t( \nabla n^R_{\eps,\ell}(s),\nabla\varphi)ds-\chi \int_{(\ell-1)T_{\eps,R}}^t\theta_R(\norm{c^R_{\eps,\ell}(s)}_{W^{1,4}})(F_{\eps}((n^R_{\eps,\ell}(s))^+)\nabla c^R_{\eps,\ell}(s),\nabla \varphi )ds\\
			&=({n^R_{\eps,\ell-1}((\ell-1)T_{\eps,R})},\varphi)+\int_{(\ell-1)T_{\eps,R}}^t(\alpha n^R_{\eps,\ell}(s)-\mu [H_{\eps}((n^R_{\eps,\ell}(s))^+)]^2,\varphi)ds\\
			&\qquad+\int_{(\ell-1)T_{\eps,R}}^t(g_1({(}n^R_{\eps,\ell}(s){)^+}),\varphi)dW_1(s),\quad\varphi\in W^{1,2}(\bo),\quad  \mathbb{P}\text{-a.s.}, \\
			&	(c^R_{\eps,\ell}(t),\psi)+\int_{(\ell-1)T_{\eps,R}}^t(A{_2} c^R_{\eps,\ell}(s),\psi)ds-\int_{(\ell-1)T_{\eps,R}}^t(H_{\eps}((n^R_{\eps,\ell}(s))^+),\psi)ds\\
			&=({c^R_{\eps,\ell-1}((\ell-1)T_{\eps,R})},\psi)+\int_{(\ell-1)T_{\eps,R}}^t(g_2(c^R_{\eps,\ell}(s)),\psi)dW_2(s),\quad\psi\in L^2(\bo),\quad\mathbb{P}\text{-a.s.},
		\end{split}
	\end{equation*}
	and
	\begin{equation*}
		\begin{split}
			&\be\norm{n_{\eps,\ell}^R(\ell T_{\eps,R})}_{L^2}^4\leq\be\sup_{t\in[(\ell-1)T_{\eps,R},\ell T_{\eps,R}]}\norm{n_{\eps,\ell}^R(t)}_{L^2}^4<\infty,\\
			&\be\norm{c_{\eps,\ell}^R(\ell T_{\eps,R})}_{W^{1,2}}^4\leq\be\sup_{t\in[(\ell-1)T_{\eps,R},\ell T_{\eps,R}]}\norm{c_{\eps,\ell}^R(t)}_{W^{1,2}}^4<\infty.
		\end{split}
	\end{equation*}
	Let us now construct our process $((n_\eps^R(t),c_\eps^R(t)))_{t\in [0,T]}$ by setting
	\begin{equation*}
		(n_\eps^R(t),c_\eps^R(t)):=(n_{\eps,\ell}^R(t),c_{\eps,\ell}^R(t))\quad \text{ if }\quad t\in [(\ell-1)T_{\eps,R},\ell T_{\eps,R}], \qquad \ell\in \{1,2,\dots, {N_{\eps,R}}\}.
	\end{equation*}
	It is not too difficult to show that $(n_\eps^R(t),c_\eps^R(t))_{t\in [0,T]}$ is an $\mathbb{F}_T$-adapted process with $\mathbb{P}$-a.s. continuous sample paths in $L^2(\bo)\times W^{1,2}(\bo)$ satisfying \re{equa3.5} and
	\begin{equation}\label{equa2.19}
		\begin{split}
			&n_\eps^R\in L^4_{\mathbb{F}_T}(\Omega ; C([0,T]; L^2(\bo)))\cap L^4_{\mathbb{F}_T}(\Omega; L^2(0,T; W^{1,2}(\bo))),\\
			&c_\eps^R\in L^4_{\mathbb{F}_T}(\Omega ; C([0,T]; W^{1,2}(\bo)))\cap L^4_{\mathbb{F}_T}(\Omega ; L^2(0,T; D(A_2))).
		\end{split}
	\end{equation}
	It remains to prove \re{equa3.6*}. Thanks to \re{equa2.19}, it suffices to prove that
	\begin{equation}\label{equa2.20}
		c_\eps^R\in L^4_{\mathbb{F}_T}(\Omega ; C([0,T]; B_{4,4,\mathbf{n}}^\frac{3}{2}(\bo)))\cap L^4_{\mathbb{F}_T}(\Omega ; L^4(0,T; D(A_4))).
	\end{equation}
	To achieve this, let us consider the stochastic parabolic evolution equation that consists in finding an $\mathbb{F}_T$-adapted, $B_{4,4,\mathbf{n}}^\frac{3}{2}(\bo)$-valued continuous process $(u(t))_{t\in [0,T]}$ such that
	\begin{equation}\label{equa2.21}
		u\in L^4_{\mathbb{F}_T}(\Omega ; C([0,T]; B_{4,4,\mathbf{n}}^\frac{3}{2}(\bo)))\cap 	 L^4_{\mathbb{F}_T}(\Omega ; L^4(0,T; D(A_4))),
	\end{equation}
	and for any $t\in [0,T]$,
	\begin{equation*}
		u(t)+\int_0^tA_4u(s)ds=v_0+\int_0^tH_\eps((n_\eps^R(s))^+)ds+\int_0^tg_2(u(s))dW_2(s),\ \text{ in } L^4(\bo) \quad \mathbb{P}\text{-a.s.},
	\end{equation*}
	meaning that for any $\varphi\in L^\frac{4}{3}(\bo)$,
	\begin{equation}\label{equa2.22}
		\begin{split}
			&	\int_\bo u(t,x)\varphi(x)dx+\int_0^t\int_\bo A_4u(s,x)\varphi(x)dxds-\int_0^t\int_\bo H_\eps((n_\eps^R(s,x))^+)\varphi(x)dxds\\
			&=\int_\bo v_0(x)\varphi(x)dx+\int_0^t\int_\bo g_2(u(s,x))\varphi(x)dxdW_2(s),  \quad \mathbb{P}\text{-a.s.}.
		\end{split}
	\end{equation}
	It is sufficient to prove that such $u$ really exists. In this way, one can take in particular $\varphi\in L^2(\bo)\hookrightarrow L^\frac{4}{3}(\bo)$ in \re{equa2.22} and derive that $u$ also satisfies the second equation of \re{equa3.5} and by pathwise uniqueness of the solution to the second equation of \re{equa3.5} (which is not difficult to prove by using \re{equa2.4}), we infer that $u=c_\eps^R$ and \re{equa2.20} follows from \re{equa2.21}.
	Since there exist constants $\rho\geq 0$ and $\varsigma\in [0,\frac{\pi}{2})$ such that $A_4+\rho$ has a bounded $H^\infty(\sum_\varsigma)$-calculus (according to \cite[Item 4 of Example 3.2]{Neerven}), thanks to \re{equa2.1}, the existence of the unique process $u$ satisfying \re{equa2.21} and \re{equa2.22} comes from \cite[Theorem 4.5]{Neerven}. In fact, let us show that the assumptions of \cite[Theorem 4.5]{Neerven} are satisfied. Setting $X_0:=L^4(\bo)$ and $X_1:=D(A_4)$, we derive from \cite[Section 4.3.3]{Neerven} that $X_\frac{1}{2}:=[X_0,X_1]_\frac{1}{2}=D(A_4^\frac{1}{2})$ and $X_{\frac{3}{4},4}:=(X_0,X_1)_{1-\frac{1}{4},4}=B^\frac{3}{2}_{4,4,\mathbf{n}}(\bo)$. Owing to \re{equa2.1}, the (HX), (HA), and (Hu$_0$) of Assumptions (H) are naturally satisfied. For item (HF), we use \re{equa3.4} to observe that for $\omega\in \Omega$, $\Int_0^T\norm{H_\eps((n_\eps^R(\omega,s))^+)}_{L^4}ds\leq \bk_\eps T<\infty$. As far as item (HB) is concerned, we use \re{equa2.4**} to infer that for any $v_1,v_2$,
	\begin{equation*}
		\norm{g_2(v_1)-g_2(v_2)}_{  \gamma(L^2,X_\frac{1}{2})}\leq L_B\norm{v_1-v_2}_{X_1} \text{ with } L_B:=\frac{1}{(\bk_4^ \Diamond+1)},
	\end{equation*}
	and by the embedding $D(A_4)\hookrightarrow D(A_4^\frac{1}{2})$, we get
	\begin{equation*}
		\norm{g_2(v_1)}_{  \gamma(L^2,X_\frac{1}{2})}\leq L_3\norm{v_1}_{D(A_4^\frac{1}{2})}\leq  \bk\norm{v_1}_{X_1}  \leq   \bk(1+\norm{v_1}_{X_1} ).
	\end{equation*}
	In addition, by using \re{equa3.4} once more, we get $\be\Int_0^T\norm{H_\eps((n_\eps^R(s))^+)}^4_{L^4}ds\leq \bk_\eps T<\infty$. 
Since $L_B\bk_4^ \Diamond=\frac{\bk_4^ \Diamond}{(\bk_4^ \Diamond+1)}<1$, the assumptions of \cite[Theorem 4.5]{Neerven} are therefore satisfied, completing the proof of Claim \ref{lem3.1}.

Let us denote by $\bk^\flat$ the constant coming from the embedding $B^\frac{3}{2}_{4,4}(\bo)\hookrightarrow W^{1,4}(\bo)$ (see \re{equa1.5}) and consider the family of $\mathbb{F}_T$-adapted, continuous $L^2(\bo)\times B_{4,4,\mathbf{n}}^\frac{3}{2}(\bo)$-valued processes $(((n^R_\eps(t),c^R_\eps(t)))_{t\in [0,T]})_{R\geq 2\chi}$ constructed in Claim \ref{lem3.1}. It follows from the well-known result \cite[Chapter I, Proposition 4.5]{Revuz} that $(\tau^R_\eps)_{R\geq 2\chi}$, where
\begin{equation*}
	\begin{split}
		\tau_\eps^R:&=
		\inf\left\{ t\in [0,T]: \Sup_{s\in[0,t]}\norm{n_\eps^R(s)}^2_{L^2}+\Sup_{s\in[0,t]}\norm{c_\eps^R(s)}^2_{B^\frac{3}{2}_{4,4}}\right.\\
		&\hspace{4cm}\left.+\int_0^{t}\left(\norm{n_\eps^R(s)}^2_{W^{1,2}}+\norm{c_\eps^R(s)}^2_{W^{2,2}}\right)ds\geq \frac{R^2}{(\bk^\flat)^2} \right\},
	\end{split}
\end{equation*}
with the convention $\inf\{\}=T$,
is a family of stopping times. In addition, $(((n^R_\eps(t),c^R_\eps(t)))_{t\in [0,T]})_{R\geq 2\chi}$ is a family of $\mathbb{F}_T$-progressively measurable processes since any $\mathbb{F}_T$-adapted and continuous process is progressively measurable. To construct a local solution, let us distinguish two cases:

\textbf{Case 1: There exists a constant $\bk^\ast>0$ such that $\norm{u_0}_{L^2}^2+\norm{v_0}^2_{B_{4,4}^\frac{3}{2}}< \bk^\ast$ almost surely.}
By choosing $R=R^\ast:={\max(2\chi,\bk^\flat\sqrt{\bk^\ast})}$, we see that $\tau_\eps^{R^\ast}>0$ $\mathbb{P}$-almost surely, and according to the definition of cut-off functions, for all $t\in [0,T\wedge\tau_\eps^{R^\ast}]$, we have $\norm{c_\eps^{R^\ast}(t)}_{W^{1,4}}\leq \bk^\flat \norm{c_\eps^{R^\ast}(t)}_{B^\frac{3}{2}_{4,4}}\leq\frac{\bk^\flat R^\ast}{ \bk^\flat}= R^\ast$, which means that $\theta_{R^\ast}\left(\norm{c_\eps^{R^\ast}}_{W^{1,4}}\right)=1$ on $[0,T\wedge\tau_\eps^{R^\ast}]$, and therefore, $(n^{R^\ast}_\eps,c^{R^\ast}_\eps, \tau^{R^\ast}_\eps)$ satisfies the system \re{equa3.5*}. In addition, by the continuity of the process $((n^{R^\ast}_\eps,c^{R^\ast}_\eps))_{t\in [0,T]}$, the embedding $W^{2,4}(\bo)\hookrightarrow W^{2,2}(\bo)$ and \re{equa3.6*}, we derive that
\begin{equation*}
	(n^{R^\ast}_\eps(\cdot\wedge\tau_\eps^{R^\ast}),c^{R^\ast}_\eps(\cdot\wedge\tau_\eps^{R^\ast}))\in C([0,T]; L^2(\bo)\times B^\frac{3}{2}_{4,4,\mathbf{n}}(\bo))\quad \mathbb{P}\text{-a.s.}.
\end{equation*}
and
\begin{equation*}
	\begin{split}
		&\be \sup_{s\in[0,T\wedge\tau_\eps^{R^\ast}]}\norm{n_\eps^{R^\ast}(s)}^4_{L^2}+\be\left(\int_0^{T\wedge\tau_\eps^{R^\ast}}\norm{n_\eps^{R^\ast}(s)}^2_{W^{1,2}}ds\right)^2<\infty,\\
		&\be \sup_{s\in[0,T\wedge\tau_\eps^{R^\ast}]}\norm{c_\eps^{R^\ast}(s)}^4_{B^\frac{3}{2}_{4,4}}+\be\left(\int_0^{T\wedge\tau_\eps^{R^\ast}}\norm{c_\eps^{R^\ast}(s)}^2_{W^{2,2}}ds\right)^2<\infty,
	\end{split}
\end{equation*}
and hence, according to Definition \ref{defi3.1}, $(n^{R^\ast}_\eps,c^{R^\ast}_\eps, \tau^{R^\ast}_\eps)$ is a local solution to the system \re{equa3.1}.

\textbf{Case 2: General case where $\be [\norm{u_0}^4_{{L^2}}+\norm{v_0}_{B_{4,4}^\frac{3}{2}}^4]<\infty$.} Here, the local solution is constructed by a localisation argument from \cite{Glatt-Holtz1}. For any integer $\kappa\geq 0$, let $\Omega_\kappa:=\{\omega\in \Omega: \kappa\leq \norm{u_0}_{L^2}+\norm{v_0}_{B_{4,4}^\frac{3}{2}}<\kappa+1\}$ and decompose $u_0$ and $v_0$ as
\begin{equation*}
	u_0=\sum_{\kappa=0}^\infty u_{0,\kappa},\quad v_0=\sum_{\kappa=0}^\infty v_{0,\kappa},\text{ where } u_{0,\kappa}:=u_0\mathbf{1}_{\Omega_\kappa}\text{ and } v_{0,\kappa}:=v_0\mathbf{1}_{\Omega_\kappa}.
\end{equation*}
Then for each $\kappa\geq 0$, one can choose $R={\max(2\chi,(\kappa+1)\bk^\flat)}$ 
 and construct as in Case 1 a local solution $(n_{\eps,\kappa},c_{\eps,\kappa}, \tau_{\eps,\kappa})$ to \re{equa3.1} with initial data $(u_{0,\kappa},v_{0,\kappa})$. Let $M>1$ be fixed. For any integer $\kappa\geq 0$, let us define the stopping time
\begin{equation*}
	\begin{split}
		\rho_{\eps, \kappa}&:=\inf\left\{ t\in [0,T]: \Sup_{s\in[0,t]}\norm{n_{\eps,\kappa}(s\wedge\tau_{\eps,\kappa})}^2_{L^2}+\Sup_{s\in[0,t]}\norm{c_{\eps,\kappa}(s\wedge\tau_{\eps,\kappa})}^2_{B^\frac{3}{2}_{4,4}}\right.\\
		&\left.\hspace{2.6cm}+\int_0^{t\wedge\tau_{\eps,\kappa}}\left(\norm{n_{\eps,\kappa}(s)}^2_{W^{1,2}}+\norm{c_{\eps,\kappa}(s)}^2_{W^{2,2}}\right)ds\geq M+\norm{{u_0}}_{L^2}^2+\norm{{v_0}}^2_{B^\frac{3}{2}_{4,4}} \right\}\wedge\tau_{\eps,\kappa}.
	\end{split}
\end{equation*}
Clearly, since $\tau_{\eps,\kappa}$ is a strictly positive stopping time, by the continuity of $(n_{\eps,\kappa}(\cdot\wedge\tau_{\eps,\kappa}),c_{\eps,\kappa}(\cdot\wedge\tau_{\eps,\kappa}))$ in $L^2(\bo)\times B^\frac{3}{2}_{4,4}(\bo)$, we are sure that $\rho_{\eps, \kappa}$ is a stopping time such that $\rho_{\eps, \kappa}>0$ almost surely and $(n_{\eps,\kappa},c_{\eps,\kappa}, \rho_{\eps,\kappa})$ is a local solution to \re{equa3.1} with initial data $(u_{0,\kappa},v_{0,\kappa})$.
We now define a triple $(n_\eps,c_\eps, \tau_\eps)$ by
\begin{equation*}
	n_\eps=\sum_{\kappa=0}^\infty {n_{\eps,\kappa}}\mathbf{1}_{\Omega_\kappa},\quad c_\eps=\sum_{\kappa=0}^\infty {c_{\eps,\kappa}}\mathbf{1}_{\Omega_\kappa},\text{ and } \ \tau_\eps=\sum_{\kappa=0}^\infty \rho_{\eps,\kappa}\mathbf{1}_{\Omega_\kappa}.
\end{equation*}
We immediately observe that $\tau_\eps>0$ almost surely. Moreover, since the filtration $\mathbb{F}_T$ is right-continuous, it follows from \cite[P. 6 and 7]{Karatzas} that $\tau_\eps$ is a stopping time.
By the definition of $\tau_\eps$ and $\rho_{\eps, \kappa}$, we have
\begin{equation*}
	\begin{split}
		&\be \left[\sup_{s\in[0,T\wedge\tau_\eps]}\norm{n_\eps(s)}^4_{L^2}+\left(\int_0^{T\wedge\tau_{\eps}}\norm{n_\eps(s)}^2_{W^{1,2}}ds\right)^2\right]\\
		&=\be\sum_{\kappa=0}^\infty\mathbf{1}_{\Omega_\kappa}\left[\sup_{s\in[0,T\wedge\tau_\eps]}\norm{n_{\eps,\kappa}(s)}^4_{L^2}+\left(\int_0^{T\wedge\tau_{\eps}}\norm{n_{\eps,\kappa}(s)}^2_{W^{1,2}}ds\right)^2\right]\\
		&\leq \be\sum_{\kappa=0}^\infty\mathbf{1}_{\Omega_\kappa}\left(M+\norm{{u_0}}^2_{L^2}+\norm{{v_0}}^2_{B^\frac{3}{2}_{4,4}}\right)^2\leq 2M^{{2}}+2\be[\norm{{u_0}}_{L^2}^2+\norm{{v_0}}^2_{B^\frac{3}{2}_{4,4}}]^2<\infty,
	\end{split}
\end{equation*}
and thanks to \re{equa1.5},
\begin{equation*}
	\begin{split}
		&\be \left[\sup_{s\in[0,T\wedge\tau_\eps]}\norm{c_\eps(s)}^4_{B^\frac{3}{2}_{4,4}}+\left(\int_0^{T\wedge\tau_{\eps}}\norm{c_\eps(s)}^2_{W^{2,2}}ds\right)^2\right]\\
		&\leq \be\sum_{\kappa=0}^\infty\mathbf{1}_{\Omega_\kappa}\bk\left(M+\norm{{u_0}}^2_{L^2}+\norm{{v_0}}^2_{B^\frac{3}{2}_{4,4}}\right)^2\leq 2\bk M^{{2}}+2{\bk}\be[\norm{{u_0}}_{L^2}^2+\norm{{v_0}}^2_{B^\frac{3}{2}_{4,4}}]^2<\infty.
	\end{split}
\end{equation*}
Moreover, since $(n_{\eps,\kappa},c_{\eps,\kappa}, \rho_{\eps,\kappa})$ is a local solution to \re{equa3.1} with initial data $(u_{0,\kappa},v_{0,\kappa})$, we infer that for any $\varphi\in W^{1,2}(\bo)$, $\psi\in L^2(\bo)$ and $t\in [0,T]$,
\begin{equation*}
	\begin{split}
		&(n_\eps(t\wedge\tau_\eps),\varphi)=\sum_{\kappa=0}^\infty\mathbf{1}_{\Omega_\kappa}(n_{\eps,\kappa}(t\wedge\tau_{\eps,\kappa}),\varphi)\\
		&	=\sum_{\kappa=0}^\infty\mathbf{1}_{\Omega_\kappa}(u_{0,\kappa},\varphi)-\sum_{\kappa=0}^\infty\mathbf{1}_{\Omega_\kappa}\int_0^{t\wedge\rho_{\eps,\kappa}}( \nabla n_{\eps,\kappa}(s),\nabla\varphi)ds+\chi \sum_{\kappa=0}^\infty\mathbf{1}_{\Omega_\kappa}\int_0^{t\wedge\rho_{\eps,\kappa}}(F_{\eps}(n^+_{\eps,\kappa}(s))\nabla c_{\eps,\kappa}(s),\nabla \varphi )ds\\
		&\qquad+\sum_{\kappa=0}^\infty\mathbf{1}_{\Omega_\kappa}\int_0^{t\wedge\rho_{\eps,\kappa}}(\alpha n_{\eps,\kappa}(s)-\mu [H_{\eps}(n^+_{\eps,\kappa}(s))]^2,\varphi)ds+\sum_{\kappa=0}^\infty\mathbf{1}_{\Omega_\kappa}\int_0^{t\wedge\rho_{\eps,\kappa}}(g_1(n^+_{\eps,\kappa}(s)),\varphi)dW_1(s)\\
		&	=({u_0},\varphi)- \int_0^{t\wedge\tau_\eps}( \nabla n_\eps(s),\nabla\varphi)ds+\chi  \int_0^{t\wedge\tau_\eps}(F_\eps(n^+_\eps(s))\nabla c_\eps(s),\nabla \varphi )ds\\
		&\qquad+ \int_0^{t\wedge\tau_\eps}(\alpha n_\eps(s)-\mu [H_\eps(n^+_\eps(s))]^2,\varphi)ds+ \int_0^{t\wedge\tau_\eps}(g_1(n^+_\eps(s)),\varphi)dW_1(s),
	\end{split}
\end{equation*}
and
\begin{equation*}
	\begin{split}
		&(c_\eps(t\wedge\tau_\eps),\psi)=	\sum_{\kappa=0}^\infty\mathbf{1}_{\Omega_\kappa}(c_{\eps,\kappa}(t\wedge\tau_{\eps,\kappa}),\psi)\\
		&=\sum_{\kappa=0}^\infty\mathbf{1}_{\Omega_\kappa}(v_{0,\kappa},\psi)-\sum_{\kappa=0}^\infty\mathbf{1}_{\Omega_\kappa}\int_0^{t\wedge\rho_{\eps,\kappa}}(A_4 c_{\eps,\kappa}(s),\psi)ds\\
		&\qquad+\sum_{\kappa=0}^\infty\mathbf{1}_{\Omega_\kappa}\int_0^{t\wedge\rho_{\eps,\kappa}}(H_{\eps}(n^+_{\eps,\kappa}(s)),\psi)ds+\sum_{\kappa=0}^\infty\mathbf{1}_{\Omega_\kappa}\int_0^{t\wedge\rho_{\eps,\kappa}}(g_2(c_{\eps,\kappa}(s)),\psi)dW_2(s)\\
		&=(v_0,\psi)-\int_0^{t\wedge\tau_\eps}(A_4 c_\eps(s),\psi)ds +\int_0^{t\wedge\tau_\eps}(H_\eps(n^+_\eps(s)),\psi)ds+\int_0^{t\wedge\tau_\eps}(g_2(c_\eps(s)),\psi)dW_2(s).
	\end{split}
\end{equation*}
Hence, the triple $(n_\eps,c_\eps, \tau_\eps)$ constructed above is a local solution to the regularised stochastic variational system \re{equa3.1} in the sense of Definition \ref{defi3.1}.


$\bullet$ \textbf{Step 2: pathwise uniqueness of local solutions.} Assume that $(n^i_\eps,c_\eps^i,\tau_\eps^i)$, $i=1,2$ are two local solutions to the regularised system \re{equa3.1}. For $M>0$, we define the stopping time $\tau_{\eps, M}:=\tau^2_{\eps, M}\wedge\tau_\eps^1\wedge\tau_\eps^2$, where
\begin{equation}
	\tau_{\eps,M}^2:=\begin{cases}
		\inf\{ t\in [0,T]: \Sup_{s\in [0,t]}\norm{c^2_\eps(s\wedge\tau_\eps^2)}_{W^{1,4}}\geq M\},\\
		T,\text{ if the above set  is empty}.
	\end{cases}
\end{equation}
Since the process $(c_2(t\wedge\tau_\eps^2))_{t\in [0,T]}$ is continuous with values in $B_{4,4,\mathbf{n}}^\frac{3}{2}\hookrightarrow W^{1,4}(\bo)$ (this comes from the definition of local solution and \re{equa1.5}), it follows from \cite[Proposition 4.5]{Revuz} that $\tau_{\eps,M}^2$ is a stopping time with respect to the filtration $\mathbb{F}_T$. In addition, for any $t\in [0,T]$, 
\begin{equation}\label{equa2.24*}
	\mathbb{P}\left\{ \omega\in \Omega:  t\wedge\tau_{\eps, M}(\omega) \xrightarrow{M\uparrow \infty}  t\wedge\tau_\eps^1(\omega)\wedge\tau_\eps^2(\omega)\right\}=1.
\end{equation}
In fact, let us fix $t\in [0,T]$. Since $(\tau_{\eps,M}^2)_{M>0}$ is a family of increasing stopping times and for any $\omega\in \Omega$, we have $t\wedge\tau^2_{\eps,M}(\omega)\wedge\tau_{\eps}^1(\omega)\wedge\tau_\eps^2(\omega)\leq t\wedge\tau_{\eps}^1(\omega)\wedge\tau_\eps^2(\omega)$ for any $M>0$, the sequence $(t\wedge\tau_{\eps,M}(\omega))_{M>0}$ converges to $\Sup_{M>0}t\wedge\tau_{\eps,M}(\omega)$ for any $\omega\in\Omega$. This, and the fact that $\Sup_{M>0}t\wedge\tau_{\eps,M}(\omega)\leq  t\wedge\tau_\eps^1(\omega)\wedge\tau_\eps^2(\omega)$ for all $\omega\in\Omega$, implies that
\begin{equation*}
	\begin{split}
		&\mathbb{P}\left\{ \omega\in \Omega:  t\wedge\tau_{\eps, M}(\omega) \xrightarrow{M\uparrow \infty}  t\wedge\tau_\eps^1(\omega)\wedge\tau_\eps^2(\omega)\right\}\\
		&=\mathbb{P}\left\{ \omega\in \Omega:  \Sup_{M>0}t\wedge\tau_{\eps,M}(\omega)= t\wedge\tau_\eps^1(\omega)\wedge\tau_\eps^2(\omega)\right\}\\
		&=1-\mathbb{P}\left\{ \omega\in \Omega:  \Sup_{M>0}t\wedge\tau_{\eps,M}(\omega)< t\wedge\tau_\eps^1(\omega)\wedge\tau_\eps^2(\omega)\right\}\\
		&\geq 1-\lim_{M\longrightarrow \infty}\mathbb{P} \left\{ \omega\in \Omega:   t\wedge\tau_{\eps,M}(\omega)< t\wedge\tau_\eps^1(\omega)\wedge\tau_\eps^2(\omega)\right\}.
	\end{split}
\end{equation*}
 In the last line, we used the fact that for any $M>0$, $\left\{ \omega\in \Omega:  \Sup_{M>0}t\wedge\tau_{\eps,M}(\omega)< t\wedge\tau_\eps^1(\omega)\wedge\tau_\eps^2(\omega)\right\}\subset\left\{ \omega\in \Omega:  t\wedge\tau_{\eps,M}(\omega)< t\wedge\tau_\eps^1(\omega)\wedge\tau_\eps^2(\omega)\right\}$. By the definition of local solution (see Definition \ref{defi3.1}) and the embedding \re{equa1.5}, we have $\be\Sup_{s\in [0,t\wedge\tau_\eps^1\wedge\tau_\eps^2]}\norm{c_2(s)}_{W^{1,4}}^4<\infty$ and using the Markov inequality, we get
\begin{equation*}
	\begin{split}
		&\mathbb{P}\left\{ \omega\in \Omega:  t\wedge\tau_{\eps, M}(\omega) \xrightarrow{M\uparrow \infty}  t\wedge\tau_\eps^1(\omega)\wedge\tau_\eps^2(\omega)\right\}\\
		&\geq 1-\lim_{M\longrightarrow \infty}\mathbb{P} \left\{ \omega\in \Omega:   \sup_{s\in [0,t\wedge\tau_\eps^1\wedge\tau_\eps^2]}\norm{c{^2_\eps}(s)}_{W^{1,4}}{\geq M}\right\}\\
		&\geq  1-\lim_{M\longrightarrow \infty}\frac{\be\Sup_{s\in [0,t\wedge\tau_\eps^1\wedge\tau_\eps^2]}\norm{c_2(s)}_{W^{1,4}}^4}{M^4}\geq  1-\lim_{M\longrightarrow \infty}\frac{\be\Sup_{s\in [0,t\wedge \tau_\eps^2]}\norm{c_2(s)}_{W^{1,4}}^4}{M^4}= 1.
	\end{split}
\end{equation*}
Let us set $u_\eps=n_\eps^1-n_\eps^2$ and $v_\eps=c^1_\eps-c^2_\eps$ and fix arbitrary $M>0$. Then $(u_\eps,v_\eps,\tau_{\eps,M})$ satisfies for all $t\in [0,T]$,
\begin{equation*}
	\begin{split}
		&	(u_\eps(t\wedge\tau_{\eps,M}),\varphi)+\int_0^{t\wedge\tau_{\eps,M}}( \nabla u_\eps(s),\nabla\varphi)ds\\
		&=\chi \int_0^{t\wedge\tau_{\eps,M}}(F_\eps((n^1_\eps(s))^+)\nabla v_\eps(s),\nabla \varphi )ds+\chi \int_0^{t\wedge\tau_{\eps,M}}((F_\eps((n^1_\eps(s))^+)-F_\eps((n^2_\eps(s))^+))\nabla c^2_\eps(s),\nabla \varphi )ds\\
		&	\hspace{0.5cm}+\int_0^{t\wedge\tau_{\eps,M}}(\alpha u_\eps(s)-\mu [H_\eps((n^1_\eps(s))^+)]^2+\mu [H_\eps((n^2_\eps(s))^+)]^2,\varphi)ds\\
		&	\hspace{0.5cm}+\int_0^{t\wedge\tau_{\eps,M}}(g_1((n^1_\eps(s))^+)-g_1((n^2_\eps(s))^+),\varphi)dW_1(s),\quad\varphi\in W^{1,2}(\bo),\quad  \mathbb{P}\text{-a.s.}, \\
		&	(v_\eps(t\wedge\tau_{\eps,M}),\psi)+\int_0^{t\wedge\tau_{\eps,M}}(A{_2} v_\eps(s),\psi)ds-\int_0^{t\wedge\tau_{\eps,M}}(H_\eps((n^1_\eps(s))^+)-H_\eps((n^2_\eps(s))^+),\psi)ds\\
		&=\int_0^{t\wedge\tau_{\eps,M}}(g_2(c^1_\eps(s))-g_2(c^2_\eps(s)),\psi)dW_2(s),\quad\psi\in L^2(\bo),\quad\mathbb{P}\text{-a.s.}.
	\end{split}
\end{equation*}
Let us now fix arbitrary $t\in [0,T]$. Since $(n_\eps^i,c_\eps^i,\tau_\eps^i)$, $i=1,2$ satisfies \re{equa3.20}, we can apply the It\^o formula (from \cite[Theorem 1.2]{Pardoux}) to obtain for any $r\in [0,t\wedge\tau_{\eps,M}]$,
\begin{equation}\label{equa2.24}
	\begin{split}
		&\norm{u_\eps(r)}_{L^2}^2+2\int_0^{r}\norm{\nabla u_\eps(s)}^2_{L^2}ds\\
		&=2\chi \int_0^{r}(F_\eps((n^1_\eps(s))^+)\nabla v_\eps(s),\nabla u_\eps(s) )ds+\int_0^{r}\norm{g_1((n^1_\eps(s))^+)-g_1((n^2_\eps(s))^+)}^2_{\mathcal{L}_2^0}ds\\
		&\hspace{0.5cm}+2\int_0^{r}(g_1((n^1_\eps(s))^+)-g_1((n^2_\eps(s))^+),u_\eps(s))dW_1(s)\\
		&\hspace{0.5cm}+2\chi \int_0^{r}((F_\eps((n^1_\eps(s))^+)-F_\eps((n^2_\eps(s))^+))\nabla c^2_\eps(s),\nabla u_\eps(s))ds\\
		&	\hspace{0.5cm}+2\int_0^{r}(\alpha u_\eps(s)-\mu [H_\eps((n^1_\eps(s))^+)]^2+\mu [H_\eps((n^2_\eps(s))^+)]^2,u_\eps(s))ds.
	\end{split}
\end{equation}
It follows from \re{equa3.3} that
\begin{equation*}
	2\chi(F_\eps((n^1_\eps)^+)\nabla v_\eps,\nabla u_\eps )\leq \frac{2\chi}{\eps}\norm{\nabla v_\eps}_{L^2}\norm{\nabla u_\eps}_{L^2}\leq \frac{2\chi^2}{\eps^2}\norm{\nabla v_\eps}_{L^2}^2+\frac{1}{2}\norm{\nabla u_\eps}^2_{L^2}.
\end{equation*}
Let us point out that, by \re{equa2.3*},
\begin{equation*}
	\norm{g_1((n^1_\eps)^+)-g_1((n^2_\eps)^+)}^2_{\mathcal{L}_2^0}\leq L_1^2\norm{(n^1_\eps)^+-(n^2_\eps)^+}_{L^2}^{{2}}\leq L^2_1\norm{u_\eps}^2_{L^2},
\end{equation*}
and by \re{equa3.2}, the Lipschitz property of the function $\max(\cdot,0)$ and \re{equa2.3},
\begin{equation*}
	\begin{split}
		2\chi((F_\eps((n^1_\eps)^+)-F_\eps((n^2_\eps)^+))\nabla c^2_\eps,\nabla u_\eps)&\leq2\chi\norm{(n^1_\eps)^+-(n^2_\eps)^+}_{L^4}\norm{\nabla c^2_\eps}_{L^4}\norm{ \nabla u_\eps}_{L^2} \\
		&\leq4\chi^2\norm{n^1_\eps-n^2_\eps}^2_{L^4}\norm{\nabla c^2_\eps}^2_{L^4}+\frac{1}{4}\norm{ \nabla u_\eps}^2_{L^2} \\
		&\leq\bk\norm{u_\eps}_{L^2}\norm{u_\eps}_{W^{1,2}}\norm{\nabla c^2_\eps}^2_{L^4}+\frac{1}{4}\norm{ \nabla u_\eps}^2_{L^2} \\
		&\leq\bk\norm{u_\eps}^2_{L^2}(\norm{ c^2_\eps}^4_{W^{1,4}}+1)+\frac{1}{2}\norm{ \nabla u_\eps}^2_{L^2}.
	\end{split}
\end{equation*}
Thanks to {\re{equa3.2} and} \re{equa3.4}, 
it holds that
\begin{equation*}
	\begin{split}
		&(\alpha u_\eps-\mu [H_\eps((n^1_\eps)^+)]^2+\mu [H_\eps((n^2_\eps)^+)]^2,u_\eps)\\
		&=\alpha\norm{u_\eps}_{L^2}^2+\mu([H_\eps((n^2_\eps)^+)-H_\eps((n^1_\eps)^+)][H_\eps((n^1_\eps)^+)+H_\eps((n^2_\eps)^+)],u_\eps)\\
		&\leq {|}\alpha{|}\norm{u_\eps}_{L^2}^2+{\frac{8}{3}\mu}\bk_\eps\norm{u_\eps}^2_{L^2},
	\end{split}
\end{equation*}
and{, by the Burkholder--Davis--Gundy inequality,}
\begin{equation*}
	\begin{split}
		&2\be\sup_{r\in[0,t\wedge\tau_{\eps,M}]}\abs{\int_0^{r}(g_1((n^1_\eps(s))^+)-g_1((n^2_\eps(s))^+),u_\eps(s))dW_1(s)}\\
		&\leq \bk\be\left(\int_0^{t\wedge\tau_{\eps,M}}\norm{g_1((n^1_\eps(s))^+)-g_1((n^2_\eps(s))^+)}^2_{\mathcal{L}_2^0}\norm{u_\eps(s)}^2_{L^2}ds\right)^\frac{1}{2}\\
		&\leq \frac{1}{2}\be\sup_{s\in[0,t\wedge\tau_{\eps,M}]}\norm{u_\eps(s)}^2_{L^2}+\bk\be\int_0^{t\wedge\tau_{\eps,M}}\norm{u_\eps(s)}^2_{L^2}ds.
	\end{split}
\end{equation*}
Taking the supremum over $r\in [0,t\wedge\tau_{\eps, M}]$ in \re{equa2.24} and the mathematical expectation, combining with these calculations and the definition of stopping times, we derive that
\begin{equation}\label{equa3.23}
	\begin{split}
		&\frac{1}{2}\be\sup_{r\in[0,t\wedge\tau_{\eps,M}]}\norm{u_\eps(r)}_{L^2}^2+\be\int_0^{t\wedge\tau_{\eps,M}}\norm{\nabla u_\eps(s)}^2_{L^2}ds\\
		&\leq \frac{2\chi^2}{\eps^2}\be\int_0^{t\wedge\tau_{\eps,M}}\norm{\nabla v_\eps(s)}_{L^2}^2ds+\bk_\eps\be\int_0^{t\wedge\tau_{\eps,M}}\norm{u_\eps(s)}^2_{L^2}\left(\norm{ c^2_\eps(s)}^4_{W^{1,4}}+1\right)ds\\
		&\leq \frac{2\chi^2}{\eps^2}\be\int_0^{t\wedge\tau_{\eps,M}}\norm{\nabla v_\eps(s)}_{L^2}^2ds+\bk_{\eps,M}\be\int_0^{t\wedge\tau_{\eps,M}}\norm{u_\eps(s)}^2_{L^2}ds.
	\end{split}
\end{equation}
On the other hand, the It\^o formula (from \cite[Theorem 4.32]{DaPrato}) gives
\begin{equation}\label{equa2.26}
	\begin{split}
		&	\norm{v_\eps(r)}_{L^2}^2+2\int_0^{r}\norm{ v_\eps(s)}^2_{W^{1,2}}ds 	\\	&=2\int_0^{r}(H_\eps((n^1_\eps(s))^+)-H_\eps((n^2_\eps(s))^+),v_\eps(s))ds+\int_0^{r}\norm{g_2(c^1_\eps(s))-g_2(c^2_\eps(s))}^2_{\mathcal{L}_2^0}ds\\
		&\qquad+2\int_0^{r}(g_2(c^1_\eps(s))-g_2(c^2_\eps(s)),v_\eps(s))dW_2(s).
	\end{split}
\end{equation}
It follows from \re{equa3.2} that
\begin{equation*}
	2(H_\eps((n^1_\eps)^+)-H_\eps((n^2_\eps)^+),v_\eps)\leq {\frac{8}{3}}\norm{(n^1_\eps)^+-(n^2_\eps)^+}_{L^2}\norm{v_\eps}_{L^2}\leq \frac{4}{3}\norm{u_\eps}^2_{L^2}+\frac{4}{3}\norm{v_\eps}^2_{L^2}.
\end{equation*}
Note that
	$\norm{g_2(c^1_\eps)-g_2(c^2_\eps)}^2_{\mathcal{L}_2^0}\leq L^2_2\norm{v_\eps}_{L^2}^2,
$ 
and
\begin{equation*}
\be\sup_{r\in[0,t\wedge\tau_{\eps,M}]}\abs{\int_0^{r}(g_2(c^1_\eps(s))-g_2(c^2_\eps(s)),v_\eps(s))dW_2(s)}\leq \frac{1}{2}\be\sup_{s\in[0,t\wedge\tau_{\eps,M}]}\norm{v_\eps(s)}^2_{L^2}+\bk\be\int_0^{t\wedge\tau_{\eps,M}}\norm{v_\eps(s)}^2_{L^2}ds.
\end{equation*}
Using these {inequalities}, we take the supremum over $r\in [0,t\wedge\tau_{\eps,M}]$ and the mathematical expectation in \re{equa2.26} to derive that
\begin{equation*}
	\begin{split}
		&\frac{1}{2}\be\sup_{r\in[0,t\wedge\tau_{\eps,M}]}\norm{v_\eps(r)}_{L^2}^2+2\be\int_0^{t\wedge\tau_{\eps,M}}\norm{ v_\eps(s)}^2_{W^{1,2}}ds\\
		&\leq \bk\be\int_0^{t\wedge\tau_{\eps,M}}\norm{v_\eps(s)}^2_{L^2}ds+7\be\int_0^{t\wedge\tau_{\eps,M}}\norm{u_\eps(s)}^2_{L^2}ds.
	\end{split}
\end{equation*}
Multiplying this inequality by $\frac{\chi^2}{\eps^2}$ and adding the result with \re{equa3.23}, we get
\begin{equation*}
	\begin{split}
		&\frac{1}{2}\be\sup_{r\in[0,t\wedge\tau_{\eps,M}]}\norm{u_\eps(r)}_{L^2}^2+\frac{\chi^2}{2\eps^2}\be\sup_{r\in[0,t\wedge\tau_{\eps,M}]}\norm{v_\eps(r)}_{L^2}^2\\
		&\qquad+\be\int_0^{t\wedge\tau_{\eps,M}}\norm{\nabla u_\eps(s)}^2_{L^2}ds+\frac{2\chi^2}{\eps^2}\be\int_0^{t\wedge\tau_{\eps,M}}\norm{v_\eps(s)}^2_{L^2}ds\\
		&\leq \bk_{\eps,M}\int_0^{{t}}\be{\sup_{r\in[0,s\wedge\tau_{\eps,M}]}}\norm{u_\eps({r})}^2_{L^2}ds+\bk\int_0^{{t}}\be{\sup_{r\in[0,s\wedge\tau_{\eps,M}]}}\norm{v_\eps({r})}^2_{L^2}ds,
	\end{split}
\end{equation*}
and the Gronwall lemma yields
$
\be\Sup_{r\in[0,t\wedge\tau_{\eps,M}]}\norm{u_\eps(r)}_{L^2}^2+\be\Sup_{r\in[0,t\wedge\tau_{\eps,M}]}\norm{v_\eps(r)}_{L^2}^2=0 \quad \forall t\in [0,T],
$
and therefore,
\begin{equation*}
	\be \norm{u_\eps(t\wedge\tau_{\eps,M})}_{L^2}^2+\be \norm{v_\eps(t\wedge\tau_{\eps,M})}_{L^2}^2=0 \quad \forall t\in [0,T].
\end{equation*}
By the Fatou lemma and the $\mathbb{P}$-almost sure convergence \re{equa2.24*}, we derive that
\begin{equation*}
	\be\norm{u_\eps(t\wedge\tau_{\eps}^1\wedge\tau_\eps^2)}_{L^2}^2+\be\norm{v_\eps(t\wedge\tau_{\eps}^1\wedge\tau_\eps^2)}_{L^2}^2=0\quad \forall t\in [0,T],
\end{equation*}
which implies that $u_\eps(t\wedge\tau_{\eps}^1\wedge\tau_\eps^2)=0$ and $v_\eps(t\wedge\tau_{\eps}^1\wedge\tau_\eps^2)=0$ for any $t\in [0,T]$, almost surely. This proves the uniqueness result.


$\bullet$ \textbf{Step 3: Existence of a maximal local solution.} The construction will be done in an adaptation of \cite[Section 4]{Glatt-Holtz1} and \cite[Theorem 3]{Mikulevicuis}, see also \cite[Section 3.4]{Brze3}. Let $\ba_\eps$ be the set of all stopping times $\tau_\eps$ for which there exist processes $(n_\eps(t))_{t\in [0,T]}$ and $(c_\eps(t))_{t\in [0,T]}$ such that the triple $(n_\eps,c_\eps,\tau_\eps)$ is a local solution to the system \re{equa3.1}. By Step 1, the set $\ba_\eps$ is not empty. By \cite[Chapter 5, Section 18]{Doob} there exists a stopping time $\tau^{\text{max}}$ such that $\tau^{\text{max}}\geq \tau_\eps$ almost surely for all $\tau_\eps\in \ba_\eps$ and a sequence of increasing stopping times $(\tau_{\eps,\kappa})_{\kappa=1}^\infty\subset\ba_\eps$ satisfying $\tau_{\eps,\kappa}\nearrow \tau^{\text{max}}$ almost surely. Let $(n_{\eps,\kappa},c_{\eps,\kappa},\tau_{\eps,\kappa})$ be the respective local solution to the system \re{equa3.1} (by definition of $\ba_\eps$). By pathwise uniqueness of local solutions proved in Step 2, we derive that for all integers $\kappa$ and $m$, the following holds almost surely
\begin{equation*}
	(n_{\eps,\kappa}(t\wedge\tau_{\eps,\kappa}\wedge\tau_{\eps,m}),c_{\eps,\kappa}(t\wedge\tau_{\eps,\kappa}\wedge\tau_{\eps,m}))=(n_{\eps,m}(t\wedge\tau_{\eps,\kappa}\wedge\tau_{\eps,m}),c_{\eps,m}(t\wedge\tau_{\eps,\kappa}\wedge\tau_{\eps,m})), \ t\in [0,T].
\end{equation*}
Hence, by the amalgamation argument in \cite[Lemma III.6.B]{Elworthy}, there exists an $\mathbb{F}_T$-adapted $L^2(\bo)\times W^{1,4}(\bo)$-valued process $(n_{\eps}(t),c_\eps(t))_{t\in [0,T]}$ such that for every $t\in [0,T]$ and almost every $\omega\in \Omega$,
\begin{equation}\label{equa2.30}
	(n_\eps(\omega,t ),c_\eps(\omega,t ))=(n_{\eps,\kappa}(t\wedge\tau_{\eps,\kappa}(\omega),\omega),c_{\eps,\kappa}(t\wedge\tau_{\eps,\kappa}(\omega),\omega)) \quad\text{ for all } \quad \kappa\geq 2\chi.
\end{equation}
It is straightforward to check that $(n_\eps,c_\eps,\tau^{\text{max}}_\eps)$ satisfies items (i) and ({ii}) of Definition \ref{defi3.3} and is therefore a local maximal solution to \re{equa3.1}. Let us point out here that the maximal local solution constructed here blows up {at time $\tau^{\max}_\eps$ on the event $\{\tau^{\max}_\eps<T\}$}. In fact, we claim that,

\begin{claim}\label{Claim2.7}
	\begin{equation}\label{equa2.9}
		\sup_{0\leq s\leq \tau_\eps^\text{max}}\left(\norm{n_\eps(s)}^2_{L^2}+\norm{c_\eps(s)}^2_{B_{4,4}^\frac{3}{2}}\right)+\int_0^{\tau_\eps^\text{max}}\left( \norm{n_\eps(s)}^2_{W^{1,2}}+\norm{A_2c_\eps(s)}^2_{L^2} \right)ds=\infty,
	\end{equation}
	almost surely on the set	$\{ \omega\in \Omega: \tau^\text{max}_\eps(\omega)<{T} \}$.
\end{claim}

 \textbf{Proof of Claim \ref{Claim2.7}.} For any integer $\kappa\geq 2\chi$, let us define
\begin{equation*}
	\begin{split}
		\rho_{\eps, \kappa}&:=\inf\left\{ t\in [0,T]: \Sup_{s\in[0,t]}\norm{n_{\eps}(s)}^2_{L^2}+\Sup_{s\in[0,t]}\norm{c_{\eps}(s)}^2_{B^\frac{3}{2}_{4,4}}\right.\\
		&\left.\hspace{5.6cm}+\int_0^{t}\left(\norm{n_{\eps}(s)}^2_{W^{1,2}}+\norm{c_{\eps}(s)}^2_{W^{2,2}}\right)ds\geq  \kappa \right\}\wedge\tau_{\eps}^\text{max},
	\end{split}
\end{equation*}
with the convention that $\inf \emptyset=T$.
Note that $(\rho_\kappa)_{\kappa\geq 2\chi}$ is a family of stopping times due to \re{equa2.30} and the fact that for any $\kappa\geq 2\chi$ the path of $(n_{\eps,\kappa}(\cdot\wedge\tau_{\eps,\kappa}), c_{\eps,\kappa}(\cdot\wedge\tau_{\eps,\kappa}))$ is continuous in $L^2(\bo)\times B^\frac{3}{2}_{4,4,\mathbf{n}}(\bo)$. Let us point out that $\{\omega\in \Omega: \tau_\eps^\text{max}(\omega)<{T}\}=\Omega_1\cup\Omega_2$ where $\Omega_i$, $i=1,2$ are disjoint sets of $\Omega$ defined respectively by
\begin{equation*}
	\begin{split}
		&\Omega_1:=\left\{\omega\in \Omega: \tau_\eps^\text{max}(\omega)<{T}\right\}\cap \left\{\omega\in \Omega :\rho_{\eps, \kappa}(\omega)< \tau_\eps^\text{max}(\omega)\text{ for all integer } \kappa\geq 2\chi\right\},\\
		&\Omega_2:=\left\{\omega\in \Omega: \tau_\eps^\text{max}(\omega)<{T}\right\}\cap \left\{\omega\in \Omega :\rho_{\eps, \kappa}(\omega)= \tau_\eps^\text{max}(\omega)\text{ for some integer } \kappa\geq 2\chi\right\}.
	\end{split}
\end{equation*}
It is not too tedious to show that on the set $\Omega_1$, the blow-up property \re{equa2.9} holds. Indeed, let $\omega\in \Omega_1$. By definition of $\rho_{\eps, \kappa}$, we have
\begin{equation*}
	\begin{split}
		&\Sup_{s\in[0,\rho_{\eps, \kappa}(\omega)]}\norm{n_{\eps}(s)}^2_{L^2}+\Sup_{s\in[0,\rho_{\eps, \kappa}(\omega)]}\norm{c_{\eps}(s)}^2_{B^\frac{3}{2}_{4,4}}\\
		&\hspace{3cm}+\int_0^{\rho_{\eps, \kappa}(\omega)}\left(\norm{n_{\eps}(s)}^2_{W^{1,2}}+\norm{c_{\eps}(s)}^2_{W^{2,2}}\right)ds\geq  \kappa \quad\forall \kappa\geq 2\chi.
	\end{split}
\end{equation*}
Since $[0,\rho_{\eps, \kappa}(\omega)]\subset [0,\tau^\text{max}_{\eps}(\omega)]$, we infer from the above inequality that
\begin{equation*}
	\begin{split}
		&\Sup_{s\in[0,\tau^\text{max}_{\eps}(\omega)]}\norm{n_{\eps}(s)}^2_{L^2}+\Sup_{s\in[0,\tau^\text{max}_{\eps}(\omega)]}\norm{c_{\eps}(s)}^2_{B^\frac{3}{2}_{4,4}}\\
		&\hspace{3cm}+\int_0^{\tau^\text{max}_{\eps}(\omega)}\left(\norm{n_{\eps}(s)}^2_{W^{1,2}}+\norm{c_{\eps}(s)}^2_{W^{2,2}}\right)ds\geq  \kappa \quad\forall \kappa\geq 2\chi.
	\end{split}
\end{equation*}
Letting $\kappa$ go to $\infty$ in this inequality, we obtain \re{equa2.9} for $\omega\in \Omega_1$. Let us show that $\mathbb{P}\{\omega\in \Omega_2\}=0$. By contradiction, let us assume that $\mathbb{P}\{\omega\in \Omega_2\}>0$. Then there exists an integer $\kappa_0\geq 2\chi$ such that
$\Omega^{\kappa_0}_2:=\left\{\omega\in \Omega: \tau_\eps^\text{max}(\omega)<{T}\right\}\cap \left\{\omega\in \Omega :\rho_{\eps, \kappa_0}(\omega)= \tau_\eps^\text{max}(\omega)\right\}$ and $\mathbb{P}\{\omega\in \Omega^{\kappa_0}_2\}>0$ because of the fact that
$\Omega_2=\bigcup_{\kappa\geq 2\chi}^\infty\left\{\omega\in \Omega: \tau_\eps^\text{max}(\omega)<{T}\right\}\cap \left\{\omega\in \Omega :\rho_{\eps, \kappa}(\omega)= \tau_\eps^\text{max}(\omega)\right\}$. In particular, we have
\begin{equation*}
	\begin{split}
		&\Sup_{s\in[0,\tau^\text{max}_{\eps}(\omega)]}\norm{n_{\eps}(s)}^2_{L^2}+\Sup_{s\in[0,\tau^\text{max}_{\eps}(\omega)]}\norm{c_{\eps}(s)}^2_{B^\frac{3}{2}_{4,4}}\\
		&\hspace{3cm}+\int_0^{\tau^\text{max}_{\eps}(\omega)}\left(\norm{n_{\eps }(s)}^2_{W^{1,2}}+\norm{c_{\eps }(s)}^2_{W^{2,2}}\right)ds<  \kappa_0 \quad\omega\in \Omega_2^{\kappa_0},
	\end{split}
\end{equation*}
and hence
\begin{equation}\label{equa2.31}
	\norm{n_{\eps}(\tau^\text{max}_{\eps}(\omega))}^2_{L^2}+\norm{c_{\eps}(\tau^\text{max}_{\eps}(\omega))}^2_{B^\frac{3}{2}_{4,4}}\leq \kappa_0 \quad\omega\in \Omega_2^{\kappa_0}.
\end{equation}
Let us define now for any $\omega\in \Omega$,
\begin{equation*}
	(\tilde{n}^0_\eps(\omega),\tilde{c}^0_\eps(\omega)):=\begin{cases}
		(n_{\eps}(\tau^\text{max}_{\eps}(\omega),\omega),c_{\eps}(\tau^\text{max}_{\eps}(\omega),\omega))\quad\text{ if }\quad \omega\in \Omega^{\kappa_0}_2,\\
		0\qquad\text{ otherwise}.
	\end{cases}
\end{equation*}
From the definition it is immediate that $\left\{\omega\in \Omega: \tau_\eps^\text{max}(\omega)<{T}\right\}\in \mathcal{F}_{\tau_\eps^\text{max}}$ (for the definition of the $\sigma$-algebra defined from a stopping time, we refer the reader to \cite[Chap.I, Definition 2.12]{Karatzas}), and hence, by \cite[Chap. I, Lemma 2.15]{Karatzas}, $\left\{\omega\in \Omega :\rho_{\eps, \kappa_0}(\omega)= \tau_\eps^\text{max}(\omega)\right\}\in \mathcal{F}_{\tau_\eps^\text{max}}$. This implies that $\Omega^{\kappa_0}_2\in  \mathcal{F}_{\tau_\eps^\text{max}}$ and by \re{equa2.31}, 
$$\be \norm{\tilde{n}^0_\eps}^4_{L^2}+\be\norm{\tilde{c}^0_\eps}^4_{B^\frac{3}{2}_{4,4}}\leq 4\kappa^2_0
$$ 
showing that $(\tilde{n}^0_\eps ,\tilde{c}^0_\eps )$
 satisfies the assumption \re{equa2.1} on the initial condition with $\mathcal{F}_{\tau_\eps^\text{max}}$ in place of $\mathcal{F}_0$. We define now the filtered probability space
\begin{equation*}
	(\tilde{\Omega}, \tilde{\mathcal{F}}, \tilde{\mathbb{P}}, \{\tilde{\mathcal{F}}_t\}_{t\in [0,T]}):=(\Omega, \mathcal{F}, \mathbb{P}, \{\mathcal{F}_{\tau_\eps^\text{max}+t}\}_{t\in [0,T]}),
\end{equation*}
and the cylindrical Wiener processes
$
\tilde{W}_i(t):=W_i(\tau_\eps^\text{max}+t)-W_i(\tau_\eps^\text{max}), \ i=1,2.
$ By Step 1, there exists a local solution $(\tilde{n}_\eps,\tilde{c}_\eps,\tilde{\tau}_\eps)$ to \re{equa3.1} with initial data $(\tilde{n}^0_\eps ,\tilde{c}^0_\eps )$, relative to the cylindrical Wiener processes $\tilde{W}_i$, $i=1,2$, on the stochastic basis $(\tilde{\Omega}, \tilde{\mathcal{F}}, \tilde{\mathbb{P}}, \{\tilde{\mathcal{F}}_t\}_{t\in [0,T]})$ with a $\mathbb{P}$-a.s. strictly positive stopping time $\tilde{\tau}_\eps$. By defining
\begin{equation*}
	\bar{n}_\eps(\omega,t ):=\begin{cases}
		n_\eps(\omega,t )\quad\text{ if }\quad \omega\in \Omega,\quad 0\leq t< \tau_\eps^\text{max}(\omega), \\
		n_\eps(\tau_{\eps,\kappa_0}(\omega),\omega)\quad\text{ if }\quad \omega\notin \Omega_2^{\kappa_0},\quad  \tau_\eps^\text{max}(\omega)\leq t\leq \tau_\eps^\text{max}(\omega)+\tilde{\tau}_\eps(\omega),\\
		\tilde{n}_\eps(t-{\tau^{\max}_\eps}(\omega),\omega)\quad\text{ if }\quad \omega\in \Omega^{\kappa_0}_2,\quad  \tau_\eps^\text{max}(\omega)\leq t\leq \tau_\eps^\text{max}(\omega)+\tilde{\tau}_\eps(\omega),
	\end{cases}
\end{equation*}
and
\begin{equation*}
	\bar{c}_\eps(\omega,t ):=\begin{cases}
		c_\eps(\omega,t )\quad\text{ if }\quad \omega\in \Omega,\quad 0\leq t< \tau_\eps^\text{max}(\omega), \\
		c_\eps(\tau_{\eps,\kappa_0}(\omega),\omega)\quad\text{ if }\quad \omega\notin \Omega_2^{\kappa_0},\quad  \tau_\eps^\text{max}(\omega)\leq t\leq \tau_\eps^\text{max}(\omega)+\tilde{\tau}_\eps(\omega),\\
		\tilde{c}_\eps(t-{\tau^{\max}_\eps}(\omega),\omega)\quad\text{ if }\quad \omega\in \Omega^{\kappa_0}_2,\quad  \tau_\eps^\text{max}(\omega)\leq t\leq \tau_\eps^\text{max}(\omega)+\tilde{\tau}_\eps(\omega),
	\end{cases}
\end{equation*}
it is straightforward {to check that $(\bar{n}_\eps,\bar{c}_\eps,\tau^{\max}_\eps+\tilde\tau_\eps\mathbf 1_{\Omega_2^{\kappa_0}})$ is a local solution to the system \re{equa3.1} with initial data $(u_0,v_0)$, and by the definition of maximal local solution (see Definition \ref{defi3.3}), we deduce that $\tau^{\max}_\eps+\tilde\tau_\eps\mathbf 1_{\Omega_2^{\kappa_0}}\leq \tau_\eps^\text{max}$ almost surely. This is a contradiction, since $\tilde{\tau}_\eps$ is a strictly positive stopping time and $\mathbb{P}\{\omega\in \Omega^{\kappa_0}_2\}>0$.} Thus $\mathbb{P}\{\omega\in \Omega_2\}=0$ and the blow-up property \re{equa2.9} is then proven.

Before moving to the next step, let us prove that 
\begin{claim}\label{lem3.6}
There exists a constant $\bk_\eps>0$, depending on $\eps$, such that for every $t\in [0,T]$,
	\begin{equation*}
		\be \sup_{0\leq s\leq t\wedge\tau_\eps^\text{max}}\left(\norm{n_\eps(s)}^4_{L^2}+\norm{c_\eps(s)}^4_{B_{4,4}^\frac{3}{2}}\right)+\be\left(\int_0^{t\wedge\tau_\eps^\text{max}}\left( \norm{n_\eps(s)}^2_{W^{1,2}}+\norm{c_\eps(s)}^2_{W^{2,2}} \right)ds\right)^2\leq \bk_\eps.
	\end{equation*}
\end{claim}
 \textbf{Proof of Claim \ref{lem3.6}.}
	Let $(\tau_{\eps,\kappa})_{\kappa\geq 2\chi}$ be a sequence of stopping times satisfying $\tau_{\eps,\kappa}\nearrow \tau_\eps^\text{max}$ almost surely. Let $\kappa\geq 2\chi$ be arbitrary fixed. Recalling that $(n_\eps,c_\eps)$ is a solution to the system \re{equa3.1}, let us fix arbitrary $t\in [0,T]$. Due to \re{equa3.20}, we can apply the It\^o formula coming from \cite[Theorem 4.32]{DaPrato} to the process $r\mapsto\norm{c_{\eps}(r)}_{L^2}^2$ to obtain for any $r\in [0,t\wedge\tau_{\eps,\kappa}]$,
	\begin{equation*}
		\begin{split}
			\norm{ c_{\eps}(r)}_{L^2}^2+2\int_0^r\norm{ c_{\eps}(s)}_{W^{1,2}}^2ds
			&=\norm{v_0}^2_{L^2}+2\int_0^r(H_{\eps}(n^+_{\eps}(s)),c_{\eps}(s))ds +\int_0^r\norm{g_2(c_{\eps}(s))}^2_{\mathcal{L}_2^0}ds\\
			&\qquad+2\int_0^r(g_2(c_{\eps}(s)),c_{\eps}(s))dW_2(s).
		\end{split}
	\end{equation*}
	Here, we use \re{equa3.4}, the H\"older inequality and the Young inequality to derive that
	\begin{equation*}
		2(H_{\eps}(n^+_{\eps}),c_{\eps})\leq 2\bk_\eps\norm{c_\eps}_{L^1}\leq2\bk_\eps\abs{\bo}^\frac{1}{2}\norm{c_\eps}_{L^2}\leq \norm{c_\eps}^2_{W^{1,2}} +\bk_\eps.
	\end{equation*}
	Next, recalling \re{equa2.4}, we get
	\begin{equation*}
		\norm{g_2(c_{\eps})}^2_{\mathcal{L}_2^0}\leq 2L^2_2\norm{c_\eps}^2_{L^2}+2\norm{g_2(0)}^2_{\mathcal{L}_2^0} .
	\end{equation*}
	From these inequalities, we derive that
	\begin{equation*}
		\begin{split}
			\norm{ c_{\eps}(r)}_{L^2}^2+\int_0^{r}\norm{ c_{\eps}(s)}_{W^{1,2}}^2ds
			&\leq \norm{v_0}^2_{L^2}+\bk_\eps T+2L^2_2\int_0^{r}\norm{c_\eps(s)}^2_{L^2}ds+2T\norm{g_2(0)}^2_{\mathcal{L}_2^0} \\
			&\qquad+2\int_0^{r}(g_2(c_{\eps}(s)),c_{\eps}(s))dW_2(s).
		\end{split}
	\end{equation*}
	Taking the supremum over $r\in [0,t\wedge\tau_{\eps,\kappa}]${, squaring and taking the expectation}, we get
	\begin{equation*}
		\begin{split}
			&\be\sup_{r\in[0,t\wedge\tau_{\eps,\kappa}]}	\norm{ c_{\eps}(r)}_{L^2}^4+\be\left(\int_0^{t\wedge\tau_{\eps,\kappa}}\norm{ c_{\eps}(s)}_{W^{1,2}}^2ds\right)^2\\
			&\leq \bk\be\norm{v_0}^4_{L^2}+\bk_\eps T+\bk T\norm{g_2(0)}^2_{\mathcal{L}_2^0}+\bk\be\left(\int_0^{t\wedge\tau_{\eps,\kappa}}\norm{c_\eps(s)}^2_{L^2}ds\right)^2 \\
			&\qquad+\bk \be \sup_{r\in[0,t\wedge\tau_{\eps,\kappa}]}\abs{\int_0^{r}(g_2(c_{\eps}(s)),c_{\eps}(s))dW_2(s)}^2.
		\end{split}
	\end{equation*}
	Thanks to \re{equa2.4}, we get
	\begin{equation*}
		\begin{split}
			&\bk\be\sup_{r\in[0,t\wedge\tau_{\eps,\kappa}]}\abs{\int_0^r(g_2(c_{\eps}(s)),c_{\eps}(s))dW_2(s)}^2\\
			&\leq \bk\be \int_0^{t\wedge\tau_{\eps,\kappa}}\norm{g_2(c_{\eps}(s))}^2_{\mathcal{L}^0_2}\norm{c_{\eps}(s)}_{L^2}^2ds\\
			&\leq \bk\be\int_0^{t\wedge\tau_{\eps,\kappa}}\norm{g_2(c_{\eps}(s))-g_2(0)}^2_{\mathcal{L}^0_2}\norm{c_{\eps}(s)}_{L^2}^2ds+ \bk\norm{g_2(0)}^2_{\mathcal{L}^0_2}\be \int_0^{t\wedge\tau_{\eps,\kappa}}\norm{c_{\eps}(s)}_{L^2}^2ds \\
			&\leq \bk\be\int_0^{t\wedge\tau_{\eps,\kappa}}\norm{c_{\eps}(s)}_{L^2}^4ds+ \bk\norm{g_2(0)}^2_{\mathcal{L}^0_2}\be \int_0^{t\wedge\tau_{\eps,\kappa}}\norm{c_{\eps}(s)}_{L^2}^4ds+ \bk T\norm{g_2(0)}^2_{\mathcal{L}^0_2}.
		\end{split}
	\end{equation*}
	The H\"older inequality yields
	\begin{equation*}
		\be\left(\int_0^{t\wedge\tau_{\eps,\kappa}}\norm{c_\eps(s)}^2_{L^2}ds\right)^2 \leq T\be\int_0^{t\wedge\tau_{\eps,\kappa}}\norm{c_\eps(s)}^4_{L^2}ds.
	\end{equation*}
	Combining these three inequalities, we infer that
	\begin{equation*}
		\begin{split}
			&	\be\sup_{r\in[0,t\wedge\tau_{\eps,\kappa}]}	\norm{ c_{\eps}(r)}_{L^2}^4+\be\left(\int_0^{t\wedge\tau_{\eps,\kappa}}\norm{ c_{\eps}(s)}_{W^{1,2}}^2ds\right)^2\\
			&\leq \bk\be\norm{v_0}^4_{L^2}+\bk_\eps T+\bk T\norm{g_2(0)}^2_{\mathcal{L}_2^0}+\bk(\norm{g_2(0)}^2_{\mathcal{L}^0_2}+T+1)\be \int_0^{t\wedge\tau_{\eps,\kappa}}\norm{c_{\eps}(s)}_{L^2}^4ds \qquad \forall t\in [0,T],
		\end{split}
	\end{equation*}
	and the Gronwall Lemma yields
	\begin{equation}\label{equa2.32}
		\begin{split}
			&	\be\sup_{r\in[0,t\wedge\tau_{\eps,\kappa}]}	\norm{ c_{\eps}(r)}_{L^2}^4+\be\left(\int_0^{t\wedge\tau_{\eps,\kappa}}\norm{ c_{\eps}(s)}_{W^{1,2}}^2ds\right)^2\\
			&\leq \left(\bk\be\norm{v_0}^4_{L^2}+\bk_\eps T+\bk T\norm{g_2(0)}^2_{\mathcal{L}_2^0}\right)e^{\bk(\norm{g_2(0)}^2_{\mathcal{L}^0_2}+T+1)T}\quad \forall  t\in [0,T].
		\end{split}
	\end{equation}
	Next, the generalised version of the It\^o formula from \cite[Proposition 3.3]{Scarpa} to the process $r\mapsto \norm{\nabla c_{\eps}(r)}_{L^2}^2$ implies that for any $r\in [0,t\wedge\tau_{\eps,\kappa}]$,
	\begin{equation*}
		\begin{split}
			&\norm{ \nabla c_{\eps}(r)}_{L^2}^2+2\int_0^{r}\left(\norm{\Delta c_{\eps}(s)}^2_{L^2}+\norm{\nabla c_{\eps}(s)}_{L^2}^2\right)ds\\
			&=\norm{\nabla v_0}^2_{L^2}+2\int_0^{r}(H_{\eps}(n^+_{\eps}(s)),-\Delta c_{\eps}(s))ds+\int_0^{r}\norm{\nabla g_2(c_{\eps}(s))}^2_{\mathcal{L}_2^0}ds+2\int_0^{r}(\nabla g_2(c_{\eps}(s)),\nabla c_{\eps}(s))dW_2(s).
		\end{split}
	\end{equation*}
	Note that
	\begin{equation*}
		2(H_{\eps}(n^+_{\eps}),-\Delta c_{\eps})\leq 2\bk_\eps\norm{\Delta c_\eps}_{L^1}\leq2\bk_\eps\abs{\bo}^\frac{1}{2}\norm{\Delta c_\eps}_{L^2}\leq \norm{\Delta c_\eps}^2_{L^{2}} +\bk_\eps,
	\end{equation*}
	and by \re{equa2.4},
	\begin{equation*}
		\norm{\nabla g_2(c_{\eps})}^2_{\mathcal{L}_2^0}\leq 2L^2_2\norm{c_\eps}^2_{W^{1,2}}+2\norm{g_2(0)}^2_{\mathcal{L}_2^1} .
	\end{equation*}
	From these inequalities, we get
	\begin{equation*}
		\begin{split}
			&\norm{ \nabla c_{\eps}(r)}_{L^2}^2+\int_0^{r}\left(\norm{\Delta c_{\eps}(s)}^2_{L^2}+2\norm{\nabla c_{\eps}(s)}_{L^2}^2\right)ds\\
			&{\leq}\norm{\nabla v_0}^2_{L^2}+2L^2_2\int_0^{r}\norm{c_\eps(s)}^2_{W^{1,2}}ds+(\bk_\eps+2\norm{g_2(0)}^2_{\mathcal{L}_2^1} )T+2\int_0^{r}(\nabla g_2(c_{\eps}(s)),\nabla c_{\eps}(s))dW_2(s).
		\end{split}
	\end{equation*}
	Taking the supremum over $r\in [0,t\wedge\tau_{\eps,\kappa}]$, {squaring}, the mathematical expectation and using \re{equa2.32}, we obtain
	\begin{equation*}
		\begin{split}
			&\be \sup_{r\in[0,t\wedge\tau_{\eps,\kappa}]}\norm{ \nabla c_{\eps}(r)}_{L^2}^4+\be\left(\int_0^{t\wedge\tau_{\eps,\kappa}}\left(\norm{\Delta c_{\eps}(s)}^2_{L^2}+2\norm{\nabla c_{\eps}(s)}_{L^2}^2\right)ds\right)^2\\
			&\leq \bk \be\norm{\nabla v_0}^4_{L^2}+\bk+\bk\be\sup_{r\in[0,t\wedge\tau_{\eps,\kappa}]} \abs{\int_0^{r}(\nabla g_2(c_{\eps}(s)),\nabla c_{\eps}(s))dW_2(s)}^2.
		\end{split}
	\end{equation*}
	Taking into account \re{equa2.4}, the H\"older inequality and \re{equa2.32}, we get
	\begin{equation*}
		\begin{split}
			&\bk\be\sup_{r\in[0,t\wedge\tau_{\eps,\kappa}]}\abs{\int_0^r(\nabla g_2(c_{\eps}(s)),\nabla c_{\eps}(s))dW_2(s)}^2\\
			&\leq \bk\be\int_0^{t\wedge\tau_{\eps,\kappa}}\norm{g_2(c_{\eps}(s))-g_2(0)}^2_{\mathcal{L}^1_2}\norm{\nabla c_{\eps}(s)}_{L^2}^2ds+ \bk\norm{g_2(0)}^2_{\mathcal{L}^1_2}\be \int_0^{t\wedge\tau_{\eps,\kappa}}\norm{\nabla c_{\eps}(s)}_{L^2}^2ds \\
			&\leq \bk\be\int_0^{t\wedge\tau_{\eps,\kappa}}\norm{c_{\eps}(s)}_{W^{1,2}}^2\norm{\nabla c_{\eps}(s)}_{L^2}^2ds+ \bk\norm{g_2(0)}^2_{\mathcal{L}^0_2}\left[\be\left(\int_0^{t\wedge\tau_{\eps,\kappa}}\norm{c_{\eps}(s)}_{W^{1,2}}^2ds\right)^2\right]^\frac{1}{2}\\
			&\leq \bk\be\sup_{r\in[0,t\wedge\tau_{\eps,\kappa}]}\norm{\nabla c_\eps(s)}^2_{L^2}\int_0^{t\wedge\tau_{\eps,\kappa}}\norm{c_{\eps}(s)}_{W^{1,2}}^2ds+\bk\\
			&\leq \frac{1}{2}\be\sup_{r\in[0,t\wedge\tau_{\eps,\kappa}]}\norm{\nabla c_\eps(s)}^4_{L^2}+\bk \be\left(\int_0^{t\wedge\tau_{\eps,\kappa}}\norm{c_{\eps}(s)}_{W^{1,2}}^2ds\right)^2+\bk\leq\frac{1}{2}\be\sup_{r\in[0,t\wedge\tau_{\eps,\kappa}]}\norm{\nabla c_\eps(s)}^4_{L^2}+ \bk.
		\end{split}
	\end{equation*}
	Combining these two inequalities, we get
	\begin{equation}\label{equa2.33}
		\begin{split}
			&\frac{1}{2}\be \sup_{r\in[0,t\wedge\tau_{\eps,\kappa}]}\norm{ \nabla c_{\eps}(r)}_{L^2}^4+\be\left(\int_0^{t\wedge\tau_{\eps,\kappa}}\left(\norm{\Delta c_{\eps}(s)}^2_{L^2}+2\norm{\nabla c_{\eps}(s)}_{L^2}^2\right)ds\right)^2\\
			&\leq \bk \be\norm{\nabla v_0}^4_{L^2}+\bk, \quad\forall t\in [0,T].
		\end{split}
	\end{equation}
	In the other side, since the second equation of \re{equa3.5*} holds in $L^2(\bo)$, it can be written in the mild form as follows for any $r\in [0,t\wedge\tau_{\eps}]$,
	\begin{equation*}
		c_{\eps}(r)=e^{r\Delta}v_0-\int_0^{r}e^{(r-s)\Delta} c_{\eps}(s)ds+\int_0^{r}e^{(r-s)\Delta} H_\eps(n^+_{\eps}(s))ds+\int_0^{r}e^{(r-s)\Delta} g_2(c_{\eps}(s))dW_2(s).
	\end{equation*}
	This, and \re{equa1.8*} imply that
	\begin{equation}\label{equa3.19*}
		\begin{split}
			\be\sup_{r\in[0,t\wedge\tau_\kappa]}	\norm{c_{\eps}(r)}^2_{B^\frac{3}{2}_{4,4}}
			&\leq 4\be\sup_{r\in[0,t\wedge\tau_\kappa]}\norm{e^{r\Delta}v_0}_{B^\frac{3}{2}_{4,4}}^2+4\be\sup_{r\in[0,t\wedge\tau_\kappa]}\left(\int_0^{r}\norm{e^{(r-s)\Delta} c_{\eps}(s)}_{B^\frac{3}{2}_{4,4}}ds\right)^2\\
			&\qquad+4\be\sup_{{r}\in[0,t\wedge\tau_\kappa]}\left(\int_0^{r}\norm{e^{(r-s)\Delta} H_\eps(n^+_{\eps}(s))}_{B^\frac{3}{2}_{4,4}}ds\right)^2\\
			&\qquad+4\be\sup_{r\in[0,t\wedge\tau_\kappa]}\norm{\int_0^{r}e^{(r-s)\Delta} g_2(c_{\eps}(s))dW_2(s)}_{B^\frac{3}{2}_{4,4}}^2\\
			&:=I_1+I_2+I_3+I_4.
		\end{split}
	\end{equation}
	Let $r\in[0,t\wedge\tau_\kappa]$ be fixed. 
{Since $\{e^{s\Delta}\}_{s\in[0,T]}$ is uniformly bounded on both $L^4(\bo)$ and $D(A_4)$, real interpolation $B^{3/2}_{4,4,\mathbf n}(\bo)=(L^4(\bo),D(A_4))_{\frac34,4}$ yields its uniform boundedness on $B^{3/2}_{4,4,\mathbf n}(\bo)$. Thus,
\begin{equation*}
	\norm{e^{r\Delta}v_0}_{B^{3/2}_{4,4}}\le\bk_T\norm{v_0}_{B^{3/2}_{4,4}},
\end{equation*}	}
and consequently,
	\begin{equation}\label{equa2.37}
		I_1
		\leq \bk_T \be \norm{v_0}_{B^\frac{3}{2}_{4,4}}^2.
	\end{equation}
{	Let us also fix $\theta\in ({\frac32},{2})$; then $\frac{\theta}{2}\in ({\frac34},{1})$ and by \re{equa1.7}, combined with \eqref{equ2.5} and \eqref{equ1.3}, we have 
	\begin{equation*}
D(A_4^\frac{\theta}{2})\hookrightarrow H^{\theta,4}(\bo)\hookrightarrow B^\theta_{4,4}(\bo)=W^{\theta,4}(\bo)\hookrightarrow W^{\frac{3}{2},4}(\bo)=B^\frac{3}{2}_{4,4}(\bo),
	\end{equation*}}
	 and therefore \re{equa2.4*} and the H\"older inequality along with the Young inequality yield,
	\begin{equation*}
		\begin{split}
			\int_0^{r}\norm{e^{(r-s)\Delta} c_{\eps}(s)}_{B^\frac{3}{2}_{4,4}}ds&\leq \bk \int_0^{r}\norm{A^\frac{\theta}{2}_4e^{(r-s)\Delta} c_{\eps}(s)}_{L^4}ds\\
			&\leq \bk\int_0^r (r-s)^{-\frac{\theta}{2}}\norm{c_{\eps}(s)}_{L^4}ds	\\
			&\leq {\frac{2\bk}{2-\theta}\,r^{1-\frac\theta2}\sup_{s\in [0,r]}\norm{c_{\eps}(s)}_{L^4}}.
		\end{split}
	\end{equation*}
	From this, and the embedding $W^{1,2}(\bo)\hookrightarrow L^4(\bo)$ combined with \re{equa2.32} and \re{equa2.33}, we deduce that
	\begin{equation}	\label{equa2.38}
		I_2\leq {\bk T^{2-\theta}\,\be\sup_{r\in[0,t\wedge\tau_\kappa]}\norm{c_{\eps}(r)}^2_{W^{1,2}}\leq\bk}.
	\end{equation}
	In addition, \re{equa3.4} helps us to have
	\begin{equation*}
		\begin{split}
			\int_0^r\norm{e^{(r-s)\Delta} H_\eps((n_{\eps}(s))_+)}_{B^\frac{3}{2}_{4,4}}ds
			&\leq \bk\int_0^r\norm{A_4^\frac{\theta}{2} e^{(r-s)\Delta} H_\eps(n^+_{\eps}(s))}_{L^4}ds\\
			&\leq \bk\int_0^r (r-s)^{-\frac{\theta}{2}}e^{-\lambda (r-s)}\norm{H_\eps(n^+_{\eps}(s))}_{L^4}ds	\\
			&\leq \bk_\eps\abs{\bo}^\frac{1}{4} \int_0^r (r-s)^{-\frac{\theta}{2}}ds	\leq \bk_\eps r^{1-\frac{\theta}{2}},
		\end{split}
	\end{equation*}
	and hence,
	\begin{equation}\label{equa2.39}
		I_3\leq \bk_\eps T^{2-\theta}.
	\end{equation}
	To control the stochastic term, {let us point out that 
		\begin{equation}\label{A32}
\begin{split}
I_4&:=4\be\sup_{r\in[0,t\wedge\tau_\kappa]}\norm{\int_0^{r}e^{(r-s)\Delta} g_2(c_{\eps}(s))dW_2(s)}_{B^\frac{3}{2}_{4,4}}^2\\
&\leq 4\be\sup_{r\in[0,t]}\norm{\int_0^{r}e^{(r-s)\Delta} 1_{[0,\tau_\kappa]}(s)g_2(c_{\eps}(s))dW_2(s)}_{B^\frac{3}{2}_{4,4}}^2\\
&\leq 2\be\sup_{r\in[0,t]}\norm{\int_0^{r}e^{(r-s)\Delta} 1_{[0,\tau_\kappa]}(s)g_2(c_{\eps}(s))dW_2(s)}_{B^\frac{3}{2}_{4,4}}^4+2.
\end{split}		
		\end{equation}
		 let define for any $z\in [0,t]$,
		\begin{equation*}
			M_\eps(z):=\int_0^ze^{(z-s)\Delta} 1_{[0,\tau_{\kappa}]}(s) g_2(c_{\eps}(s))dW_2(s).
	\end{equation*}
Then, $M_\eps$ satisfying the following stochastic differential equation in $L^4(\bo)$,
\begin{equation}\label{A33}
M_\eps(z)+\int_0^zA_4M_\eps(s)ds=\int_0^zM_\eps(s)ds+\int_0^z1_{[0,\tau_{\kappa}]}(s) g_2(c_{\eps}(s))dW_2(s), \quad z\in [0,t],
\end{equation}
where we have used the fact that $-\Delta=A_4-I$. Setting
$
X_0:=L^4(\bo)$ and $
X_1:=D(A_4),
$
we have
$
X_{\frac12}:=[X_0,X_1]_{\frac12}=D(A_4^{\frac12})
$
and
$
X_{\frac34,4}:=(X_0,X_1)_{\frac34,4}
=B^{\frac32}_{4,4,\mathbf n}(\bo).
$
Moreover, since
$
D(A_4^{\frac12})=H^{1,4}(\bo)=W^{1,4}(\bo),
$
by \eqref{equa2.4**}, the Sobolev embedding
$
W^{1,2}(\bo)\hookrightarrow L^4(\bo),
$
and the Gagliardo--Nirenberg interpolation inequality \eqref{equa2.3}, together with \eqref{equa2.32} and \eqref{equa2.33}, we obtain
\begin{equation*}
\begin{split}
&\be\int_0^t\norm{1_{[0,\tau_{\kappa}]}(s) g_2(c_{\eps}(s))}^4_{\gamma(L^2;D(A_4^\frac{1}{2}))}\\
&\leq L_3^4\be\int_0^{t\wedge\tau_{\kappa}}\norm{c_{\eps}(s)}_{L^4}^4ds+L_3^4\be\int_0^{t\wedge\tau_{\kappa}}\norm{\nabla c_{\eps}(s)}_{L^4}^4ds\\
&\leq \bk\be\int_0^{t\wedge\tau_{\kappa}}\norm{c_{\eps}(s)}_{W^{1,2}}^4ds+\bk\be\int_0^{t\wedge\tau_{\kappa}}\norm{\nabla c_{\eps}(s)}_{L^2}^2\norm{ c_{\eps}(s)}_{W^{2,2}}^2ds\\
&\leq \bk\be\int_0^{t\wedge\tau_{\kappa}}\norm{c_{\eps}(s)}_{W^{1,2}}^4ds+\bk\be\int_0^{t\wedge\tau_{\kappa}}\norm{\nabla c_{\eps}(s)}_{L^2}^2\norm{ \Delta c_{\eps}(s)}_{L^2}^2ds+\bk\be\int_0^{t\wedge\tau_{\kappa}}\norm{\nabla c_{\eps}(s)}_{L^2}^2\norm{ c_{\eps}(s)}_{W^{1,2}}^2ds\\
&\leq \bk t\be\sup_{0\leq s\leq t\wedge\tau_\kappa}\norm{c_{\eps}(s)}_{W^{1,2}}^4+\bk\be\sup_{0\leq s\leq t\wedge\tau_\kappa}\norm{\nabla c_{\eps}(s)}_{L^2}^2\int_0^{t\wedge\tau_{\kappa}}\norm{ \Delta c_{\eps}(s)}_{L^2}^2ds\\
&\leq \bk t\be\sup_{0\leq s\leq t\wedge\tau_\kappa}\norm{c_{\eps}(s)}_{W^{1,2}}^4+\bk\be\sup_{0\leq s\leq t\wedge\tau_\kappa}\norm{\nabla c_{\eps}(s)}_{L^2}^4+\bk \be\left(\int_0^{t\wedge\tau_{\kappa}}\norm{ \Delta c_{\eps}(s)}_{L^2}^2ds\right)^2<\infty.
\end{split}
\end{equation*}
Due to this, it follows from the stochastic maximal $L^4$-regularity \cite[Theorem 4.5 (ii)]{Neerven} (with $\tilde{L}_F=1$) that  there exist a constant $\bk_\eps>0$ such that 
\begin{equation*}
	\be \sup_{r\in[0,t]}\norm{M_\eps(r)}^4_{B^\frac{3}{2}_{4,4}}\leq \bk_\eps,
\end{equation*}}
which, combined with \eqref{A32}, yields
\begin{equation*}
	I_4\leq 2\bk_\eps+2.
\end{equation*} 
	Inserting this, with the estimates \re{equa2.37}, \re{equa2.38} and \re{equa2.39} in \re{equa3.19*}, we get
	\begin{equation}\label{equa2.36}
		\be\sup_{r\in[0,t\wedge\tau_\kappa]}	\norm{c_{\eps}(r)}^2_{B^\frac{3}{2}_{4,4}}\leq \bk_\eps.
	\end{equation}
	As far as the estimate of $n_\eps$ is concerned, let us apply the It\^o formula in version of \cite[Theorem 1.2]{Pardoux} to the process $r\mapsto\norm{n_{\eps}(r)}^2_{L^2}$ (which is possible due to \re{equa3.20} and the first equation of \re{equa3.5*}) to obtain for any $r\in [0,t\wedge\tau_{\eps,\kappa}]$,
	\begin{equation*}
		\begin{split}
			\norm{n_{\eps}(r)}_{L^2}^2+2\int_0^{r}\norm{\nabla n_{\eps}(s)}^2_{L^2}ds
			&=\norm{u_0}_{L^2}^2+ 2\chi \int_0^{r}(F_{\eps}(n^+_{\eps}(s))\nabla c_{\eps}(s),\nabla n_{\eps}(s) )ds\\
			&\hspace{0.5cm}+2\alpha\int_0^{r}\norm{n_{\eps}(s)}^2_{L^2}ds-2\mu\int_0^{r} ([H_{\eps}(n^+_{\eps}(s))]^2,n_{\eps}(s))ds\\
			&\hspace{0.5cm}+\int_0^{r}\norm{g_1(n^+_{\eps}(s))}^2_{\mathcal{L}_2^0}ds+2\int_0^{r}(g_1(n^+_{\eps}(s)),n_{\eps}(s))dW_1(s).
		\end{split}
	\end{equation*}
	Thanks to \re{equa3.3} and the H\"older inequality, we get
	\begin{equation*}
		2\chi(F_{\eps}(n^+_{\eps})\nabla c_{\eps},\nabla n_{\eps} )\leq \frac{2\chi}{\eps}\norm{\nabla c_{\eps}}_{L^2}\norm{\nabla n_{\eps} }_{L^2}\leq \norm{\nabla n_{\eps} }_{L^2}^2+\frac{\chi^2}{\eps^2}\norm{\nabla c_{\eps}}^2_{L^2}.
	\end{equation*}
	Noting that $n_\eps=n_\eps^+-n_\eps^-$, where $n_\eps^-:=\max(-n_\eps,0)$, it holds that
	\begin{equation*}
		-2\mu([H_{\eps}(n^+_{\eps})]^2,n_{\eps})=-2\mu([H_{\eps}(n^+_{\eps})]^2,n^+_{\eps})+2\mu([H_{\eps}(n^+_{\eps})]^2,n^-_{\eps})=-2\mu([H_{\eps}(n^+_{\eps})]^2,n_{\eps}^+)\leq 0,
	\end{equation*}
	where we have used the fact that $([H_{\eps}(n^+_{\eps})]^2,n^-_{\eps})=0$. By \re{equa2.3*}, we have
	\begin{equation*}
		\norm{g_1(n^+_{\eps})}^2_{\mathcal{L}_2^0}\leq L_0^2.
	\end{equation*}
	From these inequalities, we deduce that
	\begin{equation*}
		\begin{split}
			\norm{n_{\eps}(r)}_{L^2}^2+\int_0^{r}\norm{\nabla n_{\eps}(s)}^2_{L^2}ds
			&\leq \norm{u_0}_{L^2}^2+L^2_0+ \frac{\chi^2}{\eps^2} \int_0^{r}\norm{\nabla c_{\eps}(s)}_{L^2}^2ds+2\alpha\int_0^{r}\norm{n_{\eps}(s)}^2_{L^2}ds\\
			&\hspace{0.5cm} +2\int_0^{r}(g_1(n^+_{\eps}(s)),n_{\eps}(s))dW_1(s).
		\end{split}
	\end{equation*}
	From this, we deduce that
	\begin{equation*}
		\begin{split}
			&\be \sup_{r\in[0,t\wedge\tau_{\eps,\kappa}]}	\norm{n_{\eps}(r)}_{L^2}^4+\be\left(\int_0^{t\wedge\tau_{\eps,\kappa}}\norm{\nabla n_{\eps}(s)}^2_{L^2}ds\right)^2\\
			&\leq \bk\norm{u_0}_{L^2}^4+\bk+ \frac{\bk\chi^4}{\eps^4} \be\left(\int_0^{t\wedge\tau_{\eps,\kappa}}\norm{\nabla c_{\eps}(s)}_{L^2}^2ds\right)^2\\
			&\qquad+\bk\be\left( \int_0^{t\wedge\tau_{\eps,\kappa}}\norm{n_{\eps}(s)}^2_{L^2}ds\right)^2 +\bk\be  \sup_{r\in[0,t\wedge\tau_{\eps,\kappa}]}\abs{\int_0^{r}(g_1(n^+_{\eps}(s)),n_{\eps}(s))dW_1(s)}^2.
		\end{split}
	\end{equation*}
	By \re{equa2.33}, we get
	\begin{equation*}
		\be\left(\int_0^{t\wedge\tau_{\eps,\kappa}}\norm{\nabla c_{\eps}(s)}_{L^2}^2ds\right)^2\leq \be\left(\int_0^{T\wedge\tau_{\eps,\kappa}}\norm{\nabla c_{\eps}(s)}_{L^2}^2ds\right)^2\leq T^{{2}}\be \sup_{r\in[0,T\wedge\tau_{\eps,\kappa}]}\norm{ \nabla c_{\eps}(r)}_{L^2}^{{4}}\leq \bk.
	\end{equation*}
	The H\"older inequality gives
	\begin{equation*}
		\be\left( \int_0^{t\wedge\tau_{\eps,\kappa}}\norm{n_{\eps}(s)}^2_{L^2}ds\right)^2\leq T \be \int_0^{t\wedge\tau_{\eps,\kappa}}\norm{n_{\eps}(s)}^4_{L^2}ds,
	\end{equation*}
	and using \re{equa2.3*}, we get
	\begin{equation*}
		\begin{split}
			\be  \sup_{r\in[0,t\wedge\tau_{\eps,\kappa}]}\abs{\int_0^{r}(g_1(n^+_{\eps}(s)),n_{\eps}(s))dW_1(s)}^2&\leq \be \int_0^{t\wedge\tau_{\eps,\kappa}}\norm{g_1(n^+_{\eps}(s))}^2_{\mathcal{L}^0_2}\norm{n_{\eps}(s)}^2_{L^2}ds\\
			&\leq L^2_0 \be \int_0^{t\wedge\tau_{\eps,\kappa}}\norm{n_{\eps}(s)}^4_{L^2}ds+L^2_0 T.
		\end{split}
	\end{equation*}
	Combining these inequalities, we derive that
	\begin{equation*}
		\begin{split}
			&\be \sup_{r\in[0,t\wedge\tau_{\eps,\kappa}]}	\norm{n_{\eps}(r)}_{L^2}^4+\be\left(\int_0^{t\wedge\tau_{\eps,\kappa}}\norm{\nabla n_{\eps}(s)}^2_{L^2}ds\right)^2\\
			&\leq \bk\be\norm{u_0}_{L^2}^4+\bk+L^2_0 T+(\bk+T)\be \int_0^{t\wedge\tau_{\eps,\kappa}}\norm{n_{\eps}(s)}^4_{L^2}ds \quad\forall t\in [0,T],
		\end{split}
	\end{equation*}
	and the Gronwall lemma yields
	\begin{equation*}
		\be \sup_{r\in[0,t\wedge\tau_{\eps,\kappa}]}	\norm{n_{\eps}(r)}_{L^2}^4+\be\left(\int_0^{t\wedge\tau_{\eps,\kappa}}\norm{\nabla n_{\eps}(s)}^2_{L^2}ds\right)^2
		\leq \left(\bk\be \norm{u_0}_{L^2}^4+\bk+L^2_0 T\right)e^{(\bk+T)T}\quad\forall t\in [0,T].
	\end{equation*}
	Combining this, with \re{equa2.32}, \re{equa2.33} and \re{equa2.36}, we get
	\begin{equation*}
		\be 	\sup_{0\leq s\leq t\wedge \tau_{\eps,\kappa}}\left(\norm{n_\eps(s)}^4_{L^2}+\norm{c_\eps(s)}^4_{B_{4,4}^\frac{3}{2}}\right)+\be\left(\int_0^{t\wedge \tau_{\eps,\kappa}}\left( \norm{n_\eps(s)}^2_{W^{1,2}}+\norm{c_\eps(s)}^2_{W^{2,2}} \right)ds\right)^4\leq \bk_\eps\quad\forall t\in [0,T],
	\end{equation*}
	and the claim follows by applying the monotone convergence theorem, using that $\tau_{\eps,\kappa}\nearrow \tau_\eps^\text{max}$ almost surely. This completes the proof of Claim \ref{lem3.6}.


$\bullet$ \textbf{Step 4: Global existence.} Let $(n_\eps,c_\eps, \tau_\eps^\text{max})$ be the maximal local solution to \re{equa3.1} obtained in Step 3. Our goal is to prove that $\tau_\eps^\text{max}=T$ almost surely. By contradiction, let us assume that $\mathbb{P}\{\omega\in\Omega_3\}>0$, where $\Omega_3:=\{\omega\in\Omega: \tau_\eps^\text{max}(\omega)<T\}$.  {On $\Omega_3$ one has $T\wedge\tau^{\max}_\eps=\tau^{\max}_\eps$, and Claim~\ref{lem3.6} shows that the quantity in \eqref{equa2.9} (with $t=T$) is finite almost surely, contradicting \eqref{equa2.9} on $\Omega_3$.} 
Hence, the maximal local solution obtained in Step 3 is a global solution to the system \re{equa3.1}, and \re{equa3.7*} follows from Claim \ref{lem3.6}, completing the proof of Proposition \ref{propo3.5}. 


\subsection*{Acknowledgments}
We acknowledge financial support provided by the Austrian Science Fund (FWF). In particular, Boris Jidjou Moghomye  and Erika Hausenblas   were supported by the Austrian Science Fund, project no.: PAT 3076923, Grant DOI: 10.55776/PAT3076923. 

\end{document}